\documentclass[a4paper,12pt,twoside]{article}
\usepackage[latin1]{inputenc}
\usepackage[T1]{fontenc}
\usepackage{amssymb}
\usepackage{dochier}
\usepackage{lpretex}
\usepackage{title}
\usepackage{a4}
\usepackage{amscd}
\usepackage[all]{xy}
\EndOfDump

\def\DHpartformat#1{(#1) }
\def\DHrefpart#1{(\DHRefpart{#1})}
\LBsetstyleof{partno}{arabic}

\DocVersion[Changed to LaTeX]{1}
\DocVersion[Added deformation theory, rearranged; Sep 09 1999]{2}
\DocVersion[Two sections 1 combined; Sep 16 1999]{3}
\DocVersion[Minor tweaks; Sep 17 1999]{4}
\DocVersion[Added a new section and minor fixes; Sep 23 1999]{5}
\DocVersion[Many small changes (NSB), One new section (TE); Wed Oct  6  1999]{6}
\DocVersion[Many small changes (NSB), Stacks (TE); Wed Oct  11  1999]{7}

\let\define\def

\def\A {{\mathbb A}}  
  \def\F {{\mathbb F}}

\def\GG {{\mathbb G}}   
  
\def\N {{\mathbb N}}  \def\P {{\mathbb P}} 
\def\Q {{\mathbb Q}} \def\R {{\mathbb R}}
 
 \def\W {{\mathbb W}} 
\def\Z {{\mathbb Z}} 

\define \n {\mathbb N}
\define \z {\mathbb Z}
\define \q {\mathbb Q}
\define \PP {\mathbb P}

\def\sA {{\Cal A}} \def\sB {{\Cal B}} \def\sC {{\Cal C}}
\def\sD {{\Cal D}} \def\sE {{\Cal E}} \def\sF {{\Cal F}}
\def\sG {{\Cal G}} \def\sH {{\Cal H}} \def\sI {{\Cal I}}
 \def\sK {{\Cal K}} \def\sL {{\Cal L}}
\def\sM {{\Cal M}} \def\sN {{\Cal N}} \def\sO {{\Cal O}}
\def\sQ {{\Cal Q}} \def\sR {{\Cal R}}
 \def\sT {{\Cal T}} \def\sU {{\Cal U}}
\def\sV {{\Cal V}} \def\sW {{\Cal W}} \def\sX {{\Cal X}}
\def\sY {{\Cal Y}}

\def\sZ {{\Cal Z}}
\define \cN {\Cal N}
\define \cf {\Cal F}
\define \cg {\Cal G}
\define \cE {\Cal E}
\define \ce {\Cal E}
\define \cc {\Cal C}
\define \cV {\Cal V}
\define \cA {\Cal A}
\define \cK {\Cal K}
\define \cO {\Cal O}
\define \cF {\Cal F}
\define \cn {\Cal N}
\define \cI {\Cal I}
\define \sP {\Cal P}
\def\a {\alpha} \def\b {\beta}   
   
\def\s {\sigma}

\define \x {\xi}
\define \y {\eta}
\define \G {\Gamma}
\define \r {\rho}
\define \w {\omega}

\def \tZ {\widetilde Z}
\def\tX {\widetilde X} 
\def \tY {\widetilde Y}

\def \tD {\widetilde D}

\def \tV {\widetilde V}
\def \trho {\tilde {\rho}}

\def \tp {\widetilde{\mathbb P}}
\define \tH {\widetilde H}
\define \tG {\widetilde{\Gamma}}
\define \tW {\widetilde W}
\define \tF {\widetilde F}
\define \tm {\tilde m}
\define \St {\widetilde S}
\define \Xt {\widetilde X}
\define \tS {\widetilde S}
\define \tpsi {\tilde \psi}
\define \tL {\widetilde L}
\define \tE {\widetilde E}
\define \tl {\tilde l}
\define \tA {\widetilde A}
\define \tom {\tilde\omega}
\define \tT {\widetilde T}
\define \tB {\widetilde B}
\define \tf {\tilde f}
\define \tsA {\widetilde{\sA}}

\define \tM {\widetilde M}
\define \tphi {\widetilde{\phi}}
\define \trho {\widetilde{\rho}}
\define \tR {\widetilde R}
\define \tp {\tilde p}
\define \tq {\tilde q}
\define \tc {\tilde c}
\define \tsF {\widetilde {\sF}}
\define \tx {\tilde x}
\define \tg {\tilde g}
\define \tw {\tilde w}

\def\pd {\partial}

\def \Dx1 {\frac{\pd}{{\pd} x_1}}
\def \Dy1 {\frac{\pd}{{\pd} y_1}}
\def \Dz1 {\frac{\pd}{{\pd} z_1}}
\def \Dx2 {\frac{\pd}{{\pd} x_2}}
\def \Dy2 {\frac{\pd}{{\pd} y_2}}
\def \Dz2 {\frac{\pd}{{\pd} z_2}}

\def\q {\quad} 

\def\mapdiagr#1{\Big\searrow\rlap{$\raise 5pt\vbox{{\hbox{$\mkern -15mu\scriptstyle#1$}}}$}}   

\def\mapdiagl#1{\llap{$\raise 5pt\vbox{{\hbox{$\scriptstyle#1\mkern
-15mu$}}}$}\Big\swarrow}              

\def\Mapdiagr#1{\nearrow\rlap{$\lower 5pt\vbox{{\hbox{$\mkern
-15mu\scriptstyle#1$}}}$}} 

\def\Mapdiagl#1{\llap{$\lower 5pt\vbox{{\hbox{$\scriptstyle#1\mkern
-15mu$}}}$}\searrow} 

\def\Mapswr#1{\swarrow\rlap{$\lower 5pt\vbox{{\hbox{$\mkern
-15mu\scriptstyle#1$}}}$}}              

\def\Mapnwl#1{\nwarrow\rlap{$\lower 5pt\vbox{{\hbox{$\mkern
-15mu\scriptstyle#1$}}}$}}

\def \inj {\hookrightarrow}

\define \Rhook {\hookrightarrow}

\def \half {\raise1pt\hbox{$\scriptstyle
        \frac{1}{2}\displaystyle$}}

\def \x{{\sl X}\llap{$\mkern -2mu {\scriptstyle -}$}}
\def \Hom {\operatorname{Hom}}
\def \Spec {\operatorname{Spec}}

\def \Symm {\operatorname{Sym}}

\def \rank {\operatorname{rank}}
\def \Bl {\operatorname{Bl}}
\def \Pic {\operatorname{Pic}}
\def \Sing {\operatorname{Sing}}
\define \Kod {\operatorname{Kod}}
\define \dimension {\operatorname{dim}}
\define \codim {\operatorname{codim}}
\define \contr {\operatorname{contr}}
\define \rk {\operatorname{rank}}
\define \im {\operatorname{im}}
\define \Mor {\operatorname{Mor}}
\define \Cl {\operatorname{Cl}}
\define \Hilb {\operatorname{Hilb}}
\define \degree {\operatorname{deg}}
\define \mult {\operatorname{mult}}
\define \Aut {\operatorname{Aut}}
\define \NS {\operatorname{NS}}
\define \Gal {\operatorname{Gal}}
\define \ch {\operatorname{char}}
\define \Jac {\operatorname{Jac}}
\define \Km {\operatorname{Km}}
\define \Sec {\operatorname{Sec}}
\define \Stab {\operatorname{Stab}}
\define \Br {\operatorname{Br}}
\define \inv {\operatorname{inv}}
\define \tr {\operatorname{tr}}
\define \Frob {\operatorname{Frob}}
\define \Symn {\operatorname{Sym}^n}
\define \Ev {\sE^\vee}
\define \ordp {\operatorname{ord}_p}
\define \Supp {\operatorname{Supp}}
\define \Ann {\operatorname{Ann}}
\define \disc {\operatorname{disc}}
\define \Lie {\operatorname{Lie}}
\define \embdim {\operatorname{embdim}}
\def \Picf {\operatorname{\underline{Pic}}}
\def \sEnd {\operatorname{\sE nd}}
\def \Num {\operatorname{Num}}
\def \diag {\operatorname{diag}}
\def \Der{\operatorname{Der}}

\def\Isom{\operatorname{Isom}}
\def\Aff{\operatorname{Aff}}
\def\gr{\operatorname{gr}}
\def\Diff{\operatorname{\Cal Diff}}
\def\lHom{\operatorname{\underline{Hom}}}
\def\lth{\operatorname{length}}
\def\Tr{\operatorname{Tr}}
\def\Div{\operatorname{Div}}

\def\ad{\operatorname{ad}}

\def\hR{\hat R}
\def\exc{\operatorname{exc}}
\def\Est{\sE}

\def\Kst{\sK}
\def\tKst{\widetilde{\Kst}}
\def\Def{\operatorname{Def}}

\def\map#1#2#3{\ensuremath{{#1}\co{#2}\to{#3}}}

\def\W{\cansymb{W}}

\newcommand\contract{\mathop{\mbox{\rule{5pt}{.6pt}\rule{.6pt}{6pt}}\,}}

\def\hod#1#2#3#4{\ensuremath{\if#30 H^{#2}({#1},{\cal O}_{#1}) \else 
 H^{#2}(#1,\Omega^{#3}\if\relax{#4}\relax_{#1}\else _{#1/#4}\fi)\fi}}

\begin{document}
\title{Moduli and periods of simply connected Enriques surfaces}

\author{T. Ekedahl${}^\dagger$}
\address{Department of Mathematics\\
 Stockholm University\\
 S-106 91  Stockholm\\
Sweden}

\author{J.M.E. Hyland}
\address{D.P.M.M.S.\\
 Centre for Mathematical Sciences\\
 Wilberforce Rd.\\
 Cambridge CB2 1SB\\
 U.K.}
\email{martin@pmms.cam.ac.uk}
\author{N. I. Shepherd-Barron}
\address{D.P.M.M.S.\\
 Centre for Mathematical Sciences\\
 Wilberforce Rd.\\
 Cambridge CB2 1SB\\
 U.K.}
\email{nisb@pmms.cam.ac.uk}

\maketitle

It is well known that the presence of automorphisms and vector fields
on a projective variety complicates the nature
of the corresponding moduli spaces. The aim of this paper is to prove some general
results that quantify this and to use them to describe various moduli spaces,
both local and global, of algebraic surfaces. 
Our approach is na\"ive: suppose given a \emph{Keel--Mori} stack,
that is, an algebraic stack $\sF$
over a noetherian base $B$ such that 
the representable ([LMB00], 4.1) 
diagonal morphism $\sF\to\sF\times_B\sF$ is finite, 
so that, by the fundamental theorem of Keel and Mori [KM97],
there is a geometric quotient $F$, which is the coarse moduli space.
Then we try to describe $F$ in elementary quasi-projective terms.
We can sometimes take a slightly more refined point of view,
which is to construct a Deligne--Mumford stack
$\sF_{DM}$ that lies between $\sF$ and $F$, so that, in particular,
$F$ is also the geometric quotient of $\sF_{DM}$, and then describe
$\sF_{DM}$ so as to make transparent the passage to $F$. We cannot, however, 
prove what ought obviously to be true, namely, that any Keel--Mori stack
has a universal Deligne--Mumford quotient, although we can
do this for stacks of Enriques surfaces. (It is easy to find examples where
there is a universal Deligne--Mumford quotient whose construction
is not local, unlike the construction of the geometric quotient.)
We are interested mainly in Enriques surfaces (for whose moduli
there do exist universal Deligne--Mumford quotients); recall that
in the classification of surfaces by Mumford and Bombieri
over an algebraically closed field $k$ of
arbitrary characteristic $p$
these form one of the simplest classes beyond the ruled surfaces.
They are defined by the properties that $2K \sim 0$ and $b_2=10$ 
and do sometimes, if $p=2$, possess vector fields. 

If $Y$ is a $G$-Enriques surface (that is, if
$G=\Pic^\tau(Y)$), then $G$ has order two, so that there is a 
$G^\vee$-torsor $X\to Y$ such that $X$ is reduced and irreducible and
the dualizing sheaf $\omega_X$ is trivial.
If $p\ne 2$, then $G^\vee\cong\Z/2$, $X$ is a smooth K3 surface,
$Y$ is the quotient of $X$ by a fixed--point--free involution and,
in consequence of the Rudakov--Shafarevich theorem on the non-existence
of vector fields on a K3 surface, $Y$ has no vector fields. If,
however, $p=2$ (this we regard as the interesting case), 
then $G$ is isomorphic to one of $\Z/2$, $\mu_2$, or $\alpha_2$, and $Y$ is
the quotient of a K3 surface by a fixed-point-free involution if and only if
$G=\mu_2$.  In the other cases $G$ is unipotent (i.e., not linearly reductive)
and the Cartier dual $G^\vee$ is not reduced, so that
the canonical $G^\vee$-torsor $X\to Y$ is purely inseparable. In this case $Y$
is simply connected; we shall also say that $Y$ is unipotent.

\begin{definition}
An \Definition{RDP--K3 surface} is a surface with rational double points (RDPs)
whose minimal resolution is K3. An Enriques surface $Y$ is a
\Definition{K3--Enriques surface} if its canonical double cover is RDP--K3.
\end{definition}

The chief global result is, in over-simplified terms, that there is a period
morphism for simply connected K3--Enriques surfaces that describes
the moduli space as essentially an open piece of a
$\P^1$-bundle over the period space. (To be a little more precise, if the
canonical double cover $X$ has
only RDPs, then its minimal resolution $\tX$ is a
supersingular K3 surface, for which Ogus, Rudakov and Shafarevich have
constructed a period map and proved it to be an isomorphism. We can construct a
period map for such Enriques surfaces, by taking the periods of $\tX$.)
At the level of geometric points the fact that the period map is
a fibration by curves of genus zero arises from the fact that $X$ 
has free tangent sheaf and every vector
field on it is $2$-closed; the genus zero curve is the set of lines of vector fields.
(This idea is that of Moret-Bailly's construction of complete $\P^1$'s of abelian surfaces
by taking quotients of products of supersingular elliptic curves, although our
construction does not give complete rational curves of Enriques surfaces; they
degenerate when the line of vector fields on $X$ specializes to vanish at a 
singular point of $X$.)
However, to make this precise, first at the level of stacks
and then at the level of coarse moduli, and to show that this fibration
really is a $\P^1$-bundle requires more. We prove some general results 
on stacks in order to deal easily with these issues
(these results first appeared in a Mittag--Leffler preprint, [EHS07]).
In particular, we can describe the difference between the number
of local moduli and the number of global moduli of a variety
in terms of the extent to which its automorphism group scheme $G$
fails to be reduced.

\begin{theorem} (= \ref{2.8})
Suppose that $f:\sX\to S$ is proper and flat and is minimally
versal at $s\in S$. Then the set of points $t\in S$ such that
$\sX_s$ is geometrically isomorphic to $\sX_s$ forms a smooth
subscheme of $S$ of dimension equal to $\dim\Lie(G)-\dim G$.
\noproof
\end{theorem}

Recently Liedtke has proved the following four basic results about moduli of Enriques
surfaces, especially in characteristic $2$. His key idea is to allow 
the surfaces $X$ to have RDPs and to put a \emph{Cossec--Verra polarization}
on them. 
\begin{theorem}
The stack $\Est$ of appropriately polarized
Enriques surfaces is
algebraic over $\Spec\Z$. All of its geometric components
are $10$-dimensional.
\noproof
\end{theorem}
\begin{theorem}
If $p\ne 2$, then $\Est\Tensor\F_p$ is absolutely irreducible.
\noproof
\end{theorem}
The remaining results concern the case where $p=2$. 
\begin{theorem}
An RDP-Enriques surface with a Cossec--Verra polarization 
can be lifted to $W(k)[{\sqrt{2}}]$.
\noproof
\end{theorem} 
\begin{theorem}
$\Est\Tensor\F_2$ has just two components, $\Est_{uni}$,
the locus where $\Pic^\tau$ is unipotent, 
and $\Est_{inf}$, where $\Pic^\tau$ is infinitesimal. Both are absolutely 
irreducible of dimension $10$. $\Est_{uni}$ is the locus of
simply connected surfaces and $\Est_{inf}$ is the closure of the locus of
$\mu_2$--surfaces.  They intersect in the locus $\Est_\alpha$
of $\alpha_2$--surfaces.
\noproof
\end{theorem}

As mentioned, these are due to Liedtke. Our contribution is as follows.

\begin{theorem} (= \ref{0.6})
Every $\a_2$-Enriques surface has a miniversal deformation
space that is a regular scheme of the form 
$\Sp W(k)[[x_1,\ldots,x_{12}]]/(FG-2)$.
\noproof
\end{theorem}

\begin{theorem} (= \ref{0.7})
Every Enriques surface has a lifting to a characteristic
zero DVR whose absolute ramification index divides $2^9N$, where
$N$ is one of the numbers $192,\ 128, \ 60,\ 56,\ 40,\ 9.$
\noproof
\end{theorem}

There is a period map defined on the open substack $\Est_{K3,uni}$
of $\Est_{uni}$ corresponding to K3--Enriques surfaces. 
(The complement of this open
substack has codimension at least $3$ everywhere.
This can be proved by counting constants
for surfaces $Y$ whose canonical double cover $X$ is normal 
but not RDP-$K3$; when $X$ is not normal
it is possible to make an estimate by finding a certain
negative definite configuration $C$ of $(-2)$-curves
on $Y$ and showing that 
$H^1(Y,T_Y(\log C))$ has high codimension
in $H^1(Y,T_Y)$. However, the details involve
a messy and unenlightening consideration
of different cases, so are omitted.)

The image of the period map lies in a quotient $\sM_N/\frak S_{12}$, where $\sM_N$ is the
period space for appropriately marked K3 surfaces [Og79].

For a morphism $X\to S$ of schemes or algebraic stacks 
in characteristic $p$, we denote by $X\to X^{(n)}$ the
$n$'th Frobenius relative to $S$.

\begin{theorem}\label{F}(= \ref{0.9}, \ref{8.12})
\part[i] There is a canonical Deligne--Mumford quotient 
of the stack $\Est_{K3,uni}$, denoted
$(\Est_{K3,uni})_{DM}$.

\part[ii] There is an open piece $M_N^0$, the complement of
an explicit divisor $D$ in a period space $M_N$, and a period morphism
$$\psi:(\Est_{K3,uni})_{DM}\to[M_N^0/\frak S_{12}]^{(1)},$$ 
where $\frak S_{12}$ is the
symmetric group on $12$ letters, that identifies $(\Est_{K3,uni})_{DM}$
with an open piece of a Zariski $\P^1$-bundle over 
the quotient Deligne--Mumford stack $[M_N^0/\frak S_{12}]^{(1)}$. 

\part[iii] In particular, the geometric quotient
$(\Est_{K3,uni})_{geom}$ is the quotient of an open piece of a 
$\P^1$-bundle over
$M_N^{0(1)}$ by an action of $\frak S_{12}$. 
\noproof
\end{theorem}

It follows from this that, since $\frak S_{12}$ acts generically
freely on $M_N$, the geometric quotient of $\Est_{uni}$ is birationally ruled
over $M_N^{(1)}/\frak S_{12}$. 

\begin{theorem} (= \ref{8.13}, \ref{8.14})
There is a closed substack
$\Est_\alpha$ of $\Est_{uni}$ whose geometric points correspond to Enriques surfaces
with $\Pic^\tau\cong\alpha_2$ that is an irreducible divisor. The intersection
$\Est_{K3,\alpha}=\Est_{K3,uni}\cap\Est_\alpha$ has a canonical Deligne--Mumford quotient
$(\Est_{K3,\alpha})_{DM}$ that is generically a section of the period map
$\psi$ of Theorem \ref{F}.
\noproof
\end{theorem}
\begin{notation}
Recall [CD89] that an $\alpha_2$-surface $Y$ has a global vector field $D$, unique
modulo scalars.  We say that $Y$ is \Definition{additive} (resp.,
\Definition{multiplicative}) if $D^2=0$ (resp., $D^2\ne 0$).
\end{notation}

\begin{theorem}(= \ref{7.7}, \ref{8.16})
For every multiplicative $\alpha_2$--surface the canonical
double cover is RDP-K3. The period map is defined
for them and the period map is an isomorphism on the appropriate stack.
\noproof
\end{theorem}

The paper is arranged as follows. In \S \ref{stacktology} there are various
constructions and theorems of a general nature concerning algebraic stacks
and particularly quotients. \S \ref{local-to-global} contains theorems relating 
miniversal deformation spaces to global moduli, particularly when there
are non-integrable vector fields (i.e., when automorphism groupschemes
are not reduced). 
\S \ref{singular foliations} contains
straightforward extensions of various well known results
about foliations on smooth varieties to the singular case, mostly
in order to deal with the singular $K3$ surfaces that arise.

Then we specialize to consider Enriques surfaces.
\S \ref{Deforming Enriques surfaces} describes the local moduli space
of an Enriques surface in a way that refines some of
Liedtke's results. \S \ref{blank} describes how to mark
and polarize Enriques surfaces, and bounds the ramification
involved in a lift to characteristic zero. \S \ref{K3 periods} 
recalls the results of Ogus and others
on periods of supersingular K3 surfaces and adapts them to our context.
\S\ref{automorphisms} contains lemmas on automorphism groupschemes.
\S\ref{stacks of surfaces} gives the main results on the existence
and structure of the period map; they are derived as easy consequences of 
the general results on stacks and the geometrical results of the other sections.
In \S \ref{12A1} we examine a particular open substack of $\sE_{uni}$, 
consisting of those surfaces where the singular locus of the $K3$ cover 
is $12\times A_1$.

\begin{acknowledgments}
The surviving authors are very grateful to Derek Holt 
for the Magma routine 
that he supplied.
\end{acknowledgments}

\begin{section}{Results about stacks}\label{stacktology}
This section is devoted to some general results about algebraic stacks.
They are of the kind found in the early chapters of [LMB00]; the
main definition that is not found there is that of a \emph{relative groupscheme}.
This construction introduces problems of coherence, 
which we solve in the context of representable
morphisms by using the notion of a ``local construction'' ([LMB00], ch. 14).
These results
will be used in the description given subsequently of the various stacks
of Enriques surfaces and their double covers. There will be a fixed noetherian
base scheme $B$; this will be excellent and have perfect
residue fields in the geometric applications. The algebraic stacks and spaces
under consideration will be over $B$; that is, they will have morphisms
to $(Aff/B)$.

We start with some remarks about fibre products. Suppose that
$F':\sY'\to\sX$ and $G':\sZ'\to\sX$ are morphisms of stacks. Then the fibre
product $\sY'\times_{\sX}\sZ'$ has no natural morphism to $\sX$, even if
$F'$ and $G'$ are representable. Rather,
there are two obvious such morphisms, namely $F'\circ pr_1$ and $G'\circ pr_2$,
and a $2$-isomorphism between them. This leads to difficulties in the definition
of group objects. However, a representable morphism $F':\sY'\to\sX$ is equivalent
to the datum of a ``local construction 
[of a sheaf of algebraic spaces] 
on $\sX$'', $\underbar{\sY}:\sX\to Sp_B$, where 
$Sp_B$ is the category of algebraic spaces over $B$ (\c[LMB00], ch. 14).
The basic properties of this are summarized in the next two lemmas.

\begin{lemma}\label{local constructions} \part[i] The representable morphism
$F':\sY'\to\sX$ defines a local construction $\underbar{\sY}:\sX\to Sp_B$.

\part[ii] A local construction $\underbar{\sY}:\sX\to Sp_B$ defines
a representable morphism $F:\sY\to\sX$.

\part[iii] If $\sY$ is obtained from $\sY'$ by composing these constructions,
then $\sY$ is $1$-isomorphic to $\sY'$. 

\part[iv] The category $LC/\sX$ of local constructions on $\sX$ is isomorphic 
(not just equivalent) to a full subcategory $Rep_{lc}/\sX$ of the category
$Rep/\sX$ of representable morphisms $\sY\to\sX$, and is equivalent to $Rep/\sX$.

\part[v] $LC/\sX$ has fibre products.

\part[vi] $Rep_{lc}/\sX$ has fibre products, whereas $Rep/\sX$ does not, as already 
remarked.

\part[vii] If $\underbar{\sY}:\sX\to Sp_B$ and $\underbar{\sZ}:\sX\to Sp_B$
are local constructions on $\sX$ and
$\underbar{\sW}=\underbar{\sY}\times \underbar{\sZ}$, leading variously to
representable morphisms
$\sY\to\sX$, $\sZ\to\sX$ and $\sW\to\sX$ in $Rep_{lc}/\sX$,
then $\sW$ is isomorphic, and not just
$1$-isomorphic, to $\sY\times_{\sX}\sZ$.

\part[viii] A morphism $\underbar{\sY}\to \underbar{\sZ}$ of local constructions,
induces a morphism $\sY\to\sZ$ in $Rep_{lc}/\sX$ that is strict; that is, the
composite $\sY\to\sZ\to\sX$ is equal to, and not just $2$-isomorphic to, the
morphism $\sY\to\sX$ coming from the local construction.
\noproof
\end{lemma}

\begin{lemma} \part[i] If $\underbar{\sG}$ is a local construction
on $\sX$, then there is an obvious notion of a group structure 
(namely, a multiplication morphism
$\underbar{m}:\underbar{\sG}\times\underbar{\sG}\to\underbar{\sG}$,
an inverse morphism
$\underbar{i}:\underbar{\sG}\to\underbar{\sG}$ and an identity morphism
$\underbar{e}:\underbar{\sX}\to\underbar{\sG}$, where $\underbar{\sX}$
is the trivial local construction on $\sX$) on $\underbar{\sG}$.
We say that $\underbar{\sG}$ is a \emph{local construction of a 
relative group scheme on $\sX$}.

\part[ii] If $m:\sG\times_{\sX}\sG\to \sG$ etc. are the morphisms obtained from
$\underbar{m}$ etc., then they define a group structure on the representable morphism
$\sG\to\sX$. In particular, the various diagrams defining associativity etc. are
commutative, not merely $2$-commutative. We say that $\sG\to\sX$ is a
\emph{relative groupscheme}.

\part[iii] If $\underbar{\sG}$ is a \emph{local construction of a 
relative group scheme on $\sX$}.
and $\underbar{\sZ}$ a local construction on $\sX$, then there is an obvious notion of 
an action of $\underbar{\sG}$ on $\underbar{\sZ}$.

\part[iv] An action of $\underbar{\sG}$ on $\underbar{\sZ}$ leads to an action
of $\sG$ on $\sZ$ in which the various diagrams are
commutative, not merely $2$-commutative.
\noproof
\end{lemma}

In short, if we restrict to the category $Rep_{lc}/\sX$,
then difficulties concerning fibre products
and coherence in the definitions of group objects and their actions evaporate.
So we shall assume, usually tacitly, that representable morphisms to $\sX$
are in $Rep_{lc}/\sX$. When we construct a representable morphism to $\sX$ it
will usually be clear that the construction lies in $Rep_{lc}/\sX$.

\begin{definition} Suppose that $\sG$ acts on $\sZ$ in $Rep_{lc}/\sX$.
Then $\sZ\to\sX$ is a \emph{pseudo-torsor} if the action induces
an isomorphism $\sG\times_{\sX}\sZ\to\sZ\times_{\sX}\sZ$. It is a
\emph{torsor} if also $\sZ\times_{\sX}X$ is a torsor under $\sG\times_{\sX}X$
for all spaces $X\to\sX$.
\end{definition}

\begin{proposition} Suppose that $F:\sX\to\sY$ is a morphism of algebraic $B$-stacks.
Then there is a relative groupscheme $\phi:\sG=\sG_F\to\sX$, the \emph{stabilizer of $F$},
where, for $U\in Ob(Aff/B)$,
the objects of $(\sG_F)_U$ are pairs $(x,x\stackrel{a}{\to}x)$ with $F(a)=1_{F(x)}$
and the morphisms $(x,a)\to(y,b)$ are those morphisms $x\stackrel{g}{\to}y$
with $bg=ga$ and $\phi(x,a)=x$.

\begin{proof} We exhibit the local construction $\underbar{\sG}$. This is
defined by taking $\underbar{\sG}(x)$ to be the algebraic space that represents 
the sheaf $\sI som(x,x)=\sA ut(x)$ on $(Aff/U)$, where $x\in Ob(\sX_U)$.
\end{proof}
\end{proposition}

\begin{remark} Without making a local construction, things are a little less tidy.
Consider, for example, this approach. 

The diagonal
$\Delta_F=\Delta_{\sX/\sY}:\sX\to\sX\times_{\sY}\sX$ is representable. So the 
morphisms
$$
pr_i:\sG':=\sX\times_{\Delta_{\sX/\sY},\sX\times_{\sY}\sX,\Delta_{\sX/\sY}}\sX\to\sX,
$$ for $i=1,2$,
are representable. It is easy to prove that $\sG'$ is $1$-isomorphic to
the stabilizer groupscheme $\sG_F$ constructed above. Putting this aside,
there is a groupoid structure on $\sG'$ over $\sX$,
and one can turn this into a group structure, in some sense, using the existence of
a $2$-isomorphism $pr_1\Rightarrow pr_2$. However, it is precisely the fact
that we have a $2$-isomorphism rather than an equality of morphisms
that leads to problems of coherence.

Note that, if $\sX=[X/R]$ for an algebraic groupoid
$R=X\times_{\sX}X\rightrightarrows X$, then $\sG\to\sX$ pulls back over $X$ to the 
stabilizer groupscheme that is the restriction of $R\to X\times_B X$
to the diagonal in $X\times_BX$.
\end{remark}

\begin{lemma} Given $\sX\stackrel{F}{\to}\sY\stackrel{G}{\to}\sZ$, 
there is an exact sequence
$$1\to\sG_{\sX/\sY}\stackrel{A}{\to}\sG_{\sX/\sZ}\stackrel{B}{\to}\sG_{\sY/\sZ}\times_{\sY}\sX$$
of relative groupschemes over $\sX$. The morphisms $A,B$ are strict.
\begin{proof} This is an immediate consequence of the definitions
in terms of local constructions.
\end{proof}
\end{lemma}

\begin{lemma}\label{algebraicity criterion}
Suppose that $F:\sX\to\sY$ is a morphism of stacks, that $\sY$ 
is algebraic and that for all
spaces $Y$ and morphisms $Y\to\sY$, the fibre product $\sX\times_{\sY}Y$
is $1$-isomorphic to
an algebraic stack. Then $\sX$ is algebraic.

\begin{proof} For relative representability of $\sX$, 
suppose that $V,W$
are algebraic spaces and $v:V\to\sX$, $w:W\to\sX$ are morphisms. There is 
a lift $\tw=(w,1):W\to\sX\times_{\sY}W$ of $w$. Then there is a $1$-isomorphism
$$
V\times_{\sX}W\to (V\times_{\sX}(\sX\times_{\sY}W))\times_{\sX\times_{\sY}W,\tw}W.$$
Since $V\times_{\sY}W$ is representable and $\sX\times_{\sY}W$ is
an algebraic stack, it follows that the right hand term is representable,
which is enough.

To find a smooth presentation of $\sX$, take a smooth presentation
$Y\to\sY$ of $\sY$. Then $\sX\times_{\sY}Y\to\sX$ is smooth and surjective,
so that a smooth presentation $X\to\sX\times_{\sY}Y$ will suffice.
\end{proof}
\end{lemma}

The next result generalizes the fact that the classifying stack of a 
flat groupscheme is algebraic (\c[LMB00], Proposition (10.13.1)).

\begin{lemma}\label{classifying stacks}
Suppose that $\sG\to\sX$ is a flat relative groupscheme. Denote
the corresponding local construction by ${\underbar{\sG}}$. Then there
is an algebraic stack $\sB\sG$, the \emph{classifying stack},
with a structural morphism
$\pi:\sB\sG\to\sX$ and a tautological section $s:\sX\to\sB\sG$, 
defined by the property that for $U\in Ob(Aff/B)$, the objects
of $(\sB\sG)_U$ are pairs $(x\in Ob(\sX_U),P\to U)$,
where $P\to U$ is a torsor under the groupscheme ${\underbar{\sG}}(x)\to U$,
$\pi(x,P)=x$ and $s(x)=(x, {\underbar{\sG}}(x)\to U)$.
Moreover, $\pi:\sB\sG\to\sX$ is smooth.
\begin{proof} 
It is clear that $\sB\sG$ is a $B$-stack and that for all morphisms
$x:X\to\sX$ with $X$ an algebraic space, $\sB\sG\times_{\sX}X\to X$
is $1$-isomorphic to $B(\sG\times_{\sX}X)\to X$. So $\sB\sG$ is algebraic,
by \ref{algebraicity criterion}.

The smoothness is local on $\sX$, so we can take $X$ to be a space.
In this case the smoothness is well known, although not to be found explicitly
(by us) in \c[LMB00]. It follows from the statement that if
$\sX\stackrel{F}{\to}\sY\stackrel{G}{\to}\sZ$ are $1$-morphisms of
algebraic stacks such that $GF$ is smooth and surjective and $F$
is flat, then $G$ is smooth. After pulling back by various smooth presentations
this reduces to the analogous statement for algebraic spaces, where it
is even more well known: to get the smoothness of $BG\to X$, take
$\sX=X=\sZ$ and $\sY=BG$, so that $\sX\stackrel{F}{\to}\sY$ is a 
$G$-torsor, so flat, and $GF$ is the identity.
\end{proof}
\end{lemma}

Various useful constructions involving groupschemes carry over to
relative groupschemes in an unsurprising way.
For example, it is well known that for a group $G$ over a space $X$,
representable morphisms
$\sW\to BG$ correspond to $G$-spaces $Z\to X$; given $Z$, $\sW$ is
$[Z/G]$ and, given $\sW$, $Z$ is $X\times_{s,BG}\sW$. 

\begin{lemma}\label{stacky quotients}\label{1.8}
Assume that $\sG\to\sX$ is a flat relative groupscheme.
Then for any representable $a:\sZ\to\sX$ with a $\sG$-action, there
is a representable morphism $b:\sW\to\sB\sG$ such that for all $X\to\sX$,
$\sW\times_{\pi\circ b,\ \sX}X$ is $1$-isomorphic to the quotient stack
$[\sZ\times_{\sX}X/\sG\times_{\sX}X]$.
Conversely, every representable $a:\sZ\to\sX$ with a $\sG$-action arises in
this way, up to $1$-isomorphism.
\begin{proof}
To construct $\sW$ it is enough to produce the corresponding local construction
$\underbar{\sW}$ on $\sB\sG$. So assume that $\sZ\in Rep_{lc}/\sX$.
Now an object $\tx$ of $(\sB\sG)_U$ is a pair
$\tx=(x\in Ob(\sX_U), P\to U)$, where $P\to U$ is a torsor under
$\underbar{\sG}(x)\to U$. We have a $U$-space $\underbar{\sZ}(x)$, and we
define $\underbar{\sW}(\tx)$ to be the quotient $U$-space
$(\underbar{\sZ}(x)\times_U P)/\underbar{\sG}(x)$, where
$\underbar{\sG}(x)$ acts diagonally.

There is a diagram
$$
\xymatrix{
{\sZ} \ar[r] \ar[d]_a & {\sW} \ar[d]^b \\
{\sX} \ar[r]^s & {\sB\sG} \ar[r]^{\pi} & {\sX}
}
$$
in which the square is clearly commutative, not just $2$-commutative.
To show that it is $2$-Cartesian, it is enough to show
that it is Cartesian after base change by every morphism $X\to\sX$, where
$X$ is a space. But this is just the statement of the lemma for classifying
stacks over a space, which is, as already stated, well known to be true.

The converse is immediate.
\end{proof}
\end{lemma}

We write $\sW=[\sZ/\sG]$. This agrees with the usual notation, in the sense
that if we take the Cartesian product of the diagram above by a morphism
$X\to\sX$ when $X$ is a space, then $[\sZ/\sG]\times_{\sX}X$ is $1$-isomorphic
to $[Z/G]$, where $Z$ is a space $1$-isomorphic to $\sZ\times_{\sX}X$
and $G$ is a groupscheme $1$-isomorphic to $\sG\times_{\sX}X$.

\begin{proposition} Suppose that $\sG\to\sX$ is a 
flat relative groupscheme and that
$F:\sW\to\sX$ is a morphism of algebraic stacks. Suppose that $f:\sP\to\sW$
is a torsor under $\sG_{\sW}=\sG\times_{\sX}\sW$. 
Then there is a morphism $H:\sW\to\sB\sG$ under which $\sP$ is $1$-isomorphic
to $\sW\times_{H,\sB\sG,s}\sX$.

\begin{proof} It is clear that there is a $2$-factorization of $f$ through
$[\sP/\sG_{\sW}]$, and then, from the fact that $\sP\to\sW$ is a torsor, 
that the resulting morphism
$[\sP/\sG_{\sW}]\to\sW$ is a $1$-isomorphism.
So, by \ref{stacky quotients}, there is a morphism $t:\sW\to\sB(\sG_{\sW})$
and a $2$-Cartesian square
$$
\xymatrix{
{\sP} \ar[r] \ar[d] & {\sW} \ar[d]^t \\
{\sW} \ar[r]^{s_{\sW}} & {\sB(\sG_{\sW})}
}
$$
where $s_{\sW}$ is the tautological section.
Now $\sB(\sG_{\sW})$ is $1$-isomorphic to $\sB\sG\times_{\sX}\sW$,
so there is another $2$-Cartesian square
$$
\xymatrix{
{\sW} \ar[r]^{s_{\sW}} \ar[d] & {\sB(\sG_{\sW})} \ar[d]\\
{\sX} \ar[r]^{s_{\sX}} & {\sB\sG;}
}
$$
putting these squares together gives the result.
\end{proof}
\end{proposition}

It is convenient to ``stackify'' the construction of geometric quotients, at 
first in the context of group actions.

\begin{proposition}\label{quotients in families}\label{summary}
Suppose that $\sG\to\sX$ is a flat relative groupscheme
acting properly on the representable morphism $\sZ\to\sX$. Then there is a
a $2$-commutative
diagram
$$\xymatrix{
\label{CD:no1}
{\sZ} \ar[r] \ar[d]_a & {[\sZ/\sG]} \ar[r]\ar[d]^b & {\sZ/\sG} \ar[d] \\
{\sX} \ar[r]_s & {\sB\sG} \ar[r]_{\pi} & {\sX}
}$$
whose left hand square is $2$-Cartesian and whose vertical arrows are
representable. The morphism $(\sZ/\sG)\to\sX$ has the property
that for all spaces $X$ and \emph{flat} morphisms
$X\to\sX$, the fibre product $(\sZ/\sG)\times_{\sX}X$ is $1$-isomorphic to the geometric
quotient $(\sZ\times_{\sX}X)/(\sG\times_{\sX}X)$.
\begin{proof} The fundamental result of [KM97] is that proper actions
of flat groupschemes have geometric quotients, and the formation of these quotients 
commutes with all \emph{flat} base change. Then for any $Y\to\sX$, 
$(\sZ/\sG)\times_{\sX}Y$ is constructed by first computing 
$(\sZ\times_{\sX}X)/(\sG\times_{\sX}X)$ for flat $X\to\sX$, 
pulling this back to $Y\times_{\sX}X$
and then making a flat descent to get an algebraic space over $Y$. 
This gives a local construction
$\underbar{\sZ/\sG}$ on $\sX$, which leads in turn to the existence of $\sZ/\sG$.
\end{proof}
\end{proposition}

Of course, if $\sG$ acts freely on $\sZ$, then $[\sZ/\sG]\to\sZ/\sG$
is a $1$-isomorphism, but not otherwise.

It is natural and useful in this context to introduce gerbes.
\begin{definition} An \emph{algebraic gerbe} is an
epimorphism $F:\sX\to\sY$ of algebraic stacks 
such that the diagonal morphism 
$\Delta_{\sX/\sY}:\sX\to\sX \times_{\sY}\sX$ is also an epimorphism.
A gerbe is \emph{neutral} if it has a $2$-section.
\end{definition}

(Cf. \c[LMB00], D\'ef. 3.15, p. 22, where ``gerbe'' is defined when $\sY$
is a space.)
The meaning of this definition is
that for all schemes $Y$ and morphisms $Y\to\sY$, there is a faithfully
flat $Y'\to Y$ such that $Y'$ lifts to $\sX$ and for all schemes $X$ and
morphisms $v:X\to\sX$, $w:X\to\sX$, the morphism
$X\times_{v,\sX,w}X\to X\times_{fv,\sY,fw}X$ is surjective as a map of sheaves on $X$,
where both sides are regarded as $X$-schemes via the first projection.

\begin{lemma} Suppose that $\sZ\to \sY$ is any morphism.
and $\sX\to\sY$ is an algebraic gerbe. 
Then $\sX\times_{\sY}\sZ\to\sZ$ is an algebraic gerbe.
\noproof
\end{lemma}

\begin{lemma}\label{gerbes with flat stabilizer are smooth} 
Suppose that $F:\sX\to\sY$ is an
algebraic gerbe and that $\sG_{\sX/\sY}\to\sX$ is flat. Then $F$ is smooth.

\begin{proof} Note first that the hypotheses on $F$ and $\sG_{\sX/\sY}$ are preserved
under any base change $\sZ\to\sY$.

Pick a flat presentation $Y\to\sY$.
By assumption on $F$, there is a flat surjective map $Y'\to Y$ such that $Y'\to\sY$
lifts to $\sX$. Put $\sX'=\sX\times_{\sY} Y'$. Then $\sX'\to Y'$ has a section, say $s$,
so that, by [LMB00], Lemme 3.21, $\sX'$ is naturally identified with
$B(G'/Y')$, where $G'=\sG_{\sX/\sY}\times_{\sX,s}Y'$. Since,
by assumption, $G'\to Y'$ is flat,
it follows that $B(G'/Y')\to Y'$, and so $\sX'\to Y'$, is smooth,
by \ref{classifying stacks}.
Now $Y'\to\sY$ is a flat presentation, and so $F$ is smooth.
\end{proof}
\end{lemma}

It is not clear whether a smooth gerbe necessarily has a flat stabilizer.
This amounts to whether a groupscheme $K\to Y$ is necessarily flat if
$BK$ is algebraic and smooth over $Y$.

The next result extends Lemme 3.21 of [LMB00].

\begin{proposition}\label{neutral gerbes are classifying stacks}
Suppose that $F:\sX\to\sY$ is a neutral gerbe, with
section $s$. Define $\sK\to\sY$ by $\sK=\sG_F\times_{\sX,s}\sY$.
Then $\sX$ is $1$-isomorphic to $\sB\sK$ over $\sY$.

\begin{proof} We want to show that $s:\sY\to\sX$ is a torsor over $\sX$
under $\sK\times_{\sY,F}\sX\to\sX$. For this, we need an action morphism
$\sigma:(\sK\times_{\sY}\sX)\times_{\sX,s}\sY\to\sY$ that gives an isomorphism
$\Sigma:(\sK\times_{\sY}\sX)\times_{\sX,s}\sY\to\sY\times_{\sX}\sY$. 
Note that $(\sK\times_{\sY}\sX)\times_{\sX,s}\sY$ is $1$-isomorphic
to $\sK$, and then the only possible choice for $\sigma$ is, up to
$2$-isomorphism, the structure
morphism. Now it is shown in [LMB00] that, when $\sY$ is a space,
there is a $1$-isomorphism $\sX\to\sB\sK$ under which
$s:\sY\to\sX$ is $2$-isomorphic to the natural morphism $\sY\to\sB\sK$,
so that $\Sigma$ is an isomorphism everywhere locally on $\sY$.
Hence $\Sigma$ is a $1$-isomorphism and $s:\sY\to\sX$ is a torsor 
under $\sK\times_{\sY}\sX$. This gives a morphism $\sX\to\sB\sK$
over $\sY$. It is shown in [LMB00] that this is an isomorphism
everywhere locally on $\sY$, and so is a $1$-isomorphism.
\end{proof}
\end{proposition}

\begin{lemma}\label{quotient lemma}
Suppose that $R\stackrel{p,q}{\rightrightarrows} S$ is a smooth (resp. flat,
resp. Cohen--Macaulay, resp.\ ...) algebraic groupoid
and that there is a closed $S$-flat normal subgroupscheme $H\to S$ of the 
stabilizer groupscheme
$R\vert_{\Delta_S}\to S$. 
Then there is a quotient $H\backslash R$ that is also a smooth (resp. flat, resp.\ ...)
algebraic groupoid over $S$.

\begin{proof} Assume that $R\stackrel{p,q}{\rightrightarrows} S$ is smooth;
the proof is the same for any of the other local properties.
Put $G=R\vert_{\Delta_S}$.
There is a left action of $G\to S$ on $p:R\to S$, defined by a morphism
$a:G\times_{S,p}R\to R$. 
Via the natural isomorphisms
\begin{displaymath}
G\times_{S,p}R\cong G\times_{S,pr_1\circ j}R\cong G\times_{S,pr_1}(S\times S)_{S\times S,j}R,
\end{displaymath}
this is identified with a left action of $pr_1^*G\to S\times S$ on $j:R\to S\times S$.
Note that this left action makes $R$ into a pseudo--torsor under $pr_1^*G$
over $S\times S$.
In particular, these actions are free
and, since $H/S$ is flat, there is a geometric quotient $\pi:R\to R_1=H\backslash R$ 
and a factorization $j=j_1\circ \pi$ [Ar74a], Cor. 6.3. 
Via the identification of the two actions just mentioned,
we regard $R_1$ indifferently as $H\backslash R$ or $pr_1^*H\backslash R$.

For smoothness, regarding $R_1$ as $H\backslash R$ shows that the
morphism $p_1:R_1\to S$ induced by $p$ is smooth. Next,
let $i:R\to R$ and $c:R\times_{p,S,q}R\to R$ denote
the inverse, resp. composition, morphisms.
To show that $i$ descends to $i_1:R_1\to R_1$,
note that for $g\in H(s)$ and $f\in R(s,t)$, we have $g(f)=a(g,f)=f\circ g^{-1}$,
so that $i(g(f))=g'(i(f))$, where $g'=c(f^{-1},g(f))$. So $i$ preserves $H$-orbits,
and so descends to $i_1$. As for $c$, we identify $c$ with a morphism
$c':R\times_{p,S,p}R\to R$ via the isomorphism
$(1,i):R\times_{p,S,q}R\to R\times_{p,S,p}R$. Then $H\times_S H$ acts on
$R\times_{p,S,p}R$, and we must check that $c'$ takes $H\times_H$-orbits to $H$-orbits.
For this, note that
$$c'((g_1,g_2)(f,h))=fg_1^{-1}g_2h = g'(c'(f,h)),$$
where $g'=hg_2^{-1}g_1h^{-1}$, and so there is a composition $c_1:R_1\times R_1\to R_1$.
Finally, define $q_1:R_1\to S$ by $q_1=p_1\circ i_1$ and $j_1=(p_1,q_1):R\to S\times B$.
So $R_1$ is a smooth groupoid over $S$.
\end{proof}
\end{lemma}

\begin{remark} One crucial step in the construction given in [KM97]
of geometric quotients is the formation of the quotient groupoid, in certain
circumstances, of a groupoid by a normal subgroupoid. However, it is
not clear that \ref{quotient lemma} is a special case of this, because
it is not clear that the subgroupscheme $H$ arises as the stabilizer
of a subgroupoid of $R$.
\end{remark}

The next result restates this lemma in terms of algebraic stacks. 

\begin{definition} Suppose that $F:\sX\to\sY$ is an algebraic gerbe and that
$\sG=\sG_{\sX/\sY}\to \sX$ is flat.
Then $F:\sX\to\sY$ is a
\emph{quotient by $\sG$} or an \emph{extension by $\sG$}.
\end{definition} 

\begin{remark} For example, if $\sH\to\sY$ is a flat relative groupscheme,
then $\pi:\sB\sH\to\sY$ is an extension by $\pi^*\sH$.
In view of \ref{neutral gerbes are classifying stacks}
and the idea that the classifying stack $BH\to Y$ of a groupscheme
$G\to Y$ is ``the quotient of $Y$ by $H$'', the definition of ``quotient''
given above might seem
to be made the wrong way round. However, it will turn out to be the right 
way round for our analysis of moduli.
\end{remark}

\begin{proposition}\label{quotient stacks}\label{extensions}
Suppose that $\sX$ is an algebraic stack and that $\sH\to\sX$
is a normal flat closed relative subgroupscheme
of the stabilizer groupscheme $\sG_{\sX/B}\to\sX$.
Then there is an essentially unique smooth
algebraic gerbe $F:\sX\to\sY$ 
that is an extension (or quotient) by $\sH$. 

\begin{proof} Pick a flat presentation $S\to X$ and put 
$R=R_S:=S\times_{\sX}S\rightrightarrows S$, so that $\sX\cong [S/R]$.
Then $\sH\times_{\sX}S$ is a normal flat subgroupscheme of the $S$-groupscheme
$(R\to S\times_B S)\vert_{\Delta_{S\to B}}$. 
By \ref{quotient lemma} $\sH_S\backslash R$ is a flat groupoid over $S$. Then define
$\sY =[S/\sH_S\backslash R]$. It is clear that the obvious morphism
$F:\sX\to\sY$ is a smooth algebraic gerbe and is a quotient by $\sH$.

If $T\to \sX$ is another flat presentation, then,
after replacing $T$ by $T\times_{\sX}S$ if necessary, we can suppose that $T\to \sX$
factors through a flat morphism $T\to S$. Then 
$R_T\cong R_S\times_{S\times_B S}T\times_B T$,
the restriction of $R_S$ to $T$, and $\sH_T=\sH_S\times_S T$. The uniqueness
of $\sX\to\sY$ is now clear.
\end{proof}
\end{proposition}

This justifies the next definition.

\begin{definition}
If $F:\sX\to\sY$ is an algebraic gerbe whose stabilizer $\sG_F$ is flat over $\sX$,
then $F:\sX\to\sY$ is \emph{the} quotient of $\sX$ by
$\sH$ or \emph{the extension by $\sG_F$}, and we write $\sY=\sX/\sH$.
\end{definition}

\begin{lemma}If also $\sX\to\sY$ is a quotient by a flat representable
algebraic relative groupscheme $\sG\to \sX$, then $\sX\times_{\sY}\sZ\to\sZ$
is a quotient by the pullback of $\sG$ to $\sX\times_{\sY}\sZ$.
\noproof
\end{lemma}

\begin{proposition}\label{universal property of extensions}
Suppose that $F:\sX\to\sY$ is the extension by $\sG_F$
and that $H:\sX\to\sZ$ is a morphism such that $\sG_F$ is $1$-isomorphic
to a normal closed subgroupscheme of $\sG_H$. Then there is a morphism
$G:\sY\to\sZ$ and a $2$-isomorphism $G\circ F\Rightarrow H$.

\begin{proof} Note first that since $\sG_F\to\sX$ and $\sG_H\to\sX$ are 
representable, we can suppose that they are objects of $Rep_{lc}/\sX$,
and then that $\sG_F$ is isomorphic to a normal closed subgroupscheme of $\sG_H$.

Pick a smooth presentation $Z_0\to\sZ$ and a smooth presentation
$X_0\to\sX\times_{\sZ}Z_0$. Then $X_0\to\sX$ is a smooth presentation, so that
$\sX$ is $1$-isomorphic to $[X_0/X_1]$ and $\sZ$ to $[Z_0/Z_1]$, where
$X_1=X_0\times_{\sX}X_0$ and $Z_1=Z_0\times_{\sZ}Z_0$. So there is a commutative
diagram
$$\xymatrix{
{X_1} \ar[r] \ar[dr]_f & {Z_1\times_{Z_0\times_B Z_0}(X_0\times_B X_0)}
\ar[r] \ar[d]^h & {Z_1} \ar[d]_g\\
& {X_0\times_B X_0} \ar[r] & {Z_0\times_B Z_0}
}$$
and $\sX\to\sZ$ factors through $\sX\to [X_0/Z_1\times_{Z_0\times_B Z_0}(X_0\times_B X_0)]$.
Also,
$f$ is a pseudo-torsor under $pr_1^*(\sG_{\sX/B}\times_{\sX}X_0)$
and $g$ is a pseudo-torsor under $pr_1^*(\sG_{\sZ/B}\times_{\sZ}Z_0)$.
So $h$ is a pseudo-torsor under $pr_1^*((\sG_{\sZ/B}\times_{\sZ}\sX)\times_{\sX} X_0)$.

Now $\sG_H$ is the kernel of $\sG_{\sX/B}\to \sG_{\sZ/B}\times_{\sZ}\sX$,
so that $X_1 \to Z_1\times_{Z_0\times_B Z_0}(X_0\times_B X_0)$ factors through
the quotient morphism $X_1\to X_1/pr_1^*(\sG_H\times_{\sX} X_0)$.
Since $\sY$ is $1$-isomorphic to $[X_0/(X_1/pr_1^*(\sG_H\times_{\sX} X_0))]$,
it follows that $\sX\to [X_0/Z_1\times_{Z_0\times_B Z_0}(X_0\times_B X_0)]$ 
factors through $\sX\to [X_0/(X_1/pr_1^*(\sG_H\times_{\sX} X_0))]$. 
Since $\sG_F$ is a closed subgroupscheme of $\sG_H$, 
$\sX\to [X_0/(X_1/pr_1^*(\sG_H\times_{\sX} X_0))]$ factors through
$\sX\to [X_0/(X_1/pr_1^*(\sG_F\times_{\sX} X_0))]$. However, this last
morphism is exactly the extension by $\sG_F$, as constructed in terms of
algebraic groupoids in \ref{extensions}.
\end{proof}
\end{proposition}

\begin{proposition} Suppose that $F:\sX\to\sY$ is an algebraic gerbe.
Then one of $\sX$ and $\sY$
has a geometric quotient if and only if the other does,
and when the quotients exist they are naturally isomorphic.

\begin{proof}
This follows from the fact that $\sX$ and $\sY$ have the same associated sheaves of
connected components. (Cf. [LMB00], Lemme 3.18.)
\end{proof}
\end{proposition}

\begin{proposition}\label{closed substacks of extensions}
\label{closed substacks of smooth gerbes}
Suppose that $F:\sX\to\sY$ is the quotient of algebraic
stacks by a flat relative groupscheme $\sG\to\sX$. Then every closed
substack $\sZ$ of $\sX$ is of the form $\sW\times_{\sY}\sX$ for a unique closed
substack $\sW$ of $\sY$ that is identified with the quotient of $\sZ$
by $\sG\times_{\sX}\sZ$.

\begin{proof} For any $S\to \sX$, put $G_S =\sG\times_{\sX}S$.
Pick a flat
presentation $f:X_0\to\sX$ and put 
$X_1=X_0\times_{\sX}X_0\stackrel{p,q}{\rightrightarrows}X_0$.
Then $\sX$ is isomorphic to $[X_0/X_1]$
and $\sY$ is isomorphic to $[X_0/(X_1/p^*G_{X_0})]$.

Now suppose that $i:\sZ\to\sX$ is a closed embedding. Then $Z_0=X_0\times_{\sX}\sZ$
is a closed algebraic subspace of $X_0$. Put
$Z_1=Z_0\times_{\sZ}Z_0\stackrel{s,t}{\rightrightarrows} Z_0$; then $Z_1$ is
isomorphic to $X_1\times_{\sX}\sZ$, since $i$ is representable.
Also $\sZ$ is isomorphic to $[Z_0/Z_1]$. Put $\sW=[Z_0/(Z_1/s^*G_{Z_0})]$,
so that $\sW=\sZ/\sG\times_{\sX}\sZ$ and there is a $2$-commutative square
$$
\xymatrix{
{\sZ} \ar[r]^h \ar[d]_i & {\sW} \ar[d]^j\\
{\sX} \ar[r]^F & {\sY}
}
$$
where $h$ is a quotient by $\sG\times_{\sX}\sZ$.
The questions of whether this square is $2$-Cartesian and $j$ is a closed embedding
are local on $\sY$. Since $F$ is
a flat epimorphism, the composite $g=F\circ f:X_0\to\sY$ is a flat presentation.
After pulling back by $g$, we can suppose that $\sY$ is identified with $X_0$
and that $F$ has a section. Then ([LMB00], Lemme 3.21) $F:\sX\to X_0$ is identified
with $B(G_{X_0}/X_0)\to X_0$. Then $X_1=G_{X_0}$ and $p,q$ are both the structure
morphism, and the groupoid $X_1/p^*G_{X_0}$ is trivial because 
$(p,q):X_1\to X_0\times X_0$ is a pseudo-torsor under $p^*G_{X_0}\to X_0\times X_0$.
Pulling back under $Z_0\times Z_0\to X_0\times X_0$ shows that
$(s,t):Z_1\to Z_0\times Z_0$ is a pseudo-torsor under $s^*G_{Z_0}\to Z_0\times Z_0$.
Then $\sW$ is identified with $Z_0$ and $j$ with the inclusion.
Finally, the groupoid $Z_1\rightrightarrows Z_0\times Z_0$ is identified with the groupoid 
associated to the groupscheme $G_{Z_0}\to Z_0$.
\end{proof}
\end{proposition}

\begin{proposition}\label{summary sequel} With the previous notation and assumptions,
suppose also that $\sH\to [\sZ/\sG]$ is a closed normal flat relative subgroupscheme
of the relative stabilizer groupscheme $\sG_{[\sZ/\sG]\to\sX}$. 
Suppose that $[\sZ/\sG]\to\sE$ is the 
extension by $\sH$. Then $[\sZ/\sG]\to\sZ/\sG$ factors uniquely, up to $2$-isomorphism,
through $\sE$.
\begin{proof} This is a special case of \ref{universal property of extensions}.
\end{proof}
\end{proposition}

We expand upon this in a special case.

\begin{lemma} Suppose that $F:\sX\to\sY$ is a morphism of algebraic stacks,
that $\sY$ is Deligne-Mumford and that whenever $H:\sX\to\sZ$ is a morphism to
a Deligne-Mumford stack there is a morphism $G:\sY\to\sZ$ and a $2$-isomorphism
$GF\Rightarrow H$. Then for any $F':\sX\to\sY'$ with the same property
there is a $1$-isomorphism $\phi:\sY\to\sY'$ and a $2$-isomorphism
$\phi\circ F\Rightarrow F'$.
\noproof
\end{lemma}

\begin{definition} \part In this case
we say that $\sY$ is a \emph{Deligne-Mumford quotient of $\sX$} and
denote $\sY$ by $\sX_{DM}$.

\part A \emph{Keel-Mori stack} is an
algebraic stack whose stabilizer groupscheme is finite.
\end{definition}

Recall that Keel-Mori stacks possess geometric quotients [KM97].
It seems plausible that any Keel-Mori stack should possess a Deligne-Mumford 
quotient, although our attempts to prove this have failed.

\begin{proposition}\label{DM quotients} Suppose that $\sX$ is an algebraic stack and
there is a flat closed subrelative groupscheme
$\sH\to\sX$ of $\sG_{\sX/B}\to\sX$ such that $\sH$ has connected geometric fibres and
$\sG_{\sX/B}/\sH\to\sX$ is unramified. Then $\sH$ is normal in $\sG_{\sX}$
and $\sX/\sH$ is a Deligne--Mumford quotient of $\sX$.

\begin{proof} It is well known that if $G\to S$ is a groupscheme and
$H$ is a closed flat subgroupscheme of $G$ with connected geometric fibres
and unramified quotient $G/H\to S$, then $H$ is normal in $G$.
The normality of $\sH$ is now (even more) immediate.
It is clear from the construction of $\sX/\sH$
in terms of smooth groupoids that it has unramified stabilizer. Then it is a
Deligne--Mumford stack, by Th. 8.1 of \c[LMB00].

For the universal property, suppose that $\sX\to\sZ$ is a morphism to a 
Deligne--Mumford stack. Pick an \'etale presentation $Z\to\sZ$ and put
$T=Z\times_{\sZ}Z\rightrightarrows Z$. 
Pick a smooth presentation $X\to \sX\times_{\sZ}Z$; in particular,
$X\to\sX$ is a smooth presentation. Put $S=X\times_{\sX}X\rightrightarrows X$
and $T\vert_X=T\times_{Z\times_B Z}(X\times_B X)$;
then there is a factorization $S\to T\vert_X\rightrightarrows X$.
Since $T\vert_X\to X\times_B X$ is unramified, $S\to T\vert_X$
factors through $S/\sH_X$, and we are done.
\end{proof}
\end{proposition}

\begin{proposition}\label{DM quotients in families}\label{1.16}
With the assumptions and notation of \ref{summary sequel}
suppose also that $\sX$ is a Deligne--Mumford stack
and that $\sG\to\sX$ is a finite flat relative groupscheme all of whose
geometric fibres are connected. 
Then the sequence of morphisms $[\sZ/\sG]\to\sE\to\sZ/\sG$ 
identifies $\sZ/\sG$ with the Deligne-Mumford quotient of each of the stacks 
$[\sZ/\sG]$ and $\sE$.
\begin{proof} This is certainly true after pulling back by a flat 
presentation of $\sX$. Since the Deligne-Mumford quotient is unique
up to $1$-isomorphism, when it
exists, the result follows.
\end{proof}
\end{proposition}

\ref{quotients in families} has a straightforward generalization; the proof
is similar, so omitted.

\begin{proposition} Suppose that $F:\sX\to\sY$ is a $1$-morphism of
algebraic stacks and that the relative stabilizer groupscheme
$\sG_F\to\sX$ is finite. Then $F$ factors as
$\sX\stackrel{\pi}{\to}{\bar{\sX}}\stackrel{\bar F}{\to}\sY$,
where $\bar F$ is representable and for all spaces $Y$ and \emph{flat}
morphism $Y\to\sY$, the morphism $\pi_Y:\sX\times_{\sY}Y\to{\bar{\sX}}\times_{\sY}Y$
is a geometric quotient.
\noproof
\end{proposition}

\begin{lemma}\label{injectivity}
Suppose that $a:\sX\to\sY$ is a morphism of algebraic stacks
and that $\sX$ is Deligne--Mumford. 
Assume that for all algebraically closed fields $k$ and points 
$x,x':\Spec k\to\sX$ the morphism
$$a_{x,x'}:\Spec k\times_{x,\sX,x'}\Spec k\to \Spec k\times_{ax,\sY,ax'}\Spec k$$
is an injection. Then $a$ is representable.

\begin{proof} 
We can suppose that $\sY$ is representable, and the nit is enough too show that
for all $U\in Ob(Aff/B)$ and for all $\xi\in Ob(\sX_U)$, the sheaf
$\sA ut_{\xi}=\sI som(\xi,\xi)$ on $U$ is trivial. Note that this sheaf
is represented by a groupscheme over $U$, which, since $\sX$ is 
Deligne-Mumford, is finite and unramified.

Pick an arbitrary geometric point $z:\Spec k\to U$. Now $\xi$ is $\xi:U\to \sX$;
put $\xi\circ z=x$. Then $\sI som(\xi,\xi)\times_{U,x}\Spec k$ is identified
with $\Spec k\times_{x,\sX,x}\Spec k$, which, by the assumption of injectivity,
has only one $k$-point. Hence $\sI som(\xi,\xi)$ is trivial.
\end{proof}
\end{lemma}

\begin{proposition}\label{isomorphism} Suppose that $\sX$ and $\sY$ are separated
algebraic stacks and that $a:\sX\to\sY$ is a proper morphism with finite fibres.
Suppose that $\sX$ is Deligne--Mumford, that $\sY$ is normal
and that both are irreducible. Assume also that $a$ is a $1$-isomorphism
over a non--empty open substack of $\sY$ and that for all geometric points
$x,x':\Spec k\to \sX$, the morphism
$$a_{x,x'}:\Spec k\times_{x,\sX,x'}\Spec k\to \Spec k\times_{ax,\sY,ax'}\Spec k$$
is an injection.
Then $a$ is a $1$-isomorphism.
\begin{proof}
By \ref{injectivity} $a$ is representable.
The statement that $a$ is a $1$-isomorphism is local on $\sY$, so we are
reduced to proving that a proper finite morphism of separated normal
irreducible algebraic spaces that is generically an isomorphism
is an isomorphism. This is well known, and trivial.
\end{proof}
\end{proposition}
\end{section}

\begin{section}{Local and global moduli: general results.}
\label{local-to-global}

It is well known that the presence of vector fields complicates questions
about the existence and structure of moduli spaces, and even their
definitions. Our aim in this section is to clarify this; more precisely,
we prove results concerning the local nature of the morphism
from a local moduli space to a global coarse moduli space.

First, recall some basic results of Artin [Ar74a] and Keel and Mori
[KM97].
Fix a Noetherian base scheme $B$ with perfect residue fields.
We shall consider proper families $f:\sX \to S=\Spec R$, where $R=\sO_S$ is an
$\sO_B$-algebra of finite type,
possibly with additional structures, such as a polarization, or cone of
polarizations, or the datum of a line in a space of derivations.
We shall not be explicit about these additional structures,
although there will always be enough to make all $\Isom$ and $\Aut$
functors representable by schemes of finite type. In other words,
the stack $\sF$ under consideration will be algebraic, in Artin's sense [Ar74a].
Fix a point $s\in S$, mapping to $b\in B$, say.
We shall abuse notation by choosing a coefficient ring for $k(s)$ and
denoting it by $W(k(s))$. Of course, if $k(s)$ is perfect (for example, if
$k(s)/k(b)$ is algebraic), then this is the 
ring of Witt vectors and the choices of rings and homomorphisms evaporate.
We shall abuse language by using the phrase ``$s$ is a closed point of $S$'' to mean
``there is a morphism $\Spec k(s) \to\Spec k\to S$'', where $\Spec k\to S$ is a closed
embedding and $k(s)/k$ is algebraic, and we shall
make the running assumption that $k(s)/k(b)$ is algebraic
whenever $s$ is closed. 

According to Artin's definition, $f$ is versal, resp. formally versal, at $s$ if,
for any Henselian, resp. 
Artin, local $\sO_B$-algebra $A$ with residue field $k(s)$ and any deformation
$Y\to\Spec A$ of $\sX_s$, there is an $\sO_B$-homomorphism $g: R\to A$ and
an $A$-isomorphism $Y\to \sX\otimes_R A$.
Artin proved that in these circumstances, the formal versality of $f$ at
$s$ implies its versality at $s$ and that versality is a Zariski open
condition on $S$.

\def \hO{\widehat{\sO}}

Let $\hO$, resp. $\hR$, denote the $b$-adic completion of $\sO_{B,b}$, 
resp. the $s$-adic
completion of $R$. There are compatible homomorphisms $W(k(b))\to \hO$
and $W(k(s))\to \hR$, which are natural if $k(s)$ is perfect. 
Put $\Lambda =\hO{\hat{\otimes}}_{W(k(b))}W(k(s))$,
so that $W(k(s))\to \hR$ factors through $\Lambda$. Note that $\Lambda$
is a complete local ring with residue field $k(s)$. Let $\sC_{\Lambda}$ denote
the category of Artin local $\Lambda$-algebras with residue field $k(s)$
and $F$ the set-valued deformation functor of $\sX_s$ defined on
$\sC_{\Lambda}$. There is a hull for $F$, say $D_{\sX_s}$, and
a classifying $\Lambda$-algebra
homomorphism $\phi:D_{\sX_s}\to\hR$ that is unique on tangent spaces.

\begin{lemma} \label{versal-versus-hull}
Suppose that $k(s)$ is algebraic over $k(b)$.
Then $f$ is formally versal at $s$ if and only if $\phi$ is formally
smooth.
\begin{proof} Assume that $f$ is formally versal at $s$.
Suppose that $A\in\sC_{\Lambda}$ and that $\xi\in F(A)$.
Then there is an $\hO$-homomorphism $\pi:\hR\to A$ such that $\sX$
induces $\xi$. Let $i:\Lambda\to \hR$ and $j:\Lambda\to A$ be 
the structural maps.
It is enough to know that $\pi\circ i=j$, that is, that $\pi$ is
a $\Lambda$-homomorphism. For this, argue by induction: pick a small
principal surjection $\sigma:A\to A_1$ with kernel $I$, put 
$\sigma\circ\pi =\pi_1$
and assume that $\pi_1\circ i = \sigma\circ j$. Then 
$\pi\circ i-j \in\Der_{\sO}(\Lambda, I)$, which vanishes since
$k(s)/k(b)$ is algebraic and separable.

The converse is immediate.
\end{proof}
\end{lemma}

\begin{definition}
With the notation as above, a family is \Definition{miniversal} at $s$
if $k(s)$ is algebraic over $k(b)$ and a classifying map 
$\phi:D_{\sX_s}\to\hR$ is an isomorphism. 
By a \Definition{small deformation} of a proper variety
$X$ we mean a proper family $f:\sX \to (S,s)$ such that $S$ is
local and $\sX_s\cong X$. 
\end{definition}

The main result
of [KM97] is that the geometric quotient $F$ of the stack $\sF$
exists, as an algebraic space of finite type over $B$,
if all $\Isom$ schemes are finite. Their idea is first to construct the
quotient locally in the \'etale topology, and then glue. The local part of
the argument goes as follows. Start with $f:\sX\to S$ that is everywhere
versal and everywhere surjective and consider the groupoid (that is,
the groupoid object in the category of schemes or algebraic spaces) 
$R=\Isom_{S\times_B S}(pr_1^*\sX,pr_2^*\sX)$, with source and target
maps $p,q:R\to S$, composition $c:R\times_{p,S,q}R\to R$ and identity map
$e:S\to R$. In particular, $e$ is a section of both $p$ and $q$. 
(We fix this notation.)
The key tautologies (due to Artin [Ar74a])
are that $R$ represents the set--valued functor $S\times_{f,\sF,f}S$ and 
that the versality of $f$ implies that both projections $S\times_{f,\sF,f}S\to S$,
and so the maps $p,q:R\to S$, are smooth.
Assume that $j=(p,q):R\to S\times_B S$
is finite. This is automatic if,
for example, $\sF$ is a stack of polarized non-ruled surfaces
with only RDPs, by the theorem of Matsusaka and Mumford. For $t\in S$,
put $P_t=p^{-1}(t)$, with projection $q_t:P_t\to S$. 
Define the orbit $O(t)$ to be the image of $q_t$. The essential local steps 
of [KM97] are
to slice $S$ by a subscheme $W$ (so that, in particular, $W$ is
regularly embedded in $S$) transverse to a given orbit in such
a way that the induced groupoid $R\vert_W$ is quasi-finite and flat
(in fact, locally complete intersection) over $W$ and then to show
that there is a connected component $P$ of $R\vert_W$ such that 
$P$ is  a subgroupoid that
contains the inverse image of the diagonal
and is finite over $W$. They then construct the quotient $W/P$
directly and check that it is, locally, the quotient $S/R$.
(They work in an apparently
more general context, where $p,q$ are only assumed to 
be flat, and then slice through a Cohen-Macaulay point, for example,
a generic point of each orbit.
When $p,q$ are Cohen-Macaulay morphisms,
the slice can be taken through any point. However, Artin showed [Ar74a]
that any flat groupoid is equivalent to a smooth one.)

In terms of an arbitrary algebraic stack $\sF/B$ with a smooth presentation
$S\to\sF$, the finiteness of $R\to S\times_B S$ is equivalent to the finiteness
of the representable diagonal morphism $\sF\to\sF\times_B\sF$.
However, we have been unable to state all of our results in this generality,
as explained below.

\begin{lemma}
\label{pseudotorsor}
Let $S$ be a scheme, $G$ a flat $S$--group scheme and $P\to S$
a pseudo--torsor under $G$. Then a geometric quotient
$P/G$ exists as an algebraic space
and the morphism $P\to S$ factors through a monomorphism
$P/G\to S$.

\begin{proof}
The existence of the quotient follows from
[Ar74a] Cor.~6.3, or from the fundamental theorem of Keel and Mori
[KM97]. It is then enough to show that for each 
geometric point $\bar s \to S$ the fibre of $P/G$ over $\bar s$ 
is either empty or isomorphic to $\bar s$. So we are reduced to 
the case where $S=\bar s$, where the result is obvious.
\end{proof}
\end{lemma}

\begin{lemma}
\label{subscheme}
Suppose that $I, S$ are algebraic spaces of
finite type over a Noetherian base, that $f:I\to S$
is a monomorphism and that $I$ is not empty.
Then there is an open subspace $U$ of $S$ such that $f^{-1}(U)$
is not empty and $f^{-1}(U)\to U$ is a closed embedding.

\begin{proof} Since $f$ is quasi--finite, there is, by Zariski's Main Theorem,
an open embedding $j:I\to J$ and a proper morphism $g:J\to S$
with $f=g\circ j$. After replacing $J$ by the closure of $j(I)$,
is necessary, we can assume that there is a closed subspace $Z$ of $J$,
not containing any irreducible component of $J$, such that
$j$ induces an isomorphism $I\to J-Z$. Since $\dim I=\dim J$
and $\dim Z < \dim J$, if we define $S_0=S-g(Z)$,
$J_0=J-g^{-1}(g(Z))$ and $I_0=I-f^{-1}(Z)$,
we see that $I_0$ is not empty, that $j$ induces an isomorphism
$I_0\to J_0$ and that $J_0\to S_0$ is proper. So $U=S_0$ will do.
\end{proof}
\end{lemma}
 
We shall use \ref{pseudotorsor} when $\sX\to S$ is proper and everywhere versal,
$t$ is a closed point in $S$, $P=P_t=\Isom_S(\sX_t\times S,\sX)$ and
$G=G_t\times S$, where $G_t=\Aut_{\sX_t}$. 
The identity of any group scheme will be denoted by $1$.

We shall use \ref{subscheme} to
show that the induced morphism $\pi_t:I_t\to S$ is an 
isomorphism to a subscheme,
so that $I_t$ is identified with the largest subscheme of $S$ over 
which all geometric
fibres are geometrically isomorphic to $\sX_t$, and this subscheme is
identified with $O(t)$.

Note that, 
via $q_t$, $P_t$ is a pseudo-torsor
under $G_t\times S$. By \ref{pseudotorsor} the quotient $I_t=P_t/G_t\times S$
exists and the natural morphism $\pi_t:I_t\to S$ is a monomorphism.

Now fix a closed point $s$ of $S$ and assume that $f:\sX\to S$ is everywhere versal,
but not necessarily surjective.

\begin{proposition} $\pi_s$ is an isomorphism to a subscheme of $S$.

\begin{proof} By \ref{subscheme} there is a closed point $u$ of $I_s$
such that $\pi_s(u)=v$, say, has a neighbourhood over which
$\pi_s$ is an isomorphism to a closed subscheme. Then $X_{\bar v}\cong X_{\bar s}$,
so that $P_s$ is identified with $P_v$ and $I_s$ with $I_v$.
Moreover, comparison of the germs $S_{\bar s}=(S,\bar s)$ and 
$S_{\bar v}=(S,\bar v)$ with hulls for $X_{\bar s}$
and $X_{\bar v}$ shows that the induced families over
$S_{\bar s}$ and $S_{\bar v}$ are smoothly equivalent.
Since the question of whether $\pi_s$ is an isomorphism to a subscheme
is unaffected by smooth base change and is local on $S$, we are done.
\end{proof}
\end{proposition}

\begin{remark} We do not know whether $I_t\to S$ is an isomorphism to a subscheme
for every smooth groupoid. This is one reason why we do not try to state all
our results in terms of general algebraic stacks. Another reason, which is related,
is that we do
not have a natural definition of what it means for a smooth chart of an
algebraic stack to be miniversal at a point.
\end{remark}

Henceforth we identify $I_s$ with its image $O(s)=q_s(I_s)$.

Assume, as above, that $f:\sX\to S$ is a proper and flat
family, maybe with further structure, and assume
also that $f$ is versal everywhere.
 
\begin{definition} The \emph{excess} of $f$ at $s$, denoted
$\exc_s(f)$, is the difference
between $\dim_sS$ and the dimension of a hull for $\sX_s$.
\end{definition}

For any group scheme, $e$ will denote the identity.

\begin{theorem}\label{strata-dimension} \label{2.7}
Assume that $S$ is connected.

\part[i] The function
$t\mapsto \dim T_e(G_t)+\exc_t(f)$ is constant on the closed points of $S$.

\part[ii] $P_t$ is everywhere smooth of dimension equal to
$\dim T_e(G_t)+\exc_t(f)$.

\part[iii] $I_t$ is everywhere smooth of dimension equal to
$\dim T_e(G_t)+\exc_t(f)-\dim G_t$.

\part[iv] If $f$ is miniversal at $s$, then $P_t$ is everywhere smooth 
of dimension equal to
$\dim T_e(G_s)$ and $I_t$ is everywhere smooth of dimension equal to
$\dim T_e(G_s) -\dim G_t$. 

\part[v] Suppose that $f$ is miniversal at $s$. Then $f$ is miniversal
everywhere if and only if the function $t\mapsto\dim T_e(G_t)$
is constant on $S$.

\begin{proof} As already remarked,
the versality of $f$ implies the smoothness of $p,q$.
So $P_t$ is smooth, and then $I_t$ is too.

Pick a closed point $v\in P_t$ and put $q_t(v)=w$. Then 
$\sX_{\bar t}\cong \sX_{\bar w}$.

Assume first that $f$ is miniversal at $w$.
Then the derivative of $q_t$ vanishes at $v$ and the derivative of
$\pi_t$ is injective. Hence the tangent space $T_v(P_t)$ is
isomorphic to $T_e(G_t)$, so that $P_t$ is smooth and 
$\dim_v(P_t) = \dim T_e(G_t)$.

In general, the morphisms $P_t\to I_t\to S$
are obtained from the corresponding morphisms where $S$ is replaced by
a miniversal deformation space $S_1$ of $\sX_w$ by a smooth base change
$S\to S_1$ of fibre dimension $\exc_w(f)$. Hence
$\dim_v P_t =\dim T_e(G_t)+\exc_w(f)$.

Now take $v=e(t)$. Then $w=t$, so that $\dim_{e(t)}P_t = \dim T_e(G_t) +\exc_t(f)$.

Since $S$ is connected, the fibre dimension of the smooth morphism
$p$ is constant along any section. So
$t\mapsto \dim T_e(G_t) +\exc_t(f)$ is constant.
Also $\dim_vP_t=\dim T_e(G_w) +\exc_w(f)$, since $\sX_{\bar t}\cong \sX_{\bar w}$,
and so $P_t$ is of constant dimension.

For (4) note that $\dim P_t = \dim T_e(G_t) +\exc_t(f)
=\dim T_e(G_s)+\exc_s(f)=\dim T_e(G_s)$, and similarly for $\dim I_t$.
Finally, (5) is an immediate consequence of (1).
\end{proof}
\end{theorem}

\begin{remark} \ref{strata-dimension} 
(5) is reminiscent of Schlessinger's result [Sch68], that 
the hull pro-represents the infinitesimal deformation functor 
if and only if $\Aut_{\sX/S}^0$ is smooth over $S$. However, the
comparison is not exact; it is true that, for example, the
flatness of $\Aut_{\sX/S}^0$ over $S$ implies the miniversality of
$f$ everywhere, but it is not clear that the converse is true.
\end{remark}

\begin{corollary}\label{2.8}
For every closed point $s\in S$ mapping to $b\in B$, 
the set of closed points $t$ such that
$\sX_t$ is geometrically isomorphic to $\sX_s$ is the support of a smooth
subscheme of $S\otimes k(b)$. If also $f$ is miniversal at $s$,
then the dimension of this subscheme is $\dim \Lie(G)-\dim G$,
where $G$ is the automorphism group scheme of $\sX_s$.
\begin{proof}
Immediate from \ref {strata-dimension}.
\end{proof}
\end{corollary}

\begin{corollary} \label{fibre-dimension}
Suppose that $\sF$ is as above and for every $T$-point $\sY\to T$
of $\sF$ the automorphism group scheme $\Aut_{\sY/T}$ is finite.
Then the fibres of the classifying
morphism $S\to F$, the geometric quotient of $\sF$,
are set-theoretically smooth. If $f$ is miniversal
at $s$ and $S$ is connected, then the fibres of
$S\to F$ are all of constant dimension $\dim T_e(\Aut_{X_s})$. 

\begin{proof} 
The fibres of $S\to F$ are, as sets, the incidence subschemes $I_t$,
and we are done by \ref{strata-dimension}.
\end{proof}
\end{corollary}

We close this section with a discussion of the local structure of a Keel-Mori stack
and the morphism to its geometric quotient in the context of the moduli of
projective varieties.
There is a well known theorem to the effect that if $X$ is a projective variety
with aa finite reduced groupscheme $\Aut_X$ of automorphisms, then $\Aut_X$ acts on the
miniversal deformation space $\Def_X$ of $X$ and the quotient $\Def_X/\Aut_X$ is
the formal germ of the coarse moduli space at the point corresponding to $X$.
In this section we give a proof of this in the algebraic (i.e., henselian) context; 
we could not locate one in the
literature. (In particular, any apparent similarity between this and Th\`eor\'eme 6.1
of \c[LMB00] is only apparent; the group $G$ of \emph{loc. cit.} is not
intrinsically related to $X$.) We also discuss how to extend this to 
cover what happens when $\Aut_X$ is more general.

Assume that we have a stack $\sF$ of projective varieties such that 
for all spaces $S$ and $T$ both projections of the fibre product $S\times_{\sF}\sT$ 
to $S$ and $T$ are finite. That is, $\Isom$ schemes are finite.
We want a description, local in the \'etale topology on the geometric quotient
$F$ of $\sF$, of $S\to F$,
where $f:X\to S$ is everywhere versal, which is more informative than the
tautological statement that it is the quotient by a smooth groupoid. However, we can
only do this under the following hypothesis:
\medskip

\noindent (*)\quad there is a closed flat normal subgroup scheme
$H/S$ of $\Aut_{\sX/S}$ such
that $\Aut_{\sX/S}/H$ is unramified over $S$. 
\medskip

Equivalently, there is a closed flat normal subrelative groupscheme $\sH$ of the stabilizer
relative groupscheme $\sG=\sG_{\sF/B}$ such that $\sG/\sH\to\sF$ is unramified.


The next definition helps in the comparison of different groupoids.

Suppose that $R\stackrel{p,q}\rightrightarrows S$ and 
$\tR\stackrel{\tp,\tq}\rightrightarrows\tS$ are groupoids, with compositions
$c,\tc$ respectively, and that $a:\tS\to S$
is a morphism. Then a morphism $b:\tR\to R$ is an \emph{equivariant morphism
of groupoids over} $a:\tS\to S$ if the squares
$$\begin{equation*}
\xymatrix{
{R} \ar[r]^p & {S} \\
{\tR} \ar[r]_{\tp} \ar[u]_b & {\tS} \ar[u]_a
}
\end{equation*}
$$
and
$$\begin{equation*}
\xymatrix{
{R} \ar[r]^q & {S} \\
{\tR} \ar[u]_b \ar[r]_{\tq} & {\tS} \ar[u]_a
}
\end{equation*}$$
are Cartesian and there is a commutative diagram
$$\begin{equation*}
\xymatrix{
{\tR\times_{\tp,\tS,\tq}\tR} \ar[r] \ar[d]_{\tm} & 
{(R\times_{p,S,q}R)\times_{p\circ pr_2,S}\tS} \ar[d]^{m\times 1_{\tS}}\\
{\tR} \ar[r] & {R\times_{p,S}\tS}
}
\end{equation*}$$
where the horizontal maps are the isomorphisms given by the preceding
Cartesian squares.

\begin{lemma} Suppose that $b:\tR\to R$ is an equivariant morphism
of groupoids over $a:\tS\to S$.

\part[i] If $j=(p,q):R\to S\times S$ is finite, then so is
${\tilde{\jmath}}:\tR\to\tS\times\tS$.

\part[ii] If $P$ is any of the properties flat, smooth, finite, \'etale,
locally complete intersection (l.c.i.) and $R$ has property $P$,
then so does $\tR$.

\part[iii] Suppose that $W\to S$ is a morphism and that
$\tW=W\times_S\tS$. Then the induced groupoid

\begin{proof} The first two parts are trivial. For the third,
the groupoid $R\vert_W$ is defined by the diagram
$$
\xymatrix{
{R\vert_W} \ar[r]^{p_W} \ar[d]_{q_W} & {R\times_{p,S}W} \ar[r] \ar[d] & {W}\ar[d]\\
{W\times_{S,q}R} \ar[r] \ar[d] & {R} \ar[r]^p \ar[d]_q & {S}\\
{W} \ar[r] & {S}
}
$$
where all squares are Cartesian.
Then the analogous diagram defining $\tR\vert_{\tW}$
lies over this in an obvious way, so that the resulting array (consisting of three
cubes) is commutative and every square is Cartesian.
\end{proof}
\end{lemma}

Fix a closed point $s\in S$ and assume that $f$ is miniversal at $s$.

\begin{theorem}\label{local structure} Assume hypothesis (*). 
Then $S\to F$ is, locally on $F$, isomorphic to the quotient of a smooth
morphism $\tW\to W$, where $W\inj S$ is a slice,
by an equivariant action of the unique finite \'etale group scheme
$\tG$ over $B$ that extends $G:=(\Aut_{\sX/S}/H)_s$. 

\begin{proof} First, localize $B$ so that it is a local scheme and
$k(s)=k(b)$. Put $\Aut_{\sX/S}=G$.

Now take an \'etale local slice $W$ in $S$ through $s$, as in [KM97]. 
Put $q_1^{-1}(W)=\tW\subset\tS$; this is a slice in $\tS$ through any point
lying over $s$. According to [KM97], the natural map $W/R_1\vert_W\to S/R_1$
is \'etale, so the same is true of $\tW/\tR_1\vert_{\tW}\to\tS/\tR_1$.
So it remains to describe the morphism $\tW/\tR_1\vert_{\tW}\to W/R_1$.

Since $R_1\to S\times S$ is unramified, the morphisms $p_{1,W},q_{1,W}$ are \'etale.
So $p_{1,W}^{-1}(s)=\coprod \Spec k_i$, with $k_i/k(s)$ finite and separable.
Since $W$ is local, it follows that 
$p_{1,W}^{-1}(s)=j_{1,W}^{-1}(s,s)\subset j_{1,W}^{-1}(\Delta_W)$.
There is [KM97] a decomposition $R_1\vert_W=P\coprod Q$, where $P$ is a subgroupoid
finite over $W$ and $P\supset j_{1,W}^{-1}(\Delta_W)$. Then $p_{1,W}^{-1}(s)\subset P$,
and so, since $W$ is local, $R_1\vert_W=P$. That is, $R_1\vert_W$ is a finite 
\'etale groupoid over $W$.
Note also that $j_{1,W}^{-1}(s,s)\cong G$.

The existence and uniqueness of $\tG$ follows from 
the fact that the \'etale covers of $\Spec k(s)$
and of $B$ form equivalent categories.
Moreover, the equivalence of categories of finite \'etale covers
give isomorphisms $\psi:R_1\vert_W \to W\times_B\tG$ and
$\tpsi:\tR_1\vert_{\tW}\to \tW\times_B\tG$ such that
$pr_1\circ \psi = p_{1,W}$ and $pr_1\circ\tpsi =p_{1,\tW}$. 
Then $q_{1,W}\circ(\psi^{-1})$,
respectively $q_{1,\tW}\circ(\tpsi^{-1})$, defines an action
of $\tG$ on $W$, respectively $\tW$ (the groupoid axioms carry over 
to the axioms describing
the action of a group scheme) and these actions are equivariant
with respect to $\tW\to W$. Finally, $W/(R_1\vert_W)\cong\Spec \sO_{F,x}^h$,
and we are done.
\end{proof}
\end{theorem}

There are two special cases of this that we make explicit.
One, when $\Aut_{X_s}$ is reduced, is the following well known folk theorem 
(see, for example,
[FC]) which lacked, apparently, a published proof except when
$k(s)$ is separably closed ([LMB00] 6.2.1). Note that, despite its
apparent similarity, this result is not equivalent to Th. 6.2 of \c[LMB00].
For one thing, the group $G$ appearing there is constructed in the course of the
proof as a symmetric group $\frak S_d$, where $d$ is not intrinsic to the 
context.

\begin{theorem}\label{local}
If $B$ is henselian, $k(b)=k(s)$, $G=\Aut_{X_s}$ is reduced
and $f:\sX\to S$ is miniversal at $s$, then $F$ is locally isomorphic 
at $\pi(s)$ to the quotient of $S$ by an action of $\tG$. 
\begin{proof} Take $H=1$ above.
\end{proof}
\end{theorem}

\begin{remark} This combines with Th. 6.2 of \c[LMB00] to permit the
definition and construction of the {\emph{henselization}} of a 
Deligne--Mumford stack $\sF$ at a point $P:T=\Sp K\to \sF$.
Consider the inverse system of $2$-Cartesian diagrams
$$\xymatrix{
{T} \ar[r] \ar[d]^{=} & {\sG} \ar[d]_{\phi} \\
{T} \ar[r]^{P} & {\sF}
}$$
where $\phi$ is \'etale and representable, and we demand also that the
natural morphism $\sG\to\sG\times_{\sF}\sG$ be an isomorphism.
The inverse limit of this system is the stack $[S/\tG]$ of \ref{local}
and is the henselization of $\sF$ at $P$.
\end{remark}

\begin{theorem} If $f:\sX\to S$ is everywhere versal and $\Aut_{\sX/S}$ is flat over
$S$, then $\pi:S\to F$ is smooth.
\begin{proof} Take $H=\Aut_{\sX/S}$ above.
\end{proof}
\end{theorem}
\end{section}


\begin{section}{Tangent sheaves and foliations on singular varieties}
\label{singular foliations}
We will need to generalise the notion of a smooth (height 1) foliation from the
case of a smooth variety discussed in [Ek86] to the case of singular varieties.
\begin{definition}
Suppose that $f\co X\to S$ is of finite type. A \Definition{smooth
$1$-foliation} on $X/S$ is a subsheaf $\cF$ of the tangent sheaf $T_{X/S}$ such
that for every closed point $x$ of $X$ the induced map 
$\cF \Tensor_{\cO_X} k(x) \to \Hom_{k(x)}(m_x/m_{f(x)}+m^2_x,k(x))$ 
is injective and such that $\cF$ is
closed under commutators (and, when $S$ has positive characteristic $p$,
under $p$'th powers).
\begin{remark}
When $X/S$ is smooth, the injectivity condition simply says that $\cF$ is a
subbundle of $T_X$ so this corresponds to the already established notion.
\end{remark}
\end{definition}
We shall see that, just as in the smooth case, a smooth $1$-foliation gives
rise to a flat infinitesimal equivalence relation on $X$.
We shall need the following lemma.
\begin{lemma}
Let $(A,m,k)$ be a local ring, $M$ a finitely generated
$A$-module and $N$ a finitely generated $A$-module provided with 
an $A$-homomorphism
$N \to M^* := \Hom_A(M,A)$. Assume that the map $N/mN \to
\Hom_k(M/mM,k)$ induced by this 
is injective. Then $N$ is free and a direct factor of
$M^*$. Furthermore, $M$ can be written as a direct sum $M_1 \bigoplus M_2$ such
that $M_1$ is free and $M_1^* \hookrightarrow M^*$ can be identified with 
$N\hookrightarrow M^*$.
\begin{proof}
Let $n := \dim_k N/mN$ and pick $m_1,\dots,m_n \in M$ such that the composite
of $N/mN \to \Hom_k(M/mM,k)$ and the map $\Hom_k(M/mM,k) \to k^n$ given by
evaluation at the residues $\overline{m_i} \in M/mM$ is an
isomorphism. Evaluation at the $m_i$ gives a map $M^* \to A^n$ and the
composite $N \to M^* \to A^n$ induces an isomorphism upon reduction modulo
$m$. By Nakayama's lemma it is then surjective and by the projectivity of $A^n$
and Nakayama's lemma again it is an isomorphism. This shows that $N$ is a direct 
factor of $M^*$ as well as isomorphic to $A^n$. By construction the basis of $N$ 
thus constructed is dual to $\{m_i\}$ and hence the $m_i$ generate a direct
summand of $M$.
\end{proof}
\end{lemma}
\begin{proposition}\label{foliation}
Suppose that $X\to S$ is a finite type morphism, that $S$ has 
positive characteristic $p$
and $\cF \subseteq T_{X/S}$ a smooth $1$-foliation. Then $\cF$ is locally free and
locally a direct summand of $T_X$. Furthermore, if \map fXY is the ``scheme of
constants'' of $\cF$, i.e., $\sO_Y=\sO_X^{\cF}$,
then $f$ is a flat map of degree $p^n$, where
$n:=\rk \cF$, $\cF=T_{X/Y}$ and $\cF^*=\Omega^1_{X/Y}$. Finally, there are
locally sections $f_1, f_2,\dots,f_n \in \cO_X$ such that they generate $\cO_X$
as $\cO_Y$-algebra. Put $f_i^p=g_i$; then $X$ is recovered
from $Y$ from the formula
$\sO_X=\sO_Y[z_1,\ldots,z_n]/(z_1^p-g_i,\ldots,z_n^p-g_n)$.
\begin{proof}
If we apply the lemma with $M=\Omega^1_{X/S}$
we immediately get that $\cF$ is locally
free and locally a direct summand of $T_{X/S}$. 
Let now $\cA$ be the subalgebra of
the ring of differential operators relative to $S$
of order $<p$ generated by $\cO_X$ and
$\cF$. If we filter $\cA$ by the order filtration then there is a map
$\Symm^{<p}\cF \to \gr \cA $. We also have an obvious map $\gr \cA \to \gr
\Diff^{<p}$ and an evaluation map $\Diff^{<p} \to
\lHom_{\cO}(\Symm^{<p}\Omega^1,\cO)$ and finally a restriction map
$\lHom_{\cO}(\Symm^{<p}\Omega^1,\cO) \to \lHom_{\cO}(\Symm^{<p}\cF^*,\cO)$. The
composite of these maps give a map $\Symm^{<p}\cF \to (\Symm^{<p}\cF^*)^*$. This
map is easily seen to be the standard map which is an isomorphism as we are only
considering degrees $<p$. This shows that $\gr \cA= \Symm^{<p}\cF$, that $\cA$
is locally free and that its dual defines a closed subscheme of $X\times_S X$
which is a finite
flat equivalence relation. Its quotient is exactly $Y$ and therefore
the map $X \to Y$ is flat of degree $p^n$. If $x \in X$ and $y:=f(x)$ we now
have an exact sequence
\begin{displaymath}
m_y/m_y^2 \to m_x/m_x^2 \to \Omega_{X/Y}^1 \Tensor k(x) \to 0
\end{displaymath}
which shows that if we let $f_1,\dots,f_n \in \cO_{X,x}$ be such that their
residues form a basis for $\Omega_{X/Y}^1 \Tensor k(x)$ then they generate
$\cO_{X,x}$ as $\cO_{Y,y}$ module and further the monomials where each $f_i$
occurs to power $<p$ form a set of module generators which have cardinality
$\le p^n$ and hence form a basis.
\end{proof}
\end{proposition}

\begin{corollary}\label{smoothness} If $X$ is a normal $k$-variety and
$T_X$ is a smooth $1$-foliation, then $X$ is smooth.
\begin{proof} 
Since $X$ is normal, the quotient morphism $X\to X/T_{X/S}$ is the relative
Frobenius. This is flat, by the previous result, and now Kunz' criterion shows
that $X$ is a regular scheme, so smooth over $S$.
\end{proof}
\end{corollary}

\begin{corollary} \label{vanishing} Suppose that $X$ is a normal $k$-variety and that
$V\subset H^0(X,T_X)$ is a finite-dimensional vector space that generates
$T_X$. Then for every singular point $P$ of $X$
there is a non-zero element $v$ of $V$ that vanishes at $P$. Moreover, if $\ch k=p>0$
and $V$ is a $p$-Lie algebra,
then we can take $v$ to be $p$-closed. 
\begin{proof} We can assume that $\dim X\ge 2$.

If $v(P)$ is never zero, then $T_X$ is a smooth $1$-foliation near
$P$. If $\ch k=p>0$, then the set $W=\{v\in V\vert v(P)=0\}$ is a non-zero
sub-$p$-Lie algebra of $V$, and any such contains a $p$-closed line.
\end{proof}
\end{corollary}

\begin{corollary}\label{foliation base change}
Given $X\to S$ and $\sF$ as above, the formation of
the quotient $X/\sF$ commutes with base change.
\noproof
\end{corollary}
\begin{corollary}
\part Let $X \to Y$ be the quotient by
a smooth $1$-foliation, where $X$ is assumed to be normal. Then
$\Omega^1_{X/Y}$ is locally free. Thus $T_{X/Y}$ is locally a direct factor of
$T_X$ and the quotient is reflexive. In particular, if it is of rank 1 it is a
line bundle.

\part\label{tangloctriv} Suppose $X$ is a normal surface and $D$ a global vector
field on it which generates a smooth $1$-foliation (i.e., $\cO_XD
\hookrightarrow T_X$ is a smooth $1$-foliation). Then $T_X/\cO_XD$ is a line
bundle and in particular $T_X$ is locally free.

\part\label{tangtriv} Let $X$ be a $\Z/2$- or $\alpha_2$-Enriques surface such that the
canonical double cover $Y$ has isolated singularities. Then the tangent bundle
$T_Y$ is trivial.
\begin{proof}
The first two parts  follows directly from the proposition. For the third we get 
from the second and the fact that $T_X/\cO_XD$ is the inverse image of the dual of
$\Omega^1/\cO_X\omega$ where $\omega$ is a non-zero form combined with the fact
that by assumption the zeros of $\omega$ are isolated that $T_Y$ is an
extension of $\cO_Y$ by itself. As $H^1(Y,\cO_Y)=0$ this extension is trivial.
\end{proof}
\end{corollary}
A smooth $1$-foliation on a smooth variety has a smooth quotient simply because
the quotient map is flat. This leads to the impression that in general the
singularities of a quotient should be simpler than that of the base variety. We
do not have a general statement to that effect but the following lemma may be
seen as further evidence of its existence.
\begin{proposition}\label{RDP-quot}
Let $R$ be a complete $2$-dimensional local $k$--algebra of
characteristic $p$ with an RDP, $\cF$ a smooth
$1$-foliation on $R$ and $S$ its ring of constants. 
Then the list of possibilities is as follows.

\part[1] $S$ is regular and $R$ has a Zariski RDP.

\part[2] \label{even} $p=2$,
$S$ has an $A_n$-singularity and 
$R$ is either of type $A_{2n+1}$ or, 
when $n=1$, $D^0_{2m+1}$ or,
when $n=2$, $E^0_6$.

\part[3] \label{odd} $p\ge 3$, $R$ is of type $A_{pr-1}$ and
$S$ is of type $A_{r-1}$.

In particular, $R$ has a Zariski RDP if and only if $S$ is regular.
\begin{proof}
$R$ has embedding dimension $3$ and by \ref{foliation} 
$R=S\bigl[[v]\bigr]/(v^p-g)$ where $g \in m_S$. So $S$ has either
embedding dimension two, in which case $S$ is regular and $R$
has a Zariski RDP,
or it has embedding dimension $3$
and then $0\ne \bar g \in m_S/m_S^2$. We may assume the latter.
Thus we can assume that $S=k[[x,y,z]]/(f)$ and $g=z$.
Then $R=k[[x,y,v]]/(F)$, where $F=f(x,y,v^p)$. Then we can write
\begin{displaymath}
f=a_2(x,y)+(x,y)^3+z^mu+z^nh(x,y,z)+z(x,y)^2,
\end{displaymath}
where $u\in k[[{z}]]^*, h=\sum h_i(x,y)z^i$, $h_i$ is 
homogeneous and linear in $x,y$
and $h_0\ne 0$.

Now suppose that $p\ge 3$. 
RDPs are characterized as the double points that are absolutely
isolated, or as the double points whose blow-up at the origin
has only RDPs. Then the quadratic term of the strict transform
of $F$ under this blow-up is $a_2(x_1,y_1)$, where $x_1=x/v$
and $y_1=y/v$, and its cubic term is zero. Hence $a_2$ has rank $2$,
so that $R,S$ are of type $A_M,A_N$, say. 
Moreover, by analyzing what happens under a blow-up, the
form of $f$ reveals that
$N=\min\{m-1,2n-1\}$. Applying this analysis to $F$ shows that
$M=\min\{pm-1,2pn-1\}$, which proves \ref{odd}.

Now suppose that $p=2$.
Since $R$ has multiplicity $2$ 
$f \in m^2\setminus m^3$, where $m:=(x,y,z)$. We now divide the
discussion into three cases according to the rank of $\bar f \in m^2/m^3$.

Assume first that $\bar f$ has rank $3$. Up to linear transformations of $m/m^2$
there are two possible configurations for $\bar g$ and $\bar f$, seen as curves
in $\P(m/m^2)$ they intersect in either one or two points. In the first case, we
may, after a formal change of coordinates, assume that
$f=z^2+zj_2(x,y)+h_2(x,y)$, where the index indicates that $h_2,j_2 \in m^2$, and
$g=x$. By assumption $\overline{h_2} \in m^2/m^3$ is the product of two distinct
linear factors and we may assume that the second is $y$. After further formal
change of coordinates we may assume that $h_2=xy+k_3(y)$ and $j_2=j_2(y)$. Then
$R$ has the form $k[[{t,y,z}]]/(z^2+zj_2+x^2y+k_3)$. This is a
$D^0_{2m+1}$-singularity (when
it has an isolated singularity). In the second
case may assume that $f=z^2+xy$ and $g=z$ which leads to an $A_3$-singularity.

Assume now that $\bar f$ has rank $2$ so that it is a product of two distinct
linear factors. When $\bar g$ does not lie in the space spanned by them we may
assume that $g=z$ and that $f=z^{n+1}+xy$ and then $S$ has an $A_n$-singularity
and $R$ an $A_{2n+1}$-singularity. If it does lie in the space spanned by these
two factors then it cannot divide $\bar f$ because if it did $R$ would have
multiplicity three which it doesn't. This leads to the form $f=z^2+yz+h_3(x,z)$
and $g=y$. This gives that $R$ has equation $z^2+y^2z+h_3(x,z)$. In order for
this to be an RDP we need that $x^3$ occur in $\overline{h_3}$. This leads to a
singularity of type $E^0_6$ for $R$ while $S$ has type $A_2$.

Assume finally that $\bar f$ has rank $1$ so that it is a square of a linear
form. If $\bar g$ divides $\bar f$ $R$ would have multiplicity $4$ which is not
possible. We may therefore assume that $f=z^2+h_3(x,y,z)$ and $g=x$ so
that $R$ is defined by $z^2+h_3(x^2,y,z)$. From [Li69] p.~268 it follows
that this does not lead to an RDP.
\end{proof}
\end{proposition}
\begin{definition-lemma}\label{non-free subsheaf}
Let $X$ be a surface over a perfect field of characteristic $p$, $x \in X$, $f
\in \sO_x$ and $Y:= \Sp\sO_x[z]/(z^p-f)$ an isolated Zariski singularity. Then
there is a unique subsheaf $\sG \subset T_Y$ of colength $1$ such that a rank
$1$ subsheaf $\sF \subset T_Y$ closed under $p$'th powers is a smooth
$1$-foliation if and only if $\sF \not\subseteq \sG$. We will call this subsheaf 
the \Definition{sheaf of non-free vector fields}.
\begin{proof}
This is just the statement that the kernel of the natural homomorphism
$T_Y\Tensor k(y) \to \Hom(m_y/m^2_y,k(y))$ 
has a $1$-dimensional kernel, where $y \in Y$ is the
point above $x$. Using the fact that a basis for $T_Y$ is given by
$\partial/\partial z$ and $f_y\partial/\partial x-f_x\partial/\partial y$ and
$f_x,f_y \in m_x$ this follows immediately.
\end{proof}
\end{definition-lemma}

An isolated surface singularity in characteristic $p>0$ will be called a
\Definition{Zariski singularity} if it has the local form $R[z]/(z^p-f)$, where
$R$ is a regular local ring.

In characteristic $2$, a \Definition{Zariski RDP} is an RDP which is also a
Zariski singularity.  In Artin's notation [Ar77] these are the ones of type
$A_1$, $D_{2n}^{0}$, $E_7^{0}$ and $E_8^{0}$.  The underlying Dynkin diagrams
are the simply laced diagrams whose root lattice is $2$-elementary. Note
that contrary to the case of characteristic $0$ there are in general several
RDPs with the same Dynkin diagram, though only one Zariski RDP.
\begin{remark}
The terminology is motivated by the Zariski surfaces, surfaces in characteristic
$p$ birational to surfaces with equation $z^p-f(x,y)$.
\end{remark}
\begin{definition-lemma}\label{Zardef}
Let $f \in k[[{x,y}]]$, where $k$ is a field of characteristic $2$, be such that
$(f'_x,f'_y)$ is of finite codimension. Suppose $X \to \Spf R$ is a deformation
of the singularity $X_0=\Spec k[[{x,y,z}]]/(z^2-f)$ and that
$X/R$ admits a $2$-closed smooth
$1$-foliation $\sF$. Then this deformation is isomorphic to one of the type
$R[[{x,y,z}]]/(z^2-F)$. When $X \to \Spf T$ is a miniversal deformation we
will call the locus of such singularities the \Definition{Zariski locus}. It is
smooth of dimension equal to $\dim k[[{x,y}]]/(f'_x,f'_y)$. A deformation of a
Zariski singularity of this type will be called a \Definition{Zariski
deformation}.
\begin{proof} By \ref{foliation base change} $X/\sF=Y$, say, is formally
smooth over $\Spf R$. Let $Z\in\sO_X$ be a lifting of $z$. Then 
$Z^2:=F\in\sO_Y$ and we get a map $\sO_Y[Z]/(Z^2-F)\to\sO_X$. As this
is an isomorphism over the closed point of $\Spf R$ and $\sO_X$ is
projective of rank 2 over $\sO_Y$ this is an isomorphism.
\end{proof}
\end{definition-lemma}
\begin{lemma}\label{Zariski recognition}
Suppose that the henselian germ $(X,0)$ is an RDP of index $n$ and characteristic $2$
admitting an equicharacteristic
deformation ${\frak X} \to B$ such that $B$ is reduced and for every
geometric generic point $\bar\eta$ of $B$ the fibre
${\frak X}_{\bar\eta}$ has only Zariski RDPs of total index $n$. 
Then $(X,0)$ 
is a Zariski RDP and ${\frak X} \to B$ is a Zariski deformation.
\begin{proof}The existence of ${\frak X} \to B$ shows that
in an algebraic representative ${\frak X}_1\to V$ of a miniversal
deformation of $(X,0)$, there is a closed point $v\in V$
such that ${\frak X}_{1,v}$ has $n$ nodes. In characteristic $2$ a Zariski RDP
of index $r$ has 
$2r$ local moduli, so that, by the openness of versality and the
formal smoothness of the natural map from the miniversal deformation
space of a normal affine surface to the product of the miniversal deformation
spaces of neighbourhoods of its singularities,
$\dim V\ge 2n$. An inspection of Artin's tables [Ar77]
shows that $(X,0)$ is Zariski, and $\dim V=2n$.

In particular, for every geometric $b\in B$, the module
$T^1_{\frak X_b}$ has length $2n$. So $T^1_{\frak X/B}$ is a
locally free $\sO_B$--module, and then \ref{Tflat} $T_{\frak X/B}$ is
$B$--flat and its formation commutes with base change. So
$T_{\frak X/B}$ is $\sO_{\frak X}$--free. We conclude by \ref{Zardef}.
\end{proof}
\end{lemma}
Suppose $X$ is a local complete intersection surface with isolated singularities
and \map \pi\tX X a minimal resolution of singularities. The \Definition{local
Chern classes} $c^l_i(X) \in A_*(Z)$, where $Z := \pi^{-1}(X^{sing})$, are
defined as $c_i(T^\cdot_{\tX/X})\cap[\tX]$, where $T^\cdot$ is the relative
tangent complex (the derived dual of the cotangent complex) $T^\cdot_{\tX/X}$
considered as a perfect complex on $\tX$ with support in $Z$ and the $c_i$ are
the localized Chern classes in $A^*(Z \hookrightarrow X)$ of [Fu84], Example
18.1.3. If $x \in X^{sing}$ we define $c_i(X)_x \in A_*(\pi^{-1}(x))$ to
be the component in $A_*(\pi^{-1}(x))$ of $c^l_i(X)$. Also we define the
\Definition{local Chern numbers} $c^2_1(X)_x$ and $c_2(X)_x$ as
$\deg(c_1(T^\cdot_{\tX/X})\cap c_1(X)_x)$ and $\deg(c_2(X)_x)$ respectively.
\begin{proposition}
Let $X$ be a local complete intersection surface with isolated singularities,
\map\pi\tX X its minimal resolution and $i\co Z := \pi^{-1}(X^{sing})\hookrightarrow\tX$.

\part[1] For every irreducible component $E$ of $Z$ we have that $c^l_i(X)\cap
i^*\cO(E)=i^*c_i(T^\cdot_{\tX/X})$. In particular, $2g(E)-2=E^2+\deg(c^l_1(X)\cap
i^*\cO(E))$ and this determines completely $c^l_1(X) \in A_0(Z)$.

\part[2] (Local Noether's formula) For each $x \in X^{sing}$ we have that
\begin{displaymath}
-\lth_x(R^1\pi_*\cO_{\tX})=\frac{c^2_1(X)_x+c_2(X)_x}{12}.
\end{displaymath}

\part[3]\label{chernlocglob} We have that $c_i(\tX)=\pi^*c_i(X)+i_*c^l_i(X)$
and, when $X$ is proper, $c^2_1(\tX)=c^2_1(X)+\sum_{x \in X^{sing}}c^2_1(X)_x$
as well as $c_2(\tX)=c_2(X)+\sum_{x \in X^{sing}}c_2(X)_x$.
\begin{proof}
We have that $i^*\cO(E)\cap c^l_i(X)=i^*\cO(E)\cap c_i(T^\cdot_{\tX/X})\cap[\tX]$
and this equals $c_i(T^\cdot_{\tX/X})\cap\cO(E)\cap[\tX]=c_i(T^\cdot_{\tX/X})\cap i_*[E]$ by
[Fu84], Prop.~17.3.2, 
and this in turn equals $c_i(Li^*T^\cdot_{\tX/X})\cap [E]$ 
by the definition of the localized Chern classes. However, localised Chern
classes with support equal to the whole space are the ordinary Chern classes by
the analogue for Chern classes of [Fu84], Prop.~18.1, which gives the first
part of \DHrefpart{1}. The second part follows from the adjunction formula and
the negative definiteness of the intersection matrix for $Z$.

Continuing with \DHrefpart{3} we have that
$c(T^\cdot_{\tX})=\pi^*c(T^\cdot_X)c(T^\cdot_{\tX/X})$ because of the usual distinguished
triangle. Furthermore, $c_i(T^\cdot_{\tX/X})$ for $i>0$ has support on $Z$ so have
zero intersection with $\pi^*c_j(T^\cdot_X)$ for $j>0$. This shows that
$c_i(T^\cdot_{\tX})=\pi^*c_i(T^\cdot_X)+c_i(T^\cdot_{\tX/X})$ and 
[Fu84], Prop.~17.3.2
shows that $c_i(T^\cdot_{\tX/X})=i_*c^l_i(X)$. The rest of \DHrefpart{3} follows
from this and the orthogonality of $\pi^*c_i(T^\cdot_X)$ and $i_*c^l_i(X)$.

Finally, for \DHrefpart{2} we may assume that $X$ is proper and that $x$ is its
only singular point. Then \DHrefpart{2} follows from \DHrefpart{3} and the
Riemann-Roch formula for $X$ and $\tX$.
\end{proof}
\end{proposition}
We will use this result for rational double points and for it it becomes
necessary to compute the local Chern classes for them.
\begin{corollary}\label{tangchern}
\part[1] A surface RDP singularity has trivial local Chern classes.

\part[2] Suppose that $X$ is a surface with only isolated singularities and
whose singularities have the local form $R[z]/(z^2-f)$ where $R$ is smooth. Then
$c_1(T_X)=c_1(T^\cdot_X)$ and 
$$c_2(T_X)\cap[X]=c_2(T^\cdot_X)\cap[X]-\sum_{x \in
X^{sing}}r_x[x],$$
where $r_x := \lth_x(T^1(X))$.
\begin{proof}
\DHrefpart{1}: That the first local Chern class is trivial follows from the fact
that the relative $c_1$ is trivial and that the second one is follows from the
local Noether's formula and the fact the degree gives an isomorphism $A_1(Z)
\riso \Z$ (using the notations of the proposition).

As for \DHrefpart{2} we have, by definition, that $T^1_X=H^1(T^\cdot(X))$ and
$H^0(T^\cdot_X)=T_X$ and the other cohomology sheaves are zero.
\end{proof}
\end{corollary}
\begin{lemma}\label{cherneven}
Let $X$ be a smooth surface in characteristic $2$ and $\cF \inj T_X$ a
$1$-foliation. Then $\cF\Tensor\omega_X$ is the square of a line bundle.
\begin{proof}
By restricting to the part of $X$ where $\cF$ is smooth (and using that the
singular locus of $\cF$ is of codimension $\ge 2$) we may assume that $\cF$ is smooth.
If $X \to Y$ is the quotient by $\cF$ we have a double inseparable cover
$Y^{(-1)} \to X$. It corresponds to an $\sL$-torsor for some line bundle $\sL$
(cf. [Ek88], Prop.~1.11). It therefore gives a vector bundle embedding
$\sL^2 \inj \Omega^1_X$ which is orthogonal to $\sF$. This gives the lemma.
\end{proof}
\end{lemma}
\begin{proposition}\label{nonfree}
Let \map\pi XS be a proper cohomologically flat morphism of noetherian schemes and $\sE$ a
coherent $\cO_X$-module which is flat over $S$ such that its fibres are locally free.

\part[i]\label{nonfreecod1} The set $N$ of points $s \in S$ where $\sE_s$ is
free is locally closed.

\part[ii] If $S$ is regular and $N$ is dense then its complement
has everywhere codimension $1$.

\part[iii] If $N=S$ then $\sE=\pi^*\pi_*\sE$ and $\pi_*\sE$ is a locally free
$\pi_*\cO_X$-module that commutes with base change.
\begin{proof}
It is clear that $\sE$ is a locally free $\cO_X$-module 
and hence its rank is locally constant on $X$. To prove \DHrefpart{i} by the
semi-continuity theorem we may assume that $h^0(\sE_s)$ equals
$h^0(\cO_{X_s}^{\rank{\sE}})$ and that $S$ is reduced. Then $\pi_*\sE$ is
locally free and commutes with base change so that the locus where $\sE$ is not
free is the image of the support of $\pi^*\pi_*\sE \to \sE$ which is closed as
$\pi$ is proper.

\DHrefpart{iii} is now standard and for \DHrefpart{ii} if it is not true we may
localize at a generic point of the complement of $N$ and then take a suitable
quotient so as to reduce to the case when $S$ is $2$-dimensional and $N$ contains the
complement of one closed point $s$ and we need to show that it contains $s$ as
well. Being flat over $S$ $\sE$ has depth $2$ at $s$ as $\cO_S$-module. This
implies that $\pi_*\sE$ is reflexive and hence locally free. The map
$\pi^*\pi_*\sE \to \sE$ is an isomorphism outside of the complement of $s$. As
both have depth $2$ at $s$ it is an isomorphism.
\end{proof}
\end{proposition}
\begin{proposition}\label{Tflat}
For a flat family \map \pi XS of local complete intersection surfaces with
isolated singularities we set $T^1_{X/S} := H^1(T^\cdot_{X/S})$, where
$T^\cdot_{X/S}$ is the tangent complex.

\part $T^1_{X/S}$ is a coherent $\cO_S$-module that commutes with base change.

\part $T^1_{X/S}$ is a locally free $\cO_S$-module if and only if $T_{X/S}$ is a 
flat $\cO_S$-module that commutes with base change.
\begin{proof}
This follows immediately from the fact that $T^\cdot_{X/S}$ commutes with base
change as $\pi$ is flat and has universally amplitude $[0,1]$ with
$H^0(T^\cdot_{X/S})=T_{X/S}$.
\end{proof}
\end{proposition}
\begin{lemma}\label{flatstrat}
Let \map \pi XS be a proper morphism of noetherian schemes and $\sF$ a
coherent $\cO_X$-sheaf that is flat as $\cO_S$-module. Then there is a
subscheme stratification, $\{S_e\}_e$, of $S$ such that for a morphism \map fTS
$f^*\sF$ is cohomologically flat with $\pi_*f^*\sF$ locally free of rank $e$ if
and only if $f$ factors through $S_e \inj S$. If $\sF$ is cohomologically flat
in degree $0$ the same result holds true for cohomological flatness in 
degree $1$.
\begin{proof}
Recall (cf. [EGAIII:2], \S 7) that there is a coherent $\cO_S$-module
$\sM$ that represents $\sB \mapsto \pi_*(\sF\Tensor\pi^*\sB)$, that it
commutes with base change and is locally free if and only if $\sF$ is
cohomologically flat. Furthermore, when locally free its rank equals that of
$\pi_*\sF$. The existence of the stratification then follows immediately from
[Mu66], p.~56. The case of cohomological flatness in degree $1$ is
completely similar.
\end{proof}
\end{lemma}
\begin{lemma}\label{smoothfix}
Suppose \map fXS is a morphism of schemes and $G$ is an $S$-group scheme acting
on $X$. Let $x \in X$ be a fixed point of $G$ such that $f$ is formally smooth
at $x$ and suppose that $H^1(G_s,(\Omega_{X/S})_x)=0$, where $s:=\pi(x)$. Then
the fixed point locus of $G$ on $X$ is formally smooth at $x$.
\begin{remark}
The only consequence of this that we will use, that if $G$ is linearly reductive 
then the fixed point locus is formally smooth, is of course well known.
\end{remark}
\begin{proof}
To prove this we may assume we have a small infinitesimal $S$-extension $T
\hookrightarrow T'$ and an $S$-map of $T$ to the fixed point locus of $G$ and we 
want to show that it is liftable to $T'$. Now, such a lifting is nothing but a
lifting to $X$ that is invariant under $G$. A lifting to $X$ always exist
because of formal smoothness and the set a liftings is a torsor under
$(\Omega_{X/S})_x$. This torsor is a $G_s$-equivariant torsor and hence the
obstruction to a $G$-fixed lifting is an element of $H^1(G_s,(\Omega_{X/S})_x)$
and by assumption that obstruction is trivial.
\end{proof}
\end{lemma}

Suppose now that $X$ is a surface with Zariski RDPs of total index $r$ and minimal
resolution $\pi\co \tX\to X$.
\begin{theorem}\label{chernchange}
The sheaf $T_X$ is locally free and $c_2(T_X)=c_2(\tX)-2r$.
\begin{proof}
The local freeness of $T_X$ is \ref{tangloctriv} and the Chern class formula is
\ref{tangchern} and \ref{chernlocglob} combined with the fact that the length of $T^1$ 
is twice the index.
\end{proof}
\end{theorem}
\begin{proposition}\label{4.4}
Suppose that $X$ is Zariski RDP--K3 and that $r=12$. Then $T_X$ has trivial Chern
classes and either $T_X$ is free or there is a decomposition
$$0\to\sO(A)\to\pi^*T_X\to\sI_Z\sO(-A)\to 0,$$
where $A$ is effective and non--zero. Moreover,
in the second case $T_X$ is $H$-unstable
for all ample $H\in\Pic X$.
\begin{proof}
The triviality of the Chern classes follows directly from Theorem
\ref{chernchange}. The Riemann-Roch theorem then shows that
$h^0(X,T_X)+h^2(X,T_X) \ge 4$ but as $T_X$ is of rank $2$ with trivial
determinant it is self-dual and as $\omega_X$ is trivial this gives that
$h^0(X,T_X)=h^2(X,T_X)$ and thus $h^0(X,T_X)\ge 2$. This gives a map $\cO^2_X
\to T_X$ and if its image has rank $2$ it must be an isomorphism by taking its
determinant. The case when the image has rank $1$ leads to the second part of
the proposition.
\end{proof}
\end{proposition}

The next result gives the key relationship between Enriques and K3 surfaces
over an algebraically closed field.

\begin{theorem}\label{key}
\part{i} Suppose that $X$ is Zariski RDP--K3, that $r=12$ and that
$T_X$ is free. Then for any $2$-closed vector field $\xi$ on $X$ that does not
vanish at a singular point, the quotient $Y=X/\xi$ is smooth and Enriques.

\part{ii} Suppose that $Y$ is a unipotent Enriques surface whose canonical 
double cover $X\to Y$ is RDP-K3. Then $X$ is Zariski RDP of total index $12$
and $T_X$ is free.
\begin{proof} For \DHrefpart{i} the smoothness of $Y$ is an immediate consequence
of the non-vanishing of $\xi$. Then the adjunction formula applied to
the quotient map $\rho:X\to Y$ shows that $K_Y$ is numerically trivial.
Since $\rho$ is an \'etale homeomorphism, it follows that $b_1(Y)=0$
and $e(Y)=12$. Then $Y$ is Enriques, by the Bombieri--Mumford classification
of surfaces.

For \DHrefpart{ii} the only thing to prove is that $T_X$ has no subsheaf
$\sO(A)$ with $A>0$. Regard $\rho:X\to Y$ as the quotient of $X$ by a rank $1$
foliation $\sF\inj T_X$; the adjunction formula for $\rho$ shows that
$c_1(\sF)$ is numerically trivial, and comparison of the saturated
subsheaf $\sF\inj T_X$
with the putative subsheaf $\sO(A)\inj T_X$ gives a contradiction.
\end{proof}
\end{theorem}
\end{section}

\begin{section}{Further deformation theory: Enriques surfaces}
\label{Deforming Enriques surfaces}
A smooth, proper and polarized
surface $X$ is \Definition{correctly obstructed}
if $h^0(T_X)=0$, there is an effective formal deformation
over its hull $D_X$ and
$\dim D_X = h^1(T_X)-h^2(T_X)$.
\begin{theorem}
Suppose that $X$ is a smooth projective surface over
$k$, that $f:\sX\to S$ is versal and is miniversal at the closed point $s$. 
Assume that $\Aut_X$ is finite and that $S$ is connected.

\part[i] Every irreducible
component of $S$ whose geometric generic fibre is correctly
obstructed has dimension $h^1(X,T_X)-h^2(X,T_X)$.

\part[ii]\label{miniversals} If the set of geometric points in $S$ with
correctly obstructed fibres is dense,
then $(S,s)$ a local complete intersection over $\W$ of dimension
$h^1(X,T_X)-h^2(X,T_X)$.

\part[iii]\label{miniversalw} Suppose that $h^2(T_X)=1$ and
$f:\sX\to S$ contains a correctly obstructed geometric
fibre. Then $\dim_sS=h^1(X,T_X)-1$.

\begin{proof} 
For any geometric point $t\in S$, there is a smooth morphism
$\alpha :(S,t)\to D_{\sX_t}$. If $h^0(T_{\sX_t})=0$, then $D_{\sX_t}$
is universal, so that the classifying map $D_{\sX_t}\to \hat M$,
where $\hat M$ is the formal completion of the coarse moduli space
$M$ at the appropriate point, is identified with a quotient by
$\Aut_{\sX_t}$. Then,
by \ref{fibre-dimension}, the fibre-dimension of $\alpha$ is
$h^0(T_X)$. Hence $\dim_t S=\dim D_{\sX_t}+h^0(T_X)$.
If $\sX_t$ is correctly obstructed, then 
$\dim D_{\sX_t}=h^1(T_{\sX_t})-h^2(T_{\sX_t})=-\chi(T_{\sX_t})$
and $\dim_tS=h^1(T_X)-h^2(T_X)$. This proves \DHrefpart{i}, and
now \DHrefpart{ii} is clear.

Finally, \DHrefpart{iii} follows from \DHrefpart{ii} since,
if the dimension is $h^1(T_X)$,
then a small miniversal deformation is irreducible at $X$.
\end{proof}
\end{theorem}

The Picard scheme plays an important role in the study of Enriques surfaces,
so it is not surprising that for the deformation theory we need to know how the
Picard scheme varies in a family of Enriques surfaces. It would seem that the
fact that we may have $h^2(\cO)\ne 0$ could cause problems but the following
result shows that this is not so.
\begin{proposition}
Suppose that \map fXS is smooth and proper and that for every geometric
point $\bar s\to S$ we have $h^2({\cO}_{X_{\bar s}})=h^1({\cO}_{X_{\bar
s}})-b_1(X_{\bar s})/2$. Then $\Picf(X/S)$ is flat over $S$ and so if $f$ is
projective, $\Picf(X/S) /\Picf^\tau(X/S)$ is locally constant.
\begin{proof}
We may assume that $S=\Spec R$, where $R$ is local and artinian, the 
closed point is $s$ and that $k=k(s)$ is
algebraically closed. The condition $h^2({\Cal
O}_{X_s})=h^1({\cO}_{X_s})-b_1(X_s)/2$ is equivalent to
$H^2(X_s,W\Ko{X s})$ being of finite length and $ H^3(X_s,W\Ko{X s})$ 
being $V$-torsion free if $\ch k >0$ and $H^2(X_s,\cO)=0$ if $\ch k =0$. 
Thus the case of characteristic zero is clear and we
may assume that $k$ is of positive characteristic.

Let us first show that $\Pic(X/S)$ is flat along the zero section. Indeed, the
completion $T$ of the local ring at $0$ of $\Pic(X/S)$ is the quotient of some
power series ring $P$ over $R$ by an ideal $\sI$ generated by $h^2(\sO_{X_k})$
elements. However, as a finite type group scheme over a field is a local
complete intersection the conditions of the proposition implies that
the minimal number of generators of
$\ker(P/mP\to T/mT)$ is $h^1(\cO_{X_s})-b_1(X_s)/2$ 
and the hypotheses then imply that $\sI$ is a complete
intersection ideal, so that $T$ is $R$-flat. When $\Pic(X/S)$ is flat along the
zero section $R^2f_\ast\mulf$ is prorepresentable ([Ra79], 2.7.5.3) and so
is zero as its restriction to $k$ is zero, again by the assumptions and
[Ek85], Prop.~8.1. To prove flatness it is sufficient to show that any
$k$-point of $\Pic(X/S)$ lifts to a $Y$-point where $Y\to S$ is flat as then
translation gives an isomorphism of local rings. However, the obstruction for
lifting this $k$-point to an $S$-point is an element of $H^2(X,\mulf)$ and as
$R^2f_\ast\mulf=0$ this element is killed by some flat extension. Finally, if
$f$ is projective then $\Pic^\tau(X/S)$ is an open sub-algebraic space and so is
flat. This means that $\Pic(X/S)/\Pic^\tau(X/S)$ is a flat group-algebraic space
and as it is always unramified, it is \'etale. As each component is proper it is
therefore locally constant.
\end{proof}
\end{proposition}

As was mentioned our intended application of this result is to a family of
Enriques surfaces.
\begin{corollary}\label{effectivity}
\part Let $X\to S$ be a family of Enriques surfaces. Then $\Pic^\tau(X/S)$ is a
flat group scheme of order two and $\Pic(X/S)/\Pic^\tau(X/S)$ is a
locally constant sheaf of torsion free finitely generated abelian groups.

\part A formal deformation of Enriques surfaces is effective.
\begin{proof}
Indeed, an Enriques surface fulfils the conditions of the proposition and over
an algebraically closed field $\Pic^\tau$ is of order 2.

As for the second part, the first part shows that the square of any line bundle
lifts over any formal deformation and thus it is projective and hence effective.
\end{proof}
\end{corollary}

\begin{proposition}
If $Y$ is an Enriques surface over an algebraically closed field, then
$h^2(Y,T_Y)=h^0(Y,T_Y)\le 1$.
\begin{proof} $h^0=h^2$ is proved in [CD89] and both vanish if $p\ne 2$.
In [CD89] it
is shown that $h^0(Y,T_Y)=0$ if $Y$ is a $\mu_2$--surface and $h^0(Y,T_Y)=1$ if
$Y$ is an $\alpha_2$--surface. In [SB96] it is claimed that $h^0(Y,T_Y)=0$ for
$\Z/2$--surfaces, but this is false; the true statement is that $h^0(Y,T_Y)\le
1$ [ESB99].
\end{proof}
\end{proposition}

There is not much to say about miniversal deformations of
$\Z/2$- and
$\mu_2$-surfaces, except that we would like to be able
to say more about the case where there are non-trivial vector fields.
\begin{theorem}\label{notalphadef}
Let $\sX \to S \ni s$ be a family of $\mu_2$- or $\Z/2$-Enriques surfaces miniversal 
(for all deformations not just modulo $2$) at $s$ and let $X:=\sX_s$.

\part If $h^0(X,T_X)=0$ then $S$ is formally smooth over $\Z$ of relative
dimension $10$ at $s$.

\part If $h^0(X,T_X)=1$ then $S$ is flat over $\Z$ of relative dimension $11$ at
$s$ and it has a hypersurface singularity there. The incidence scheme of
$\sX\times_\Z k(s) \to S\times_\Z k(s)$ and $S\times_\Z k(s)\times X$ is
formally smooth of dimension $1$ over $k(s)$.
\begin{proof}
The first part is clear. For the second, 
note that $S\Tensor\F_2$ is formally smooth over 
$k$ outside a subset $S_{bad}$ of codimension at least $3$ and
every geometric point of $S-S_{bad}$ corresponds to a surface
with $h^0(T)=h^2(T)=0$.
The result, except the flatness, then follows from \ref{miniversalw}.
Flatness also follows from this, since $S-S_{bad}$ is $\Z$-flat.
\end{proof}
\end{theorem}
We now direct our attention towards the deformation theory of
$\alpha_2$-Enriques surfaces which, as we will see, is much richer. First some
preliminaries. A \Definition{$\Pic$-rigidification} of a family \map
\pi XS of Enriques surfaces is a trivialisation of the line bundle
$\sO_{\Picf^\tau}/\sO_S$. Note that for a family of $\mu_2$- or
$\Z/2$-surfaces there is a natural such identification coming from the unique
isomorphism of $\Picf^\tau$ with $\mu_2$ resp.~$\Z/2$ but for
$\alpha_2$-surfaces this is not the case.

If $\sG \to S$ is a finite flat group scheme of order $2$ and $t$ a generator of
$\sL:=\sO_\sG/\sO_S$ then there is a unique element, also denoted $t$, of the
augmentation ideal of $\sO_\sG$ which maps to $t$ and $\{1,\ t\}$ is a basis for
$\sO_\sG$. There then are unique elements $f$ and $g$ of $\sO_S$ such that
$t^2=ft$ and the coproduct takes $t$ to $1\tensor t-g(t\tensor t)+t\tensor 1$;
one sees immediately that this defines a group scheme if and only if
$fg=2$. Clearly if $t$ is replaced by $\lambda t$ then $f$ is replaced by
$\lambda f$ and $g$ by $\lambda^{-1}g$. This means that 
we have maps $\sL \to \sO_S$ and $\sL^{-1} \to \sO_S$ taking
$t$ to $f$ and $t^{-1}$ to $g$ respectively such that the map induced by
multiplication $\sO_S =\sL\Tensor\sL^{-1} \to \sO_S$ takes $1$ to $2$. The
subschemes defined by the images of $\sL$ resp.~$\sL^{-1}$ will be called the
\Definition{locus of infinitesimality} and the \Definition{locus of unipotence}
respectively. 
\begin{proposition} \label{0.6}
Let $X$ be an $\alpha_2$-Enriques surface over $k$. Its
miniversal deformation space $S=\Spec R$
is isomorphic to $\Spec W(k)[[{x_1,\dots,x_{12}}]]/(FG-2)$ where
$F,G$ lie in the ideal 
$(2,x_1,\dots,x_{12})$. In particular, $S$ is a regular scheme.  The
isomorphism may be chosen such that there is a $\Pic$-rigidification of a
miniversal deformation over $S$
such that $F$ and $G$ map to
the $f$ and $g$ just defined.
\begin{proof}
Since $h^1(X,T_X)=12$,
we can write a hull $R$ as a quotient of the power series ring
$W(k)[[{x_1,\dots,x_{12}}]]$.
If we choose a $\Pic$-rigidification we get elements $f$ and
$g$ in $R$ such that $fg=2$. This means that we can write $R$ as
a quotient of $W(k)[[{x_1,\dots,x_{12}}]]/(FG-2)$. However, as $f$ and 
$g$ lie in
the maximal ideal of $R$, $W(k)[[{x_1,\dots,x_{12}}]]/(FG-2)$ is
regular and as $h^2(X,T_X)=1$, $R$ has dimension $\ge 11$. This shows that the
quotient map must be an isomorphism.
\end{proof}
\end{proposition}
In order to get more information we will now determine the image of $f$ and $g$
in the cotangent space of the hull. Note that as $\omega_X$ for an
$\alpha_2$-surface is trivial we can identify $H^1(X,T_X)$ with
$H^1(X,\Omega^1_X)$ and by Serre duality we can thus identify the cotangent
space of a miniversal deformation with $H^1(X,\Omega^1_X)$. The images 
sought will therefore be cohomology classes in $H^1(X,\Omega^1_X)$. They are
only defined as elements when a $\Pic$-rigidification has been chosen and we
will start by describing certain elements in the cohomology of $X$ that can be
constructed from such a rigidification.
\begin{definition}\label{multelts}
Let $(X,t)$ be a $\Pic$-rigidified $\alpha_2$-Enriques surface over a field
$k$. Associated to it are the following elements of Hodge cohomology.

\part 
$\beta \in \hod X10 k$ which corresponds to $t$ under the identifications
between $\hod X10 k$ and the augmentation ideal of $\Cal O_{\Pic^\tau(X/k)}$. 

\part 
$\check\beta \in \hod X12 k$ which is dual to $\beta$ (i.e., $\Tr(\beta\check\beta)=1$).

\part
$\eta \in \hod X01 k$ which is the image of the inclusion $H^0(X,B_1)
\hookrightarrow \hod X01 k$ of the (uniquely determined) element of $H^0(X,B_1)$
which maps to $\eta$ under the boundary map of the long exact sequence of
cohomology associated to $\shex{\sO_X}{\sO_X}{B_1}$.

\part $\beta^2 \in \hod X20 k$ which is the cup square of $\beta$.

\part $\check\beta_2 \in \hod X02 k$ which is the dual of $\beta^2$.

\part $D \in H^0(X,T_X)$ is the global vector field defined by
$D(f)\check\beta_2=\eta\wedge df$.
\end{definition}
Some of the statements we will make about these elements will be independent on
the rigidification chosen. In such situations we will use these elements without
always explicitly choosing a rigidification. Let us also note that with respect
to changing a rigidification $\beta, \beta^2, \eta, \check\beta,
\check\beta^2, D$ will be homogeneous of degrees $1, 2, 2, -1, -2, -4$
respectively.

Using the computation of the cohomology of Enriques surfaces it is clear that
all these elements are non-zero. From them we can construct two elements of
$H^1(X,\Omega^1_X)$ namely $d\beta$ and $\eta\beta$.
\begin{proposition}
Let \map\pi \sX S be a family of Enriques surfaces.

\part[i] The infinitesimal locus of $\Picf^\tau(\sX/S)$ equals the locus where
$\cO_\sX$ is cohomologically flat in degree $1$ and $R^1\pi_*\cO_\sX$ has rank
$1$ (cf., \ref{flatstrat}).

\part[ii] The unipotent locus of $\Picf^\tau(\sX/S)$ equals the locus where
$\Omega^1_\sX$ is cohomologically flat (in degree $0$) and $\pi_*\Omega^1_\sX$
has rank $1$ (cf., \ref{flatstrat}).

\part[iii] Let $X$ be an $\alpha_2$-Enriques surface. Then the tangent space of the
infinitesimal locus of $\Picf^\tau(X)$ in a miniversal deformation of $X$ is the
orthogonal complement of $d\beta$, where $\beta$ corresponds to a
$\Pic$-rigidification as above.

\part[iv] Let $X$ be an $\alpha_2$-Enriques surface. Then the tangent space of the
unipotent locus of $\Picf^\tau(X)$ in a miniversal deformation of $X$ is the
orthogonal complement of $\eta\beta$, where $\beta$ and $\eta$ correspond to a
$\Pic$-rigidification as above.
\begin{proof}
We start with \DHrefpart{ii}: We are immediately reduced to assuming that
$S=\Spec R$ is
local artinian with the fibre $X$ of $\pi$ over the closed point of $S$ an
$\alpha_2$-surface and we want to show that $\Omega^1_\sX$ is cohomologically flat
if and only if $\Picf^\tau(\sX/S)$ is unipotent.

We start by showing the equality modulo $2$. In general, for a family \map
\pi\sX S of smooth and proper varieties in characteristic $p$ we have a morphism
${}_p\Picf(\sX/S) \to \pi_*\Omega^1_{\sX/S}$ defined as follows (following the
lines of [Od69]). If we represent a line bundle on $\sX$ by a cocycle
$f_{ij}$, then if the $p$'th power of it is trivial there are functions $f_i$
such that $f^p_{ij}=f_if^{-1}_j$. Applying the logarithmic derivative gives us a
global $1$-form $d\log(f_i)$. This gives a map upon sheafification and it is
natural for all base changes in the sense that for any map \map fTS the
following diagram commutes
$$\xymatrix{
{f_*({}_p\Picf(\sX\times_ST/T))} \ar[r] & {f_*\pi_*\Omega^1_{\sX\times_ST/T}}\\
{{}_p\Picf(\sX/S)} \ar[u] \ar[r] & {\pi_*\Omega^1_{\sX/S}.}\ar[u]
}$$
Returning to the case at hand we note that this map is injective for 
$X \to k$. Note that for any rigidified order $2$ finite flat group scheme $\sG$ over
$S$ with functions $f$ and $g$ as above, $\Hom(\sG,\add)$ equals the
annihilator of $g$ in $R$. Assume first that $\Omega^1$ is cohomologically flat
and choose a $\Pic$-rigidification of $\sX$ (which gives us $f$ and $g$) as well
as a generator of $^0(\sX,\Omega_{\sX/S})$. In particular this means that if
\map f k S is the inclusion of the closed point then $H^0(\sX,\Omega^1_{\sX/S})
\to H^0(X,\Omega^1_X)$ is surjective but also that the map ${}_2\Picf(\sX/S) \to
\pi_*\Omega^1_{\sX/S}$ can be seen as a map ${}_2\Picf(\sX/S) \to \add$. By the
commutativity of the diagram this means that the reduction map $\Ann_R g \to k$ 
takes an element of $\Ann_R g$ to a non-zero element of $k$ which means that by 
Nakayama's lemma that $\Ann_R g = R$ and so $g=0$.

Conversely, assume, with the same notations as previously, that $g=0$ so that
${}_2\Picf(\sX/S)$ is isomorphic to $\alpha_2$. If $\Omega^1_{\sX/S}$ is not
cohomologically flat then the map $H^0(\sX,\Omega_{\sX/S})
\to H^0(X,\Omega^1_X)$ is not surjective and thus zero. Restricting to the flat
site of $S$ it then remains zero. Looking at the
commutative diagram and using the injectivity of ${}_2\Picf(X/k) \to
H^0(X,\Omega_X)$ we see that this means that ${}_2\Picf(\sX/S) \to f_*({}_2\Picf(\sX/S))$ 
is zero on the flat site of $S$. However, this is false as ${}_2\Picf(\sX/S)$ is 
isomorphic to $\alpha_2$ and $\alpha_2 \to f_*\alpha_2$ is non-zero on the flat
site of $S$ as can be seen by considering it for $S[\delta]$.

This proves the statement when restricted to characteristic $2$. For the general
case we may assume that the base is a deformation hull and note that the ideal
generated by $g$ (using notation as before) contains $2$ as $fg-2=0$. As the
ideal are equal modulo $2$ this means that the ideal $I$ defining the locus where
$\Omega^1$ is cohomologically flat with $\pi_*\Omega$ of rank $1$ is contained
in the ideal generated by $g$. Now by what we have just proved there is an
element $g'$ in $I$ which modulo $2$ is $g$. As $g'=hg$ for some element we see
immediately that $h$ is a unit so the two ideals are the same.

The proof of \DHrefpart{i} is similar but simpler.

Turning to \DHrefpart{iii} we need by \DHrefpart{i} to show that for a
deformation $\sX$ of $X$ over $k[\delta]$ the boundary map $H^1(X,\cO_X) \to
H^2(X,\cO_X)$ is zero if and only if the deformation class in
$H^1(X,T_X)=H^1(X,\Omega^1_X)$ is orthogonal to $d\beta$. In general if $\sE$ is
a sheaf of tensors on $X$, we have an action $L_D$ by a vector field $D$. The
boundary map $H^i(X,\sE) \to H^{i+1}(X,\sE)$ given by a deformation over the
dual numbers is given by cupping with $L_\nu$, where $\nu \in H^1(X,T_X)$ is the
class of the deformation. When $\sE=\cO_X$, $L_D$ is simply $D$ itself. Now, we
have to look into the identification between $T_X$ and $\Omega^1_X$. If $\omega$
is a non-zero global $2$-form then the derivation $D$ corresponds to the form
$\nu$ when $D(f)\omega=\nu\wedge df$ for all $1$-forms $\nu$. This implies that
$(D \contract \psi)\omega = \nu \wedge \psi$ for all $1$-forms $\nu$ and
$\psi$. Now choose a $\Pic$-rigidification and let $\omega=\check\beta_2$. Then
if we let $\nu \in H^1(X,\Omega^1_X)$ correspond to $D \in H^1(X,T_X)$ then if
$\alpha \in H^1(X,\cO_X)$ we have that $L_D(\alpha)\check\beta_2=\nu \wedge
d\alpha$. As $H^1(X,\cO)$ is spanned by $\beta$ and $H^2(X,\cO_X)$ by $\beta^2$
which is dual to $\check\beta_2$ this shows that the boundary map is zero if and
only if $\nu$ is orthogonal to $d\beta$.

To prove \DHrefpart{iv} similarly requires an elucidation of $L_D(\eta)$. Here
we use Cartan's formula $L_D\omega=d(D\contract\omega)+D\contract d\omega$. As
$\eta$ is closed we get that $L_D\eta=d(D\contract \eta)$. As \map
d{H^1(X,\cO_X)}{H^1(X,\Omega^1_X)} is injective this shows that $L_D\eta$ is
zero if and only if $D\contract\eta$ is. Now multiplication by $\beta$ on
$H^1(X,\cO_X)$ is injective so this is zero if and only if
$\beta(D\contract\eta)=0$ and this in turn is true if and only if
$\beta(D\contract\eta)\check\beta_2=0$ and this equals $\beta\eta\nu$.
\end{proof}
\end{proposition}
This proposition has several interesting consequences.
\begin{corollary}\label{alphadef}
Let $X$ be an $\alpha_2$-Enriques surface over $k$
and let $S=\Spf R$ be a deformation hull of
$X$. Choose a $\Pic$-rigidification of the miniversal deformation of $X$ over
$S$ so that the elements $f$, $g$, $\beta$, $\eta$ and so on are defined.

\part[i] The formal subscheme of $S$ where $\Picf^\tau$ is infinitesimal is
formally smooth of dimension $11$ over $k$. In particular 
$R\cong W(k)[[{f,x_2,\dots,x_{12}}]]/(fg-2)$.

\part[ii] The formal subscheme of $S$ where $\Picf^\tau$ is unipotent is
formally smooth over $k$ precisely when $\eta\beta \ne 0$. It always has
dimension $11$ and embedding dimension $\le 12$.

\part[iii] The formal subscheme of $S$ where $\Picf^\tau$ is isomorphic to
$\alpha_2$ is formally smooth over $k$ precisely when $d\beta$ and
$\eta\beta$ are linearly independent. It is reduced and $10$-dimensional
and has embedding dimension $\le 11$. In particular, $g$ is 
divisible by neither $f$ nor $2$ and $R$ is flat over $W(k)$.

\part[iv] 
$R\cong W(k)[[{f,g,x_3,\dots,x_{12}}]]/(fg-2)$ when $d\beta$ and $\eta\beta$ are
linearly independent. In particular the map from the deformation problem of $X$
to that of deforming $\Picf^\tau(X)$ as a flat group scheme is formally smooth
in this case.

\part $X$ can be lifted to characteristic $0$. It cannot be lifted
to $W(k)$.

\part If $d\beta$ and $\eta\beta$ are non-proportional then $X$ can be
lifted to any DVR of characteristic zero that is ramified over $W(k)$.

\part If $d\beta$ and $\eta\beta$ are proportional but $\eta\beta\ne 0$ then
$X$ has a lifting over some characteristic zero DVR  with absolute 
ramification index $2$.

\part If $\eta\beta=0$ then $X$ cannot be lifted to any DVR with
absolute ramification index $\le 2$.
\begin{proof}
Recall that $f,g\in \frak m_R$, the maximal ideal of $R$, and that
always $d\beta \ne 0$. 
This means that $f$ is non-zero
modulo $\frak m_R^2$, where $\frak m_R$ is
the maximal ideal of $R$. Hence we can chose $f$ as
one of the parameters when writing $R$ as a quotient of a power series ring
over $W(k)$, so
that $R$ is isomorphic to the quotient ring
$W(k)[[{f,x_2,\dots,x_{12}}]]/(fg-2)$. This implies
that $R/Rf$ is isomorphic to $k[[{x_2,\dots,x_{12}}]]$ which proves
\DHrefpart{i}.

To continue with \DHrefpart{ii} we get similarly that $g$ modulo 
$\frak m_R^2$ is non-zero
if and only if $\eta\beta$ is. This proves the first part. To prove the second
we need to show that $g \ne 0$. If $g=0$ then $R$ has characteristic $2$ but
this is not possible as that would force any versal deformation of $X$ to have
the base killed by $2$ and by \DHrefpart{i}, there are $\mu_2$-surfaces in any
neighbourhood of any versal deformation. By openness of versality a versal
deformation of such $\mu_2$-surface would also have base killed by $2$ which is
absurd.

As for \DHrefpart{iii} $d\beta$ and $\eta\beta$ are linearly independent if and
only if the images of $f$ and $g$ modulo $\frak m_R^2$ are, 
and the $\alpha_2$-locus is smooth if and only if these images are linearly
independent, since it is defined by the vanishing of $f$ and $g$.
In general, $g$ cannot be divisible by $f$, because if so $R/2R$ is
non-reduced. However, a versal deformation is generically reduced and $R/2R$ is a
local complete intersection. Similarly, if $g$ is divisible by $2$ then 
$g$ would vanish in $R/2R$ and then
any mod $2$-deformation of $X$ would have unipotent
$\Picf^\tau$, in contradiction to the fact that,
as we have just seen, there are $\mu_2$-deformations. 
Flatness over $\W(k)$ now follows as $fg-2$ is not divisible by $2$.
Further, if $f$ and $g$ modulo $\frak m_R^2$ are linearly independent, 
then they may be chosen as
parameters when writing $R$ as a quotient of $W(k)[[{x_1,\dots,x_{12}}]]$,
which proves \DHrefpart{iv}. 

Next, we consider the question of
liftability. Liftability follows from the flatness over
$W(k)$ of a versal deformation. For any lifting over a complete DVR $V$ we get
elements $\ovl f$ and $\ovl g$ in $\frak m_V$ such that $\ovl f\ovl g=2$ which forces
$V$ to be ramified. If $f$ and $g$ can be chosen to be parameters 
in $R$ then it is clear that
we can find a lifting over any ramified $V$. If
$d\beta$ and $\eta\beta \ne 0$ are proportional then we can write $g$ as
$\lambda f+O(2)$ where $\lambda$ is a unit in $W(k)$,
so that $R/(x_2,\ldots,x_{12})$ is
of the form
$\pow[W(k)]{f}/(h(f))$, where $h(f)$ is of the form $\ovl\lambda f^2+O(2)$ modulo
$2$. By the Weierstrass preparation theorem $\pow[W(k)]{f}/(h(f))$ is then free of 
rank $2$ over $W(k)$ and there is a lifting over it. 
Pull back to the normalization of this
ring to get the result.

Finally, if $\eta\beta =0$, then $g\in\frak m_R^2$, so that
$2\in\frak m_R^3$, and we are done.
\end{proof}
\end{corollary}

We will now spend some time studying multiplicative $\alpha_2$-surfaces. 
\begin{proposition}
Let $X$ be a multiplicative $\alpha_2$-Enriques surface.

\part $D$ acts non-trivially on $H^1(X,\cO_X)$.

\part $\eta\beta d\beta \ne 0$ and in particular $d(\eta\beta)\ne 0$.

\part $d\beta$ and $\eta\beta$ are linearly independent.
\begin{proof}
Rescale $D$ so that $D^2=D$. This means that $\cO_X$ is the direct
sum $\sE\Dsum\sF$ where $D$ acts by zero on $\sE$ and by $1$ on $\sF$. If $D$
acts trivially on $\hod X10 k$ then $h^1(X,\sE)=1$, but $\sE$ is the structure
sheaf of the quotient of $X$ by $D$ and this quotient has only rational double
points with minimal resolution rational or K3 so $h^1(X,\sE)=0$. This means that
$D\beta\ne 0$ and as $D^2=D$ we get $D\beta=\beta$.  Now there is a $0\ne
\lambda\in k$ s.t.~for $f\in \cO_X$, $D(f)\check\beta_2=\lambda\eta\wedge df$,
as this is true for any non-zero global vector field. Hence $0\ne
\beta\check\beta_2=D(\beta)\check\beta_2=\lambda\eta\wedge d\beta$;
i.e., $\eta\wedge d\beta$ and so $d\beta\eta\beta$ is different from zero. By
duality we have that multiplication by $\beta$ on $H^1(X,\Omega^2_X)$ is
injective so that $\eta d\beta=d(\eta\beta) \ne 0$. Finally $d\beta d\beta=0$
so that $d\beta$ and $\eta\beta$ are linearly independent.
\end{proof}
\end{proposition}
\begin{corollary}
\part A versal deformation of $\alpha_2$-Enriques surfaces is smooth at a surface of
multiplicative type.

\part The moduli problem of $\Pic$-rigidified Enriques surfaces is
prorepresentable and smooth of dimension $11$ at a multiplicative
$\alpha_2$-Enriques surface.
\begin{proof}
The first part follows immediately from the proposition and \ref{alphadef}. For the 
second we note that as the global non-trivial vector fields act non-trivially on 
$H^1(X,\cO_X)$ we see that two different $\Pic$-rigidifications of the same
deformation of $X$ over a local artinian base $S$ that are the same when
restricted to a subscheme $T \inj S$ which is defined by an ideal of length $1$
are isomorphic. Hence, by decomposing a general map into such substeps we see
that over a local artinian base isomorphism classes of deformations of $X$ are in 
bijection with isomorphism classes of $\Pic$-rigidified deformations of
$X$. Similarly, using small steps as above one shows that the automorphism group 
of a $\Pic$-rigidified deformation of $X$ is trivial, hence that deformation
problem is prorepresentable.
\end{proof}
\end{corollary}
\begin{remark}
This gives, in the multiplicative type case, a clearer picture of the rather
strange non-separation of the moduli stack of Enriques surfaces at
$\alpha_2$-surfaces: The moduli stack of $\Pic$-rigidified $\mu_2$-, $\Z/2$-, or 
multiplicative $\alpha_2$-Enriques surfaces has an unramified diagonal and the
stack of un-rigidified surfaces is the quotient of it by the natural
$\mul$-action.
\end{remark}
Even though we are not going to use it we finish with a description of an
invariant of a multiplicative $\alpha_2$-Enriques surface which gives an \'etale
map from the moduli stack.
\begin{proposition}
Let \map\pi XS be a $\Pic$-rigidified family of multiplicative
$\alpha_2$-Enriques surfaces. Then there is a canonical map $\Picf(X/S)\to \add$
such that
$\Picf^\tau(X/S)$ maps isomorphically to $\alpha_2$. Dividing by
$\Picf^\tau(X/S)$ gives an additive map $NS\to \add$. Suppose that $S$ is of
finite type over a field $k$ and that we have a marking $\Z^{10} \cong NS$ so
that the map $NS\to \add$ can be interpreted as a map \map\rho
S{Hom(\Z^{10},\add)}. At any point at which $\pi$ is miniversal, $\rho$ is
\'etale.
\begin{proof}
We can define $\beta$ and $\eta$ as before even though we work in a
family. Define $\Picf(X/S)\to \add$ by $\sL\mapsto \Tr(\eta\beta
c_1(\sL))$ where $c_1(\sL)$ is the de Rham Chern class. Certainly
$\Picf^\tau(X/S)$ maps into $\alpha_2$ and to show that it maps onto it suffices
to show that the tangent map is surjective. Hence we may suppose that
$S= k[\delta]$ for some field $k$ and we want to show that $\eta\beta
c_1(\sL)\ne 0$ for $\sL$ given by $1+\delta\beta \in
H^1(X,1+\delta\cO_X)$. However, $c_1(\sL)=d\log(1+\delta\beta)=\delta d\beta$
and we have just seen that $\eta\beta d\beta\ne 0$. Furthermore, as a miniversal
deformation is smooth, to show that $\rho$ is \'etale it is sufficient to show
that it is an isomorphism on tangent vectors. As the tangent spaces have the
same dimension it will be enough to show that the tangent map is injective.
Suppose therefore that we have a deformation $Y\to \Sp k[\delta]$, the pullback
of $X$ by a map $k[\delta]\to S$, for which the map $NS\to \add$ is
constant. As we have seen, for any deformation, the extension
\begin{displaymath}
\shex{\Picf^\tau}{\Picf}{NS}
\end{displaymath}
is isomorphic to the pull back of the extension
\begin{displaymath}
\shex{\alt p}\add\add
\end{displaymath}
by the map $NS\to \add$. Hence $\Picf(Y/\Sp k[\delta])=
\Picf^\tau(Y/\Sp k[\delta])\bigoplus NS$ and in particular $Pic(Y)\to Pic(Y_k)$
is surjective. This means that the element $\alpha$ in $H^1(Y_k,T^1_{Y/k})$
classifying $Y\to \Sp k[\delta]$ has cup product zero with
$c_1(Pic(Y_k))$. However, we have seen that the fact that $Y$ is an
$\alpha_2$-surface means that $\alpha$ has cup product zero with $\eta\beta$ and
$d\beta$. Now, as $d(\eta\beta)\ne 0$ it is easy to see that $\eta\beta$,
$d\beta$ and $c_1(Pic(Y_k))$ span \hod {Y_k}11 k\ as a $k$-vector space and so by
duality and the isomorphism $\omega_{Y_k}\liso\ko{Y_k}$, $\alpha=0$.
\end{proof}
\end{proposition}
\end{section}

\begin{section}{Marking and polarizations of Enriques surfaces}\label{blank}
\label{Enr marking}

We start this section with a definition that provides a solution to
the problem of how to
mark an Enriques surface and simultaneously polarize it. Then we discuss
the morphism from the corresponding stack to the stack considered by
Liedtke [L10] of Enriques surfaces with a \emph{Cossec--Verra}
polarization and derive an upper bound on the ramification
involved in lifting an Enriques surface to characteristic zero. 
After that we go on to consider
those (singular) K3 surfaces that arise as canonical double covers of
an Enriques surface over an algebraically closed field of characteristic $2$.

Denote by $E$ the lattice $E_{10}(-1)$, $O(E)$
its orthogonal group and $O^+(E)$ the index $2$ subgroup of
$O(E)$ consisting of elements that preserve the two
cones of positive vectors in $E_\R$. It is known [CD]
that $O^+(E)$ is the Weyl group $W(E)$, the group
generated by reflections in the roots of $E$. For any
Enriques surface $Y$ over an algebraically closed field,
$\Num(Y)\cong E$.  Fix, once and for all, a chamber $\sD_0$ defined by the roots
(that is, the $(-2)$-vectors)
in the positive cone of $E\Tensor\R$. (We shall not always be scrupulous in
distinguishing between $\sD$ and its closure.)  This defines a root basis
$\alpha_1,\ldots,\alpha_{10}$ of $E$ that in turn defines a Dynkin diagram of
type $E_{10}=T_{2,3,7}$. We recover $\sD_0$ as $\sD_0=\sum\R_{\ge 0}\varpi_i$, where
$\varpi_1,\ldots,\varpi_{10}$ are the fundamental dominant weights defined by
the root basis. That is, $\varpi_i.\alpha_j=\delta_{ij}$.
We label the simple roots $\alpha_1,\ldots,\alpha_{10}$ according
to the diagram
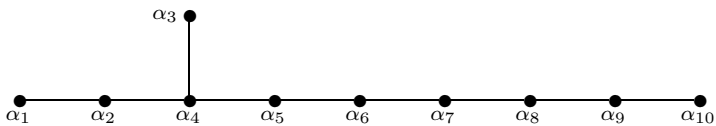
{

\catcode`\@=11
\let\e@\expandafter
\def\d@nkinsize{1cm}
\newcount\d@nkincount
\d@nkincount"00
\newif\if@ne

\def\LINE{\hbox to \hsize}
\def\setdynkin#1=#2 {\csname d@nkin@#1\endcsname{#2}}
\def\d@nkin@size#1{\def\d@nkinsize{#1}}
\def\d@nkin@slope#1{\d@nkincount#1}
\def\d@nkin@lines#1{\ifnum#1=1\@netrue\else\@nefalse\fi}
\newif\ifd@dotchoice
\newif\ifd@tags
\def\d@nkin@dotchoice#1{\if#1y\d@dotchoicetrue\else\d@dotchoicefalse\fi
\d@setdots}
\def\d@nkin@tags#1{\if#1y\d@tagstrue\else\d@tagsfalse\fi\d@setdots}
\def\d@setdots{%
\ifd@dotchoice
\def\d@nkinball{\e@\dynkinball\e@{\the\t@@@}}%
\ifd@tags
 \let\d@loopargs\d@loopargsa
 \let\count@rg\count@rga
 \let\init@rg\init@rga
 \let\pr@pare\pr@parea
\else
 \let\d@loopargs\d@loopargsbc
 \let\count@rg\count@rgbc
 \let\init@rg\init@rgbc
 \let\pr@pare\pr@pareb
\fi
\else
\def\d@nkinball{\dynkinball}%
\ifd@tags
 \let\d@loopargs\d@loopargsbc
 \let\count@rg\count@rgbc
 \let\init@rg\init@rgbc
 \let\pr@pare\pr@parec
\else
 \let\d@loopargs\d@loopargsd
 \let\count@rg\count@rgd
 \let\init@rg\init@rgd
 \let\pr@pare\pr@pared
\fi
\fi}

\def\d@nkin@ball#1{\e@\let\e@\dynkinball\csname d@nkin#1\endcsname}

\def\spl@ttoksa#1\@t#2\@t#3\t@end{\t@@@={#1}\t@@@@={#2}\global\t@@={#3}}%
\def\spl@ttoksb#1\@t#2\t@end{\t@@@={#1}\t@@@@={}\global\t@@={#2}}%
\def\spl@ttoksc#1\@t#2\t@end{\t@@@={}\t@@@@={#1}\global\t@@={#2}}%

\def\pr@parea{\e@\spl@ttoksa\the\t@@\t@end}%
\def\pr@pareb{\e@\spl@ttoksb\the\t@@\t@end}%
\def\pr@parec{\e@\spl@ttoksc\the\t@@\t@end}%
\def\pr@pared{\t@@{}\t@@@{}}%

\def\n@metop{\hbox to 0pt{\kern-\d@nkinrightdim\vbox to
\wd0{\vss\offinterlineskip\hsize\ht0\pr@pare
\centerline{\the\t@@@@}\kern2pt\d@nkinball}\hss}}%
\def\n@medown{\hbox to 0pt{\kern-\d@nkinrightdim%
\vtop to \wd0{\offinterlineskip\hsize\ht0\pr@pare%
\d@nkinball\kern2pt\centerline{\the\t@@@@}\vss}\hss}}%
\def\n@meright{\hbox to 0pt{\kern-\d@nkinrightdim\vtop to
\wd0{\offinterlineskip\hsize\ht0\pr@pare
\hbox{\d@nkinball\kern2pt\vbox to 
\wd0{\parindent0pt\vss\rlap{\the\t@@@@}\vss}}\vss}\hss}}%

\def\n@meleft{\hbox to 0pt{\hss\kern-\d@nkinrightdim\vtop to
\wd0{\offinterlineskip\hsize\ht0\pr@pare
\hbox{\vbox to
\wd0{\parindent0pt\vss\llap{\the\t@@@@\kern2pt}\vss}\d@nkinball}\vss}}}

\def\Tl@ne#1{\kern\d@nkinrightdim
\vbox{\hbox to #1{}\hrule}\kern\d@nkinrightdim}%
\def\tl@ne#1{\kern\d@nkinrightdim
\vbox{\hbox to #1{\vrule width1.5mm height.4pt
\if@ne\dimen0 1mm\else\dimen0 2mm\fi
\leaders\hbox to \dimen0{\hss.\hss}\hfill
\vrule width1.5mm height.4pt}}\kern\d@nkinrightdim
\let\hl@ne\Tl@ne}%
\def\b{\global\let\hl@ne\tl@ne}%
\def\vl@ne#1{\raise\dimen0\hbox{\vbox to #1{}\vrule}\kern-.4pt}%
\def\d@nkinone{\dimen9\ht0\advance\dimen9-.4pt
\raise-0.5\dimen9\n@medown}%
\def\d@nkinrightfork{%
\hbox{\kern\d@nkinrightdim\kern-0.1464466094\wd0
\raise-\f@rkraise
\vbox{\offinterlineskip\parindent0pt
\dimen0 0.70710678108\ht0
\if@ne
\advance\dimen0 \ht1
\else
\advance\dimen0 2\ht1
\fi
\hbox{\raise0.85355339059\wd0
\@pline%
\kern0.35355339059\wd0
\raise\dimen0
\n@meright}%
\kern-0.1464466094\ht0
\hbox{\raise0.85355339059\wd0
\d@wnline%
\kern-0.1464466094\wd0
\kern\d@nkinrightdim\n@meright}}}}%
\def\d@nkinleftfork{\kern.5\wd0\raise-\f@rkraise
\vbox{\offinterlineskip\parindent0pt
\dimen0 0.70710678108\ht0
\hbox{\d@wnline
\kern-2\wd1
\if@ne\dimen9 \ht1\else\dimen9 2\ht1\fi
\advance\dimen9-0.1464466094\ht0\if@ne\kern\wd1\fi
\raise\dimen9
\n@meleft}%
\kern0.70710678108\ht0
\hbox{\raise0.85355339059\wd0
\@pline%
\if@ne\kern-\wd1\else\kern-2\wd1\fi
\n@meleft}}\kern0.35355339059\wd0
\if@ne\kern\wd1  \else\kern2\wd1\fi
}%
\def\d@nkinright{\hl@ne\d@nkinstep\d@nkinone}%
\def\d@nkinup{{\vl@ne\d@nkinstep\kern0.5\wd0\raise\dimen2\n@meleft}%
\kern-0.5\wd0}%
\def\d@nkinuptwo{\vl@ne\d@nkinstep\kern0.5\wd0\raise\dimen2\n@meleft%
\kern-0.5\wd0\dimen9\dimen2
\advance\dimen9 0.5\ht0
\raise\dimen9\hbox{\vl@ne\d@nkinstep\kern0.5\wd0\raise\dimen2\n@meleft}%
\kern-0.5\wd0}%
\toksdef\t@0\toksdef\t@@1\toksdef\t@@@2%
\toksdef\t@@@@3%

\def\dynkininit#1{%
\t@{}\t@@{}\t@@@{}\t@@@@{}%
\dimendef\d@nkinstep4
\e@\d@nkinstep\d@nkinsize
\setbox0\hbox{\ifd@dotchoice\dynkinball{#1}\else\dynkinball\fi}%
\dimendef\d@nkinraise6
\dimendef\d@nkinrightdim7
\d@nkinrightdim\wd0
\divide\d@nkinrightdim2
\font\l@net=line10
\dimendef\f@rkraise3
\setbox1\hbox{\l@net\char\d@nkincount}%
\f@rkraise 1.20710678108\ht0
\advance\f@rkraise-.2pt
\if@ne
\advance\f@rkraise \ht1
\else
\advance\f@rkraise 2\ht1 
\fi
\dimen0\ht0
\divide\dimen0 2
\dimen2\d@nkinstep
\advance\dimen2\dimen0
\let\hl@ne\Tl@ne
\edef\@pline{\vbox{\offinterlineskip\parindent0pt\l@net
\if@ne\hsize \wd1\char\d@nkincount
\else
\hsize 2\wd1
\LINE{\hskip\wd1\char\d@nkincount\hss}\par
\char\d@nkincount\fi}}%
\edef\d@wnline{\vbox{\advance\d@nkincount"40\l@net
\offinterlineskip\parindent0pt\if@ne\char\d@nkincount\hsize \wd1
\else
\hsize 2\wd1
\char\d@nkincount\par
\LINE{\hskip\wd1\char\d@nkincount\hss}\fi}}%
}

\def\count@rga#1 #2,{\count255 #1\count253 #1\dynkininit{#2}\next@rg#2,}
\def\count@rgbc#1 #2 {\count255 #1\count253 #1\dynkininit{#2}\next@rg#2 }
\def\count@rgd#1 {\count255 #1\count253 #1\dynkininit{}\next@rg}
\def\init@rga#1,{\dynkininit{#1}\next@rg#1,}
\def\init@rgbc#1 {\dynkininit{#1}\next@rg#1 }
\def\init@rgd{\dynkininit{}\next@rg}

\def\d@loopargsa#1{\def\next@rg##1,##2 {\t@@=\e@{\the\t@@##1\@t##2\@t}\next}
\def\next{\ifnum\count255>1\advance\count255 -1\e@\next@rg\else#1\fi}}
\def\d@loopargsbc#1{\def\next@rg##1 {\t@@=\e@{\the\t@@##1\@t}\next}%
\def\next{\ifnum\count255>1\advance\count255 -1\e@\next@rg\else#1\fi}}
\def\d@loopargsd#1{\def\next@rg{#1}}

\def\read@rgs#1\afterread{
\d@loopargs{#1}
\ifrepe@t
\e@\count@rg
\else
\count255 \count254
\count253 \count254
\e@\init@rg
\fi}

\newif\ifrepe@t
\def\d@nkrepeat{rightrepeat}
\def\d@nkstop{st@p}
\def\d@one{1}
\def\d@right{1}
\def\d@rightrepeat{0}
\def\d@up{1}
\def\d@uptwo{2}
\def\d@rightfork{2}
\def\d@leftfork{2}

\def\countd@nkin#1,{\def\d@rep{#1}%
\ifx\d@nkrepeat\d@rep\repe@ttrue\fi\ifx\d@nkstop\d@rep\else
\advance\count254 \csname d@#1\endcsname
\e@\countd@nkin\fi}

\def\dod@nkin#1,{\def\d@rep{#1}%
\ifx\d@nkstop\d@rep\else
\csname d@nkin#1\endcsname
\e@\dod@nkin\fi}

\def\d@nkinrightrepeat{\loop
\ifnum\count253>0
\advance\count253 -1
\d@nkinright%
\repeat}

\d@setdots
\def\e@ttok#1#2\@nd{\t@@={#2}}

\def\rawdynkin#1;{\count254 0
\repe@tfalse
\countd@nkin#1,st@p,%
\bgroup
\read@rgs
\advance\count253 -\count254
\hbox{
\dod@nkin#1,st@p,%
\kern\d@nkinrightdim}%
\egroup\afterread
}
\def\dynkinA{\rawdynkin one,rightrepeat;}
\def\dynkinD{\rawdynkin one,rightrepeat,rightfork;}
\def\dynkinDtilde#1 {\count252 #1\advance\count 252 1
\rawdynkin leftfork,one,rightrepeat,rightfork;{\count252} }
\def\dynkinE{\rawdynkin one,right,right,up,rightrepeat;}
\def\dynkinEtilde#1 {\let\d@it\relax\ifcase#1\or\or\or\or\or\or
\def\d@it{\rawdynkin one,right,right,uptwo,right,right;}\or
\def\d@it{\rawdynkin one,right,right,right,up,right,right,right;}\else
\count252 #1\advance\count 252 1
\def\d@it{\rawdynkin one,right,right,up,rightrepeat;{\count252} }\fi
\d@it}

\def\d@nkinbullet{\vbox to 3.85pt{\kern-.1pt\hbox to 3.7pt
{\kern-.43pt$\bullet$\hss}\vss}}

\def\d@nkinbulletandcirc#1{\vbox to 3.85pt{\kern-.1pt\hbox to 3.7pt
{\kern-.43pt\ifx #1c$\circ$\else$\bullet$\fi\hss}\vss}}

\def\d@nkinbigcirc{\vbox to 8.4pt{\kern-.1pt\hbox to 8.8pt
{\kern-.4pt$\bigcirc$\hss}\vss}}

\setdynkin ball=bullet
\catcode`\@=12


\setdynkin tags=y \setdynkin slope="25
\def\x#1{$\scriptstyle#1$}
\def\present#1{\hbox{\hbox to2cm{$#1$:\hfil}\hskip1.5cm}}
\par\medskip
\dynkinE 10 \x{\alpha_1} \x{\alpha_2} \x{\alpha_4} \x{\alpha_3} \x{\alpha_5} \x{\alpha_6} 
\x{\alpha_7} \x{\alpha_8} \x{\alpha_9} \x{\alpha_{10}}
\par
\bigskip
}

If $f:Y\to S$ is a family of smooth Enriques surfaces, then the chambers
in the N{\'e}ron-Severi groups of the geometric fibres define a local system
of sets on $S$; that is, an {\'e}tale
covering $\tS_f\to S$. 
\begin{definition}
An $E$-\Definition{marking} of a family $Y\to S$ of Enriques surfaces
is an isometry $\phi\co E_S\to
\Num(Y/S)$. A \Definition{$\sD$--polarization} of $Y\to S$ is a choice of nef classes
$[L_1],\ldots,[L_{10}]$ in $\Num(Y/S)$ for which there is a marking $\phi$ with
$\phi(\varpi_i)=[L_i]$. Equivalently, a $\sD$-polarization of $Y\to S$ is a choice
of section $\s$ of the covering $\tS_f\to S$ such that $\s(s)$ lies in the nef cone
of every geometric fibre $f^{-1}(s)$.
\end{definition}

\begin{lemma}
The datum of a $\sD$--polarization is equivalent to the datum of an $E$-marking
$\phi$ such that $\phi(\sD_0)$ lies in the ample cone of $Y$. In particular,
a $\sD$-polarization determines an $E$-marking.
\begin{proof} The Dynkin diagram $E_{10}$ has no symmetries, so the
same is true of the chamber $\sD_0$.
\end{proof}
\end{lemma}

That is, the stack $\Est_\sD$ of $\sD$-polarized Enriques surfaces is isomorphic
to the stack of smooth Enriques surfaces $Y$ with an $E$-marking $\phi$ such that
$\phi_\R(\sD_0)$ lies in the nef cone. (To be precise,
the morphisms in the latter stack are defined as follows:
if $(Y,\phi:E\to \NS(Y))$ and $(Z,\psi:E\to \NS(Z))$ are objects, so that
$\phi_\R(\sD_0)$ and $\psi_\R(\sD_0)$ lie in the relevant nef cones, then a
morphism from $(Y,\phi)$ to $(Z,\psi)$ is an isomorphism 
$f:Y\to Z$ such that $f^*\circ\psi=\phi$.)

Liedtke [L10] considers the stack $\Est_{CV,\ l.b.}$ of RDP-Enriques surfaces $X$ that
are equipped with a line bundle $\sL$, not just a polarization, such that

\noindent (i) $c_1(\sL)^2=4$,

\noindent (ii) $\sL$ is ample and

\noindent (iii) for every elliptic half-fibre $\Phi$ on the minimal resolution
$\tX$ of $X$, we have $c_1(\sL).\Phi\ge 2$.

One of his main results is that, via the morphism sending
an RDP-Enriques surface to its $\Pic^\tau$, the stack
$\Est_{CV,\ l.b.}$ is smooth over the stack of
groupschemes of order $2$, so smooth over $\Sp\Z[x,y]/(xy-2)$. 
For us it is more convenient to consider the stack $\Est_{CV}$
of RDP-Enriques surfaces with an ample Cossec-Verra polarization
rather than a line bundle. Since the forgetful morphism
$\Est_{CV,\ l.b.}\to \Est_{CV}$ is a torsor under $\Pic^\tau$,
it follows that $\Est_{CV}$ is also smooth over $\Sp\Z[x,y]/(xy-2)$.
Liedtke goes on to point out that this reduces the question
of controlling the ramification in a characteristic zero lift of a smooth
Enriques surface to controlling the ramification in a simultaneous 
resolution of RDPs (the existence of such a resolution, even in
bad characteristics, was proved by Artin [Ar74a]).

\begin{lemma} Suppose that $X$ is a smooth $\sD$-polarized
Enriques surface. Then the weight $\varpi_1$ is a Cossec-Verra polarization.
\begin{proof} By assumption, each $\varpi_i$ is nef.
The calculation of each $\varpi_i$ is simplified by first calculating
the discriminant of $\alpha_i^\perp$ in the unimodular lattice $E$.
In particular, it is easy to verify that
$$\varpi_1=4\alpha_1+9\alpha_2+14\alpha_4+7\a_3+12\a_5+10\a_6+8\a_7+6\a_8+4\a_9+2\a_{10}$$
and that
$$\varpi_{10}=2\a_1+4\a_2+3\a_3+6\a_4+5\a_5+4\a_6+3\a_7+2\a_8+\a_9.$$
So $\varpi_1^2=4$, $\varpi_{10}^2=0$ and $\varpi_1.\varpi_{10}=2$.

Now suppose that $\Phi$ is an elliptic half-fibre.
Then $\varpi_1.\Phi\ge 1$; suppose that $\varpi_1.\Phi=1$.
Note that also $\Phi.\varpi_{10}\ge 1$. Then
$\varpi_1.(2\Phi+\varpi_{10})=4$ and $(2\Phi+\varpi_{10})^2\ge 4$;
it follows from the index theorem 
that $\varpi_1=2\Phi+\varpi_{10}$, so that
$\varpi_1-\varpi_{10}$ is even. But this contradicts the expressions
above.
\end{proof}
\end{lemma}

\begin{lemma} Every Cossec--Verra polarization
$H$ is $O^+(E)$-equivalent to $\varpi_1$.
\begin{proof} $H$ is equivalent to a C-V polarization in $\sD$,
so we can assume that $H\in\sD$. Then $H=\sum_1^{10} n_i\varpi_i$
with $n_i\in\N$, and $H^2\ge\sum n_i^2\varpi_i^2$.
Since $\varpi_{10}^2=0,\varpi_9^2=2, \varpi_1^2=4$ 
and $\varpi_j^2>4$ otherwise, the result follows at once.
\end{proof}
\end{lemma}

There is an obvious forgetful morphism
$\Est_\sD\to \Est_{CV}$: given a family $X\to S$ of $\sD$-polarized
Enriques surfaces, take the corresponding family $X'\to S$ of RDP-Enriques
surfaces that arises from the polarization $\varpi_1$. The singularities
in the fibres of $X'\to S$ are configurations whose root lattices embed
in the root lattice $D_9$, since $D_9$ is what arises by deleting the
vertex $\a_1$ from the diagram $E_{10}$. These are classified by
Oshima [Os].

Artin's extension [Ar74a] of  Brieskorn's
construction [Br68] of simultaneous resolution for deformations of RDPs 
in the complex case to all contexts 
(bad characteristic, mixed characteristic) came at the cost of losing
some of the explicit information available in good characteristic, 
when the corresponding simple algebraic group tells us everything
[Br70], [Sl80], [SB01]. Now $2$ is not a good characteristic
for RDPs of type $D$ or $E$, so the next result
is not an immediate consequence of the group theory.
It depends much more upon Brieskorn's construction.

Suppose that $\sO_S\to\sO_T$ is a finite extension of mixed
characteristic DVRs, $S=\Sp\sO_S,T=\Sp\sO_T$ and $s,t$ are the
generic points of $S,T$.
\def\bt{\bar{t}}

\begin{proposition} Assume that $X\to S$ is a $1$-parameter deformation
of an RDP $(X_0,x)$ and that the absolute Galois
group of the fraction field $k(t)$
acts trivially on the Picard group of the geometric
generic fibre $X_{\bt}$ of $X_T=X\times_S T\to T$. 
Then $X_T\to T$ has a simultaneous resolution, with
no further base change.
\begin{proof} By Artin's result, there is a base change $V\to T$
and a resolution $\pi:{\widetilde{X_V}}\to X_V$. According
to [Br68], p. $257$, there is a
divisorial ideal $J$ of $X_V$ such that $\pi$ is the blow-up of $J$.
Now the restriction map from the class group $\Cl(X_V)$ of $X_V$ 
to the Picard group $\Pic (X_v)$ of the generic fibre
is an isomorphism (this holds for any $1$-parameter
deformation of a normal variety) and, by the assumption
on the Galois action, the natural homomorphism
$\Pic(X_t)\to\Pic(X_v)$ is an isomorphism.
Hence the natural homomorphism 
$\Cl(X_T)\to\Cl(X_V)$ is an isomorphism. 

According
to Theorem $6.2$ of [Sa64], this homomorphism
is defined by $I\mapsto I.\sO_{X_V}$. Therefore there is a 
divisorial ideal $I$ of $\sO_{X_T}$ such that
$I.\sO_{X_V}$ is the product of $J$ with a principal
fractional ideal. Since $X_V\to X_T$ is flat, we
have $I^n.\sO_{X_V}=I^n\otimes\sO_{X_V}$ for all $n$,
so that 
$${\widetilde{X_V}}=\Bl_J X_V=(\Bl_I X_T)\times_T V.$$
Since ${\widetilde{X_V}}\to V$ is smooth, so is
$\Bl_I X_T\to T$, and we are done.
\end{proof}
\end{proposition}

\begin{corollary} 
\part Every smooth Enriques surface $X$ over
an algebraically closed field $k$ of characteristic $2$ has
a lifting over a DVR $A$ that is a finite Galois
extension of $\W(k)[{\sqrt{2}}]$ whose
Galois group $G$ is a subgroup of the Weyl group $W(D_9)$.

\part $G$ is a semi-direct product
$G=H\rtimes C$, where $H$ is a $2$-group and $C$
is a cyclic group of odd order.

\begin{proof} Choose a $\sD$-polarization on $X$. This defines
a Cossec-Verra polarization on $X$, and so an RDP model
$X'$ that can, by Liedtke's result,
be lifted over $\W(k)[{\sqrt{2}}]$. The previous proposition
on simultaneous resolution,
combined with the fact that the stabilizer
of $\varpi_1$ in the orthogonal group $Q(E)$
is isomorphic to $W(D_9)$, 
then gives \DHrefpart{i}.
Standard results about local fields then give \DHrefpart{ii}.
\end{proof}
\end{corollary}

Recall that $W(D_9)$ is a semi-direct product
$$W(D_9)=(\Z/2)^8\rtimes\frak S_9.$$
Let $G'$ be the image of $G$ in the quotient group $\frak S_9$ of $W(D_9)$.
Then $G'=H'\rtimes C'$, with $H'$ a $2$-group and $C'$ cyclic of odd order.
MAGMA reveals that there are just $171$ conjugacy classes in $\frak S_9$
of such subgroups, and their orders divide an element of the set
$$\Sigma=\{2^63,\ 2^7,\ 2.3.5,\ 2^37,\ 2^35,\ 3^2\}.$$

\begin{corollary} \label{0.7} Every smooth Enriques surface over
an algebraically closed field $k$ of characteristic $2$ can
be lifted to a DVR of characteristic zero
whose absolute ramification index
divides $2^9N$, where $N\in\Sigma$.
\noproof
\end{corollary}
\end{section}

\begin{section}{Periods, automorphisms and moduli for singular K3 surfaces.}
\label{K3 periods}

Now move to K3 surfaces.
Let $M$ denote the lattice
$\bigoplus_1^{12}\Z.e_i+\Z.\half\sum e_i$, where $e_i.e_j=-2\delta_{ij}$ and
$N=M\oplus E(2)$. Fix a chamber $\sD$ in $E(2)\Tensor \R$ defined by the roots
in $E$. Note that $N$ is 2-elementary (that is, $2N^\vee\subset N$, where
$N^\vee$ is the dual lattice), $\sigma(N)=10$ and $x^2\equiv 0\pmod 4$
for all $x\in 2N^*$. Put $N_0=2N^\vee/2N$.
Then $\dim_{\F_2}N_0=2\sigma(N)$.

Fix, in the positive cone $N_+$ of
$N\Tensor\R$, a chamber $\sC$ of the decomposition of $N_+$ defined by the
roots. Fix also a chamber $\sD$ in the decomposition of $E(2)\Tensor\R$ defined
by the roots in $E$ that lies in the closure of $\sC$.

\begin{definition}
An \Definition{$N$--marked K3 surface over $B$} is a pair $(S\to B,\phi)$, where
$S\to B$ is a K3 surface and $\phi\co N_B\to \Pic_{S/B}$ is an embedding.
\end{definition}
\begin{definition}
$(S\to B,\phi)$ as above is \Definition{polarized} if $\phi_b(\sC)$ contains the ample
cone of $S_b$ for every geometric point $b$.
\end{definition}
\begin{definition}
The \Definition{period space} $\sM_N$ is the $k$--variety classifying maximal
totally isotropic (that is, of dimension $\sigma(N)$) subspaces $U$ of
$N_0\otimes_{\F_2}k$ such that $\dim(U\cap F(U))=\sigma(N)-1$. (Note that when
we are dealing, as we do here, with quadratic spaces, a subspace is totally isotropic
when the restriction to it of the quadratic function -- and not just the associated
scalar product -- is zero.)
\end{definition}
$\sM_N$ is smooth and projective over $\F_2$. $\sM_N\Tensor\k$ has two
irreducible components, both unirational and $9$-dimensional. 
\begin{definition}
The \Definition{period point} of an $N$--marked K3 surface $(X,\phi)$ is
$\ker(\phi\tensor 1\co N\Tensor k\to H^2_{dR}(X/k))$. (This space is a point
of $\sM_N$.)
\end{definition}
A point $s$ in $\sM_N$ determines a lattice $N(s)$
with $N\subset N(s)\subset N^*$, so a root system
$\Delta(s)$ and a decomposition of $N_+$ into chambers
[Ogus, pp. 373--374].

For the rest of this section, we shall assume that $Y$ is a unipotent Enriques surface
over an algebraically closed field of characteristic $2$
whose canonical double cover $\rho:X\to Y$ is RDP-K3. Denote by $\pi:\tX\to X$
the minimal resolution; it is well known that $\tX$ is supersingular.
Put $\rho\circ\pi=\alpha$ and let $\Gamma$ denote the exceptional locus of $\pi$.

\begin{lemma}\label{rootemb} Suppose that $(Z,P)$ is a Zariski RDP of rank $r$
with minimal resolution $f:\tZ\to Z$
and that $\xi$ is a $2$-closed vector field on $Z$ with $\xi(P)\ne 0$. Denote by
$\sF$ the $1$-foliation on $\tZ$ that is generically generated by $\xi$. Then
$\sF$ is smooth and $c_1(\sF)=-A$, where $A$ is given by the following diagram.
Moreover, if $L$ is the root lattice corresponding to $(Z,P)$, then there is an
embedding $A_1^r\inj L$ such that the sum of the positive roots is $A$.
{
\setdynkin tags=y \setdynkin slope="25
\def\x#1{$\scriptstyle#1$}
\def\present#1{\hbox{\hbox to2cm{$#1$:\hfil}\hskip1.5cm}}
\par\medskip
\present{A_1} \dynkinA 1 \x1
\par\medskip\medskip
\present{D_{2n}}
\dynkinD 9 \x2 \x2 \x4 \x4 \b\x6 \x{2n-2} \llap{\x{2n-2}} \x{n} \x{n} 
\par\medskip
\present{E_7}
\dynkinE 7 \x2 \x6 \x8 \x5 \x7 \x4 \x3
\par\medskip\bigskip
\present{E_8}
\dynkinE 8 \x4 \x{10} \x{14} \x8 \x{12} \x8 \x6 \x2 
\par
\bigskip
}
\begin{proof}
Let $\rho:Z\to (W,Q)=Z/\xi$ be the quotient, so that $W$ is a smooth germ.
Then there is a sequence $g:\tW\to W$ of $r$ blow-ups and a commutative diagram
$$\xymatrix{
{\tZ} \ar[r]^{\trho} \ar[d]_f & {\tW}\ar[d]_g\\
{Z} \ar[r]^{\rho} & {W}
}$$
where $\trho$ is the quotient by $\sF$. Then $\sF$ is smooth and
for every exceptional curve $E$ on $\tZ$, either $\sF$ is tangent to
$E$, in which case there is an isomorphism $\sF\vert_E\to T_E$ and $\trho(E)$
is a $(-1)$-curve, or $\sF$
is normal to $E$, there is an isomorphism $\sF\vert_E\to N_{E/\tZ}$
and $\trho(E)$ is a $(-4)$-curve. Since the geometry of $\tW$ is determined by $Z$,
the intersection numbers $c_1(\sF).E$ are now determined.
Define $A$ to be the cycle in the diagram above. It is then trivial to check
that $-A.E = c_1(\sF).E$ for every $E$, and so $c_1(\sF)=-A$.

Finally, $-\trho^*K_{\tW}=c_1(\sF)$.
Since $K_{\tW}$ is the sum of the $(-1)$-classes in 
each of the blow-ups comprising $g$, it follows that $-c_1(\sF)$ is the
sum of the positive roots in an embedding of $A_1^r$.
\end{proof}
\end{lemma}
\begin{lemma}\label{canpic}
Suppose that $\Pic^{\tau}(Y)=G^\vee$.
Then $\pi^*\undl{\Pic}(Y)=\undl{\Pic}(X)^G$ and if $G$ is infinitesimal
it acts trivially on $\Pic(X)$ so that in that case $\pi^*\Pic(Y)=\Pic(X)$.
\begin{proof} Since $X\to Y$ is a $G$-torsor, there is
a Hochschild-Serre spectral sequence
\begin{displaymath}
E^{ij}_2=H^i(G,H^j(X,\mul)) \Rightarrow H^{i+j}(Y,\mul).
\end{displaymath}
As $H^i(G,\mul)=0$ for $i=1,2$ (cf. [Mu70] Lemma 23.1 ii) this gives that
$\pi^*\Pic(Y)=\Pic(X)^G$. As $H^1(X,\cO_X)=0$, $\Pic(X)$ is \'etale and so $G$
acts trivially on it when it is infinitesimal.
\end{proof}
\end{lemma}

Each of the finitely many lattices $L$ between $N$ and $N^\vee=M^\vee + E(2)^\vee$
determines a locally finite collection of walls in the positive cone $C$ of
$N\otimes\R$. 
Say that $L$ is \Definition{good} if $L\cap E(2)^\vee$ contains no
roots. Take the union, over the good lattices, of the corresponding sets of
walls. This is still locally finite, and decomposes $C$ and $H$ into
not necessarily congruent regions. Such a region is a \Definition{good subchamber}
if it meets $E\otimes\R$.
\begin{lemma}\label{5.4} The chamber $\sD$ in $E\otimes\R$ 
lies in a unique good subchamber $\sR$ of $N\Tensor \R$. 
\begin{proof} This is equivalent to the statement that, for every
root $\delta$ in a good lattice, the wall
$\delta^\perp$ does not meet the interior $\sD^o$ of $\sD$.
So suppose $\delta^\perp$ meets $\sD^o$. Choose
$x$ in $\sD^o$ with $x.\delta=0$. Note that $x^2>0$.
Write $\delta=\delta_1+\delta_2$, with $\delta_1\in M^\vee$
and $\delta_2\in E(2)^\vee\cong E(\half)$.
Then $\delta_2\ne 0$ and
$\delta_1.\delta_2=x.\delta_1=x.\delta_2=0$,
so that $\delta_2^2<0$ and $\delta_1^2\le 0$.
Also, $\delta_2^2\in\Z$, from the nature of $E(\half)$,
so that either
\begin{enumerate}
\item $\delta_1^2=\delta_2^2=-1$ or

\item $\delta_1=0$ and $\delta_2^2=-2$.
\end{enumerate}
In case 1), $2\delta_2=\eta$, say, is a root in $E$, so that by construction
$\sD^o$ is disjoint from $\eta^\perp$, and so from $\delta^\perp$. In case 2)
$\delta=\delta_2$, and $L$ is not good.
\end{proof}
\end{lemma}

From now on, we fix this good subchamber $\sR$.
An \Definition{$\sR$--polarization} on a K3 surface $S\to B$
is an $N$--marking $\psi:N_B\to NS_{S/B}$
such that $\psi_b(\sR^o)$ meets the ample cone in $NS_{S_b}$
for every geometric point $b$ of $B$.

\begin{lemma}
\part If $NS_{S_b}$, regarded as a lattice between 
$N$ and $N^\vee$,
is good, then $\phi_b(\sR^o)$ meets the ample cone in $NS_{S_b}$
if and only if it is contained in it.

\part $NS_{S_b}$ is good if and only if the period point of $S_b$ lies in the
open subscheme $\sM_N^0$ (see p.~\pageref{M0def}).
\noproof
\end{lemma}

\begin{proposition}\label{polarizations}
Fix a $\sD$-polarization $\phi$ of $Y$. Then $\phi$ extends
to an $\sR$-polarization $\psi$ of $\tX$.
Moreover, $\alpha^*NS(Y)$ is a saturated sublattice of $\NS(\tX)$
and the image of the sublattice
$\bigoplus\Z e_i$ of $M$ is the sublattice of the lattice generated by the
exceptional $(-2)$--curves given by \ref{rootemb}.
\begin{proof}
There is a commutative diagram
$$\xymatrix{
{\tX} \ar[r]^{\pi} \ar[d]_{\trho} & {X}\ar[d]^{\rho}\\
{\tY} \ar[r]^{\tau} & {Y}
}$$
where $\tau$ is the composite of $12$ blow-ups and $\trho$ is the quotient by a
$1$-foliation $\sF\inj T_{\tX}$. Suppose that $C_i$ is the class in $\Pic(\tY)$
of the exceptional curve of the $i$th blow-up and put $\trho^*C_i=D_i$.
Define $\psi_1:\bigoplus\Z e_i\oplus E(2)\to \NS(\tX)$ by
$\psi_1(e_i)=D_i$ and $\psi_1(x)=\alpha^*\phi(x)$ for $x\in E$.
Next, note that $-c_1(\sF)=\sum D_i$, 
so is the sum of the positive roots
in $12\times A_1$
embedded in the sum of the root
lattices generated by $\Gamma$ as in \ref{rootemb}. 
Since $c_1(\sF)$ is even in $NS(\tX)$, by \ref{cherneven}, $\psi_1$ extends
to $\psi:N\to\NS(\tX)$. That $\psi$ is an $\sR$-polarization follows from
the fact that $\phi(\sD)$ meets the ample cone of $Y$.

As $\pi$ is a resolution of RDPs we know that
$\pi^*NS(X)$ is the orthogonal complement of the exceptional set and as such it is
saturated. Hence to show that $\alpha^*NS(Y)$ is saturated it is enough to show
that $NS(X)=\rho^*NS(Y)$ but this follows directly from \ref{canpic}.
\end{proof}
\end{proposition}

We have seen that the period point of a K3 surface $\tX$ as above
has the particular property that for any of the $N$-markings obtained
by choosing an embedding of $M$ into the saturation of the lattice generated
by the exceptional $(-2)$--curves the $E(2)$-part is 
saturated. In the following we will only use that there are no $-2$-curves in
the saturation. The final conclusion then implies that this implies that there
is no even lattice in the saturation. As a motivation for the considerations to
follow we prove this directly.
\begin{lemma}\label{E(2)-superlattices}
Any even lattice containing $E(2)$ as a proper sublattice of finite index is
generated by $E_2$ and $(-2)$-vectors.
\begin{proof}
It is enough to prove that such a lattice that contains $E(2)$ as a sublattice
of index $2$ contains a $(-2)$-vectors.
However, the group of orthogonal
transformations of $E$ maps onto the orthogonal group of $E/2E$ (cf.
[Ni80], Thm.~1.14.2) so it is enough to find one such lattice containing a
$(-2)$-vector but that is obvious.
\end{proof}
\end{lemma}
\begin{lemma}
The locus $D_1$ in $\sM_N$ whose points correspond to $N$--marked K3 surfaces
where there is a root in the saturation of $\phi(E(2))$ is a divisor.
\begin{proof}Suppose that $\phi\co N\to\NS_S$ is an $N$-marking
of $S$. Then $S$ corresponds to a point of $D_1$ if and only
if there is a root $\delta$ in the saturation $\tE$ of
$\phi(E)$. Since $\NS_S$ is 2-elementary, 
we have $\phi(E)\subset\tE\subset (\phi(E))^\vee=\frac{1}{2}\phi(E)$,
so that $2\delta\in\phi(E)$. 

Put $E_1=E+\Z.\phi^{-1}(\delta)$;
then $\disc(E_1)=-2^8$ and $E_1$ is $2$--elementary,
so is unique up to isomorphism. Moreover, there are only
finitely many copies of $E_1$ between $E$ and $E^\vee$,
so only finitely many copies of $F$ inside $E_1$.
Since $S$ is a point in $D_1$ if and only
if $\phi$ extends to $\phi_1\co M\oplus E_1\to\NS_S$,
we have identified $D_1$ with the union of
a finite number of copies of the period space
of $(M\oplus E_1)$--marked K3 surfaces.
\end{proof}
\end{lemma}
Put $\sM_N^0=\sM_N-D_1$\label{M0def}, $B_N^0=\pi^{-1}(\sM_N^0)$ and let $S^0\to
B_N^0$ be the induced universal family.

Recall that a morphism of finite type is \emph{almost proper} if it
satisfies the surjectivity part of the valuative criterion for properness.

Recall that $\varpi_1,\ldots,\varpi_{10}$ denote the fundamental
dominant weights for the labelling of the nodes of the Dynkin diagram
$E_{10}$ fixed earlier; they generate the cone $\sD$.

The next result is Ogus' global Torelli theorem for supersingular K3 surfaces,
adapted very slightly for our purposes.
\begin{theorem}\label{5.8}\label{contraction}\label{Torelli}
\part[i] There is an algebraic space $B_N$ and a
universal $\sR$--polarized $N$--marked K3 surface
$S\to B_N$. 

\part[ii] The period map $\pi:B_N\to \sM_N$ is \'etale, surjective,
of degree one and almost proper.

\part[iii] $B_N^0\to\sM_N^0$ is an isomorphism.

\part[iv] There is a contraction $S^0:=S\mid_{B_N^0}\stackrel{g}{\to}
X\stackrel{q}{\to} B_N^0$ of the induced family \map f{S^0}{B_N^0} such that for
each geometric $b\in B_N^0$ the exceptional locus of $S^0_b=S_b\to \sX_b$ is the
unique rank $12$ configuration $\Gamma_b$ of $(-2)$--curves lying in
$\phi_b(E(2))^\perp$. In particular, the configuration of RDPs on $X_b$ has
rank $12$.

\part[v] There are line bundles $\sA_1,\ldots,\sA_{10}$ on $X$
such that $c_1(g^*\sA_i)=\varpi_i$.

\begin{proof}
\DHrefpart{i} is a matter of observing that the stack of 
$\sR$--polarized $N$--marked K3 surfaces is algebraic, as follows from the
criteria of [Ar74b], and that the automorphism groupschemes of these
objects are trivial.

\DHrefpart{ii} is a very slight variant of Ogus' global Torelli
theorem, with the additional datum of an $\sR$--polarization.
His argument shows, in fact, that $\pi$ is an isomorphism over the locus
where $\phi_b(\sR^o)$ lies in the ample cone, which is \DHrefpart{iii}. 

For \DHrefpart{iv} and \DHrefpart{v}, 
note that $\phi_b(\sR)$ lies in the nef cone of $S_b$ for
every $b\in B_N^0$. So the same holds for $\phi_b(\sD)$.  
Then the line bundles $\sL_i=\phi(\varpi_i)$ on $S^0$ define the 
contraction $g:S^0\to X$;
that is, there are line bundles $\sA_1,\ldots,\sA_{10}$ on $X$
with $g^*\sA_i=\sL_i$.

On $S_b$ the exceptional locus of $g$ 
is the union of the $(-2)$--curves orthogonal to
every $L_i$. These are clearly the curves whose classes span 
the saturation in $NS_{S_b}$ of the sublattice $\phi_b(M)$.
\end{proof}
\end{theorem}

We shall be careless in distinguishing between $B_N^0$ and
$\sM_N^0$.

\begin{definition} Suppose that $Y\to S$ is a family of RDP-K3 surfaces
whose geometric fibres are of RDP index $12$. Then a
\Definition{$\sD$--polarization} of $Y\to S$ is an embedding
$\phi\co E(2)_S\to  \Pic_{Y/S}$ such that $\phi_s(\sD)$ lies
in the ample cone of $Y_s$ for all geometric $s\in S$.
\end{definition}

So, for example, the line bundles $\sA_i$ provide a $\sD$--polarization
of $q:\sX\to B_N^0$. 

Fix a geometric point $b\in B_N^0$.
\begin{lemma}\label{components}
The connected components of $\Gamma_b$ are of type
$A_1,D_{2n},E_7$ and $E_8$.
\begin{proof}
Let $\Lambda_b$ be the $\Z$--span of $\Gamma_b$. Then $\Lambda_b$ is
$2$--elementary and the result follows.
\end{proof}
\end{lemma}
Define $\Gamma_b$ to be \Definition{even} if the divisor $\phi_b(\sum e_i)$
is divisible by $2$ as a divisor, not just as a divisor class.

\begin{lemma}\label{even nodes}
If there is an even set of $r$ disjoint
$(-2)$-curves on $S$, then $r\equiv 0\pmod 4$ and $r\ge 8$.
\begin{proof}If $\sum E_i\sim 2M$, then $-2r=4M^2\equiv 0\pmod 8$.
If $r=4$, then $M^2=-2$, so that $\pm M$ is effective.
This is absurd.
\end{proof}
\end{lemma}
\begin{corollary}\label{6.4}
If $\Gamma$ is odd, then it is one of $12\times A_1$, $8\times A_1+D_4$,
$6\times A_1+D_6$, $5\times A_1+E_7$.  If $\Gamma$ is even, then it is one of
$3\times D_4$, $D_4+D_8$, $D_4+E_8$, $D_{12}$.
\begin{proof} Combine \ref{components}, \ref{even nodes} and \ref{rootemb}.
\end{proof}
\end{corollary}
\begin{lemma}
\part The set $\{e_1,\dots,e_{12}\}$ can be split into disjoint subsets, one for
each connected component of $\Gamma_b$, such that the sum of the roots in the
subset is the vector described in Lemma \ref{rootemb}.

\part $\Gamma_b$ is even if and only if its connected components are of type
$D_{4n}$ or $E_8$.
\noproof
\end{lemma}
\begin{definition}
$\sX_b$ has \Definition{even singularities} if the configuration 
$\Gamma_b$ on $S_b$ is even.
\end{definition}
\begin{lemma}
For a generic $N$--marked K3 surface $S$, the configuration $\Gamma$ given by
Proposition \ref{5.8} is $12\times A_1$.
\begin{proof}We have $\disc(\NS_S)=-2^{20}=\disc(N)$,
so $N=\NS_S$. If $\Gamma\ne 12\times A_1$, then 
$|\disc(\Z.\Gamma)|\le 2^{8-2}\times 2^2$ (equality being
achieved when $\Gamma =8\times A_1+D_4$, in which case
the $8\times A_1$ configuration is even), which gives
$|\disc(\NS_S)|\le 2^{18}$.
\end{proof}
\end{lemma}

\begin{corollary}\label{isPi-RDP}
\part
Every geometric fibre of $q\co \sX\to \sM_N^0$
has Zariski RDPs and its tangent bundle has trivial Chern classes.

\part For each singular point in a fibre, the induced deformation
is a Zariski deformation.

\part $T_{\sX/B_N^0}$ is $\cO_{B_N^0}$-cohomologically flat and
$T^1_{\sX/B_N^0}$ is a locally free $\cO_{B_N^0}$-module of rank $24$.
\begin{proof} \DHrefpart{i} and \DHrefpart{ii} follow from \ref{4.4}
and \ref{Zariski recognition}.

To verify \DHrefpart{iii} it is
enough, by \ref{Tflat}, to show that $T^1_{\sX/B_N^0}$ is flat of rank $24$. As
$B_N^0$ is reduced it is enough to show that the dimension of $T^1$ is $24$ for
each fibre. However, for a Zariski RDP this 
dimension is twice the index and the sum of all the indices is $12$.
\end{proof}
\end{corollary}
\begin{proposition}\label{3 x D4 locus}
The locus $V$ in $\sM_N^0$ corresponding to surfaces with even singularities is
irreducible and of codimension $2$.
\begin{proof}There is a unique embedding $M\inj D_4^{\oplus 3}$, and the 
closure of $V$ in $\sM_N$ is then identified with the period space
$\sM_{(D_4^{\oplus 3}\oplus E)}$.  This is $7$-dimensional.
\end{proof}
\end{proposition}

We shall refer to $V$ as the \emph{$3\times D_4$ locus}.
\begin{theorem}\label{5.13}
$T_X$ is free for every fibre $X$ of $q$.
\begin{proof}
The idea is to show that the locus $D_2$ of surfaces $X$ such that $T_X$ is not
free is a divisor in $\sM_N^0$ and is empty in codimension 1.

First, $D_2$ is not all of $\sM_N^0$, since, for example, there are canonical
double covers of Enriques surfaces giving points in $\sM_N^0-D_2$. The fact that 
$D_2$ is everywhere of codimension $1$ is \ref{nonfreecod1} and proposition \ref{5.8}.

\begin{lemma}\label{5.13.1}
The locus in $B_N^0$ (or, equivalently, $\sM_N^0$) where $\Gamma$ is of type
neither $12\times A_1$ nor $8\times A_1+D_4$ is everywhere of codimension 
$\ge 2$.
\begin{proof} If $\Gamma\ne 12\times A_1,8\times A_1+D_4$,
then $|\disc(\NS_S)|\le 2^{16}$.
\end{proof}
\end{lemma}
So it is enough to prove \ref{5.13} in the two cases 
$\Gamma= 12\times A_1,8\times A_1+D_4$. We assume that $T_X$
is not free and that $\pi\co\tX\to X$ is the minimal resolution,
so that, by \ref{4.4}, $\pi^*T_X$ has a 
destabilizing subsheaf, and deduce that there is a 
$(-2)$--curve in the smooth locus of $X$.
The existence of this curve means that the period point of $\tX$
lies in the forbidden divisor $D_1$.

\medskip

\noindent{\bf{Case (1):}} $\Gamma= 12\times A_1$. Let $E_1,\ldots,E_{12}$ be the
$(-2)$-curves in $\Gamma$, so that $\sum E_i$ is even.  Then the existence of
the corresponding square root fibration gives a $2$-integrable $1$-foliation
$\sF=\sO(B-\sum E_i)$ in $T_{\tX}$, 
where $B\ge 0$ and $B$ is disjoint from $\sum
E_i$.  By comparing $\sF$ with the destabilizing subsheaf of $\pi^*T_X$, we
see that $B>0$. We also know (Lemma \ref{cherneven}) that $c_1(B)$ is even in
$\NS_{\tX}$; say $B\sim 2C$.  Let $\rho\co \tX\to Y=\tX/\sF$ be the quotient and
$\alpha\co \tY\to Y$ the minimal resolution. Then, arguing as in \S 1, $\tY$ is
rational and $Y$ has only RDPs, of total index $r$, say. Computing $c_2(\tX)$
gives $24=-(B-\sum E_i)^2+r$, while $\rho^*K_Y\sim -B+\sum E_i$. Then Noether's
formula applied to $\tY$ gives $12=\half (B-\sum E_i)^2 +24+r$. Hence $r=0$ (and
$B^2=0$), so that $Y$ is smooth. Put $\rho(E_i)=F_i$, so that $F_i$ is a
$(-1)$-curve on $Y$ and $\rho^*F_i=E_i$.

There is a commutative diagram
$$\xymatrix{
{\tX} \ar[r]^{\rho} \ar[d]_{\pi} & {Y}\ar[d]^{\tau}\\
{X} \ar[r]^{\rho_1} & {Y_1}
}$$ 
where $\tau\co Y\to Y_1$ contracts the curves $F_i$ to points $P_i$.  Note that
any curve $C$ on $Y$ with $C^2<0$ has $C^2=-1$ or $C^2=-4$, according as
$\rho^*C=D$ or $\rho^*C=2D$.

Suppose that $l_1$ is a $(-1)$-curve on $Y_1$, with strict transform $l$ on
$Y$. If $l^2=-4$, then $\rho^*l=2m$ and $\tau^*l_1=l+\sum_1^3F_i$, so that
$l.\sum_1^{12}F_i=3$.  Hence $m.\sum E_i=3$.  But $\sum E_i$ is even, so that
$l^2=-1$. Then $\rho^*l$ is the $(-2)$-curve we sought.

\medskip

\noindent{\bf{Case (2):}} $\Gamma =8\times A_1+D_4$.
Then, if $E_1,\ldots,E_8$ form the $8\times A_1$
configuration, $\sum E_i$ is even.
Suppose that $\sum G_i$ is the $D_4$-configuration,
with $G_1$ the central curve.
As in (1), there is a $2$-integrable 1-foliation
$\sF=\sO(B-\sum E_i)$ with $B\ge 0$ and 
$B\cap\sum E_i=\emptyset$, and comparison
with a destabilizing subsequence of $\pi^*T_X$
shows that $B.\pi^*H>0$ for ample $H\in\Pic X$,
so that $B>0$. A calculation similar to that in (1)
shows that $B^2=-8$ and $\tX/\sF=Y$ is smooth.
Moreover, $B$ is even; say $B\sim 2C$. Then $C^2=-2$
and $C.\pi^*H>0$, so that, by Riemann-Roch,
$C$ is effective. Moreover, $\Supp C$ is not
contained in $\sum G_i$.

Consider the commutative diagram
$$\xymatrix{
{\tX} \ar[r]^{\rho} \ar[d]_{\phi} & {Y} \ar[r]^{\sigma} \ar[d]_{\psi} & 
{\tX^{(1)}} \ar[d]^{\phi^{(1)}}\\
{X_2} \ar[r]^{\rho_2} & {Y_2} \ar[r]^{\sigma_2} & {X_2^{(1)}}
}$$
where $\phi$ contracts the $E_i$ and $\psi$
contracts their images $F_i=\rho(E_i)$ to points $P_i$. 
So $Y_2$ is a smooth rational surface
and $\rho_2^*c_1(Y_2)\sim B\sim 2C$.
So $-K_{Y_2}\sim \sigma_2^*C^{(1)}$ and is thus effective.
Let $H_i$ denote both $\rho(G_i)$ and its image in $Y_2$.

\begin{lemma}
Suppose that $D$ is an effective
anti-canonical divisor on a smooth rational surface.

\part $D$ is connected.

\part Either $D$ is reduced and $h^1(\sO_D)=1$ or $D_{red}$ is a tree of
$\P^1$s.

\part If $E\le D$ with $h^1(\sO_E)=1$ and $E$ is reduced, then $E=D$.
\begin{proof}This is easy and well known.
\end{proof}
\end{lemma}
As in Case (1), no $(-1)$-curve on $Y_2$ passes through
any $P_i$.

\begin{lemma}
$C.G_i=\pm 1$.
\begin{proof}$C.G_i=1$ if and only if $\sF$ is tangent to $G_i$
and $C.G_i\le -1$ otherwise. If $C.G_i\le -2$, then the non-zero
map $\sF|_{G_i}\to\sN_{G_i/Y_2}$ is not an isomorphism,
so that $H_i$ is singular. Also, $G_i\le C$, so that
there exists $H\in|-K_{Y_2}|$ with $H\ge H_i$. Then $H=H_i$,
so that $C$ is supported on $\sum G_i$.
\end{proof}
\end{lemma}
\begin{lemma}
$C.(G_1+G_i)=0$ for $i\ne 1$.
\begin{proof}Suppose that $C.G_1=C.G_2=-1$. Then
$G_1,G_2\le C$ and $H_1.H_2=2$. So there exists
$H\in|-K_{Y_2}|$ with $H\ge H_1+H_2$, so that
$H=H_1+H_2$ and again $C$ is supported on $\sum G_i$.
\end{proof}
\end{lemma}
So there are two possibilities for the configuration
of intersection numbers $C.G_i$.

{\bf{Subcase (2.1):}} $C.G_1=-1$ and $C.G_i=1$ for $i\ne 1$.
Then $\sum H_i$ contracts to a smooth point, say via
$Y_2\to Y_1$. Suppose that $l_1$ is a $(-1)$-curve
on $Y_1$ and $l$ its strict transform on $Y$. If $l^2=-1$,
then $\rho^*l$ is the $(-2)$-curve we sought. So
suppose that $l^2=-4$. Then, as in Case (1), $l.\sum F_i$
is even. It is then easy to verify that $l$ meets
just two of the $F_i$, say $F_1,F_2$, and that $l.H_1=1$,
$l.H_i=0$ for $i\ne 1$. Moreover, $\rho^*H_1=2G_1$
and $\rho^*l=2m$, so that $4m.G_1=2$, which is absurd.

{\bf{Subcase (2.2):}} $C.G_1=1$ and $C.G_i=-1$ for $i\ne 1$.
Then $G_i\le C$ for $i\ne 1$ and $H_1$ is a $(-1)$-curve.
Let $Y_1\to Y_1'$ be the contraction of $H_1$ and $H_i'$ the image
of $H_i$ in $Y_1'$. Then there exists $H\in|-K_{Y_2}|$ with 
$H\ge\sum_2^4 H_i'$, so that $H=\sum_2^4 H_i'$ and
$C$ is supported on $\sum G_i$. This contradiction completes the
proof of \ref{5.13}.
\end{proof}
\end{theorem}
\begin{corollary}\label{5.14}
The formation of $T_{\sX/\sM_N^0}$ commutes with base change and
$q_*T_{\sX/\sM_N^0}$ is locally free of rank $2$, its formation
commutee with base change and
the map $q^*q_*T_{\sX/\sM_N^0} \to T_{\sX/\sM_N^0}$ is an isomorphism.
\begin{proof}
The first part is \ref{isPi-RDP} and the rest \ref{nonfree}.
\end{proof}
\end{corollary}
\end{section}

\begin{section}{Some automorphism group schemes.}
\label{automorphisms}
Here we make some calculations concerning
automorphism group schemes of certain unipotent Enriques surfaces and
RDP-K3 surfaces. These will be useful in the analysis in \S \ref{stacks of surfaces}
of the relationships between various stacks of surfaces. 

We fix the following notation to be used in this section.

\begin{notation} 
The base is an algebraically closed field $k$ of characteristic $2$.
For a group scheme $H$, denote by ${}_nH$ the kernel
of $\Frob^{(n)}:H\to H^{(n)}$ (the height $n$ part of $H$). Denote by
$\Aff$ the group of affine transformations
of the affine line. $(S,\psi)$ will be
an $\sR$-polarized $N$--marked K3 surface 
whose period point lies in
$\sM_N^0$ and $\pi\co S\to S_1$ the contraction of the configuration
$\Gamma$ given by Proposition \ref{5.8}. Then $T_{S_1}$ is free, so that
$H^0(T_{S_1})$ is a $2$-dimensional
$2$-Lie algebra $\frak g$; it is the Lie algebra of 
$G=G_{S_1}:={}_1(\Aut^0_{S_1})$. Recall that there is a bijection between
connected height $1$ groupschemes $H$ over $k$ and $2$-Lie algebras $\frak h$
over $k$; if $\frak h=\Lie(H)$, then we write $H=\exp(\frak h)$.
\end{notation}

\begin{proposition}\label{6.3}
\part Suppose that $\xi\in\frak g$ is $2$--closed
and that $Y=S_1/\xi$ is smooth. Then $Y$ is an Enriques
surface and $\rho\co S_1\to Y$ is the canonical double cover.

\part There is a unique $\sD$-polarization $\phi:E\to \NS(Y)$
such that $\rho^*\phi:E(2)\to\NS(S)$ extends to 
$\psi$ and $\phi(\sD)$ lies in $\sR$.

\begin{proof} \DHrefpart{i} $V$ generates a copy of $\sO_{S_1}$ in $T_{S_1}$,
so that $\rho^*c_1(Y)\sim 0$. Then $2K_Y\sim 0$ and $b_2(Y)=10$,
so that $Y$ is Enriques.
Since $\Pic^\tau_{S_1}=0$, $\rho$ factors through the
canonical double cover $Z$ of $Y$. So $S_1$ is the normalization
of $Z$. Since $\omega_Z$ and $\omega_{S_1}$ are trivial
the conductor is empty and $S_1\to Z$ is an isomorphism.

\DHrefpart{ii} Since $\psi(E(2))$ is orthogonal to the exceptional curves
of $\pi$ it is contained in $\rho^*(\NS(Y))$. Then $\psi\vert_{E(2)}$
is an isometry $\phi:E\to\NS(Y)$. Moreover,
$\phi(\sD)$ lies in $\psi(\sR)$, which meets the ample cone of $S$,
so that $\phi(\sD)$ meets the ample cone of $Y$ and therefore lies in it.
\end{proof}
\end{proposition}

\begin{proposition} Suppose that $(Y,\psi)$ is a unipotent $\sD$-polarized
Enriques surface and that its canonical double cover $X$ is RDP-K3 with minimal
resolution $\pi:\tX\to X$. Then $X\cong S_1$ for some $S_1$
as above.
\begin{proof} This is just \ref{polarizations}
\end{proof}
\end{proposition}

We know that $Y\cong X/\sF$ for some smooth $1$-foliation $\sF$,
that $X\cong S_1$ for some $S_1$ as described, and that $c_1(\sF)$ is numerically
trivial in codimension one. It follows at once that $\sF$ is globally generated
by a $2$-closed vector field.

\begin{corollary} \label{2-closed}
Every vector in $\frak g$ is $2$-closed.
\begin{proof} We know that an Enriques surface has at least $10$ moduli (that is,
every component of the coarse moduli space is at least $10$-dimensional),
while a supersingular K3 surface has $9$ moduli. Hence there is at least one
degree of freedom in the choice of $2$-closed line in $H^0(T_{S_1})$.
\end{proof}
\end{corollary}

\begin{lemma}\label{6.1}
$G$ is either ${}_1\Aff$ or $\alpha_2\times\alpha_2$.

\begin{proof}
We shall solve the equivalent problem of determining $\frak g$.

Suppose first that $G$ is not commutative.
Let $H$ denote the derived subgroup scheme $[G,G]$. Then
$\frak g$ has a basis $\{\xi,\eta\}$ such that $\xi$ generates $\Lie H$
and $[\xi,\eta]=\xi$. Rescale $\xi$ so that $\xi^2=\xi$ or $0$. Say
$\eta^2=b\eta$. Then, since $(\lambda\xi+\eta)^2$ is proportional to
$\lambda\xi+\eta$ for all $\lambda\in k$, it follows at once that
$\xi^2 =0$. Then rescale $\eta$ so that $\eta^2=\eta$ or $0$.
In the first case $G\cong {}_1\Aff$ and in the second the proportionality
of $(\lambda\xi+\eta)^2$ and $\lambda\xi+\eta$ for all $\lambda\in k$
gives a contradiction.

Now suppose that $G$ is commutative. Then pick a basis $\{\xi,\eta\}$
and consider $(\lambda\xi+\eta)^2$, as before. It is immediate that
$\xi^2=\eta^2=0$, so that $G\cong \alpha_2\times\alpha_2$.
\end{proof}
\end{lemma} 

\begin{remark} The list of all connected height 1 group schemes of
order 4 is rather longer. It includes a non-commutative 
family parametrized by $\P^1$,
of which ${}_1\Aff$ is one member, whose general member has Lie algebra
generated by $\xi,\eta$, with $\xi^2=0$, $[\xi,\eta]=\xi$ and 
$\eta^2=a\xi +\eta$.
\end{remark}

\begin{lemma}\label{A_1 behaviour}\label{alpha free}
If $0\ne \xi\in H^0(T_{S_1})$ and $\xi^2=0$, then
$\xi$ does not vanish at any $A_1$ singularity of $S_1$.

\begin{proof} At an $A_1$ singularity $P$, $S_1$ is locally given by
an equation $xy+z^2=0$. Then $\{D_1=x\pd_x+y\pd_y,D_2=\pd_z\}$
is a local basis of $T_{S_1}$. Since $\xi$ is an element of a global
basis, we can write $\xi =fD_1+gD_2$, where $f,g\in\sO_{S_1,P}$ and at
least one of them is a unit. Assume that $\xi(P)=0$; then $g\in\frak m_P$,
so that $f$ is a unit. 
Since $D_1^2=D_1$ and $D_2^2=0$, taking the coefficient of $D_1$ in
the equation $\xi^2=0$ gives
\begin{displaymath}
fD_1(f) + f + gD_2(f)=0.
\end{displaymath}
However, $D_1(f),g\in\frak m_P$, which is absurd.
\end{proof}
\end{lemma}

\begin{theorem}\label{odd implies Aff}
If $\Gamma$ is even, then $G\cong\alpha_2\times\alpha_2$.
Otherwise, $G\cong {}_1\Aff$.
\begin{proof} Since $\frak g$ is $2$-dimensional, there is, for any singular
point $P$ on $S_1$, a vector field $\xi\in\frak g$ with $\xi(P)=0$. 
If $\Gamma$ is odd, then $S_1$ has an $A_1$ singularity, and the result
follows from \ref{A_1 behaviour}. 

For the converse, suppose that $G\cong {}_1\Aff$.
Denote the radical of $\frak g$ by $\frak a$
and put $Y_1^{(1)}=S_1/\frak a$. Then $Y\to S$ is a
principal $\mu_2$-bundle over $S_1-\Sing(S_1)=S_0$, say.
The exact sequence 
\begin{displaymath}
0\to\oplus_{i\in I}\Z.E_i\to\Pic(S)\to\Pic(S_0)\to 0
\end{displaymath}
shows that then there is a subset $J\subset I$
such that $\sum_JE_j$ is even. Since $\# J=8$ or $12$, by \ref{even nodes},
it is easy to check that this is impossible for
each of the even configurations listed in Corollary
\ref{6.4}. So $\Gamma$ is odd.
\end{proof}
\end{theorem}

\begin{proposition}\label{multiplicative}\label{7.7}
A multiplicative $\alpha_2$-Enriques surface $X$ is $12A_1$-Enriques.
\begin{proof}
The multiplicative vector field corresponds to an action of $\mu_2$ on 
$X$. We know by \ref{smoothfix} that the fixed point scheme of such an
action is smooth. We first want to show that this fixed point scheme is
non-reduced at all of its codimension $1$ points. We choose a $\Pic$-rigidification 
of $X$ and employ the notations of Definition \ref{multelts}. As we have
$D(f)\check\beta_2=\eta \wedge df$ we see that the zero-scheme of $D$ equals
that of $\eta$. On the other hand, $\eta$ comes from $H^0(X,B_1)$ so is locally
of the form $df$. However, in characteristic $2$ the divisorial part of the zero 
scheme of such a form is twice a divisor which can be seen either by make a
calculation in coordinates at a point where the zero set of $f$ is smooth or by
using [Ek88], Prop.~1.11. Hence the multiplicity of the codimension
$1$-part of the zero set of $D$ is even and as it is also smooth it must be empty.

We now note that if $\eta$ locally can be written $df$ then, over the same open
subset, the canonical double cover can be written $\Sp\sO_X[z]/(z^2-f)$ so that
to begin with the canonical double cover has only isolated
singularities. Furthermore, under a singularity, $df$ has a smooth zero scheme
which means that it must have the form $\lambda+xy$ in suitable local
coordinates and hence the double cover has an $A_1$-singularity and $X$ is
K3-Enriques.
\end{proof}
\end{proposition}

Define $S_1$ to be of type $12\times A_1$ if it has $12$ singularities, 
each of type $A_1$.

\begin{theorem}\label{height 1}
$G=\Aut^0_{S_1}$ if $S_1$ is of type $12\times A_1$.
\begin{remark} In Corollary \ref{7.9}
this will be extended
to all surfaces $S_1$ arising from a K3 surface $S$ whose period
point lies in $\sM_N^0$.
\end{remark}

\begin{proof} Put $A=\Aut^0_{S_1}$ and
assume the result false. So
$G={}_1A\ne A$. Then $Frob :A\to A^{(1)}$ kills $G$ but not all of $A$.
Then there is an order $2$ subgroupscheme in the image; define
$B$ to be its inverse image. 
Then $B$ is an order $8$ subgroupscheme of $A$ with $G={}_1B$ and $G$ normal in $B$,
since $G$ is normal in $A$.

By \ref{odd implies Aff},
$G={}_1\Aff$. Then the unique copy $\alpha$ of $\alpha_2$ in $G$ is normal
in $B$. Conjugation gives a homomorphism 
$\phi:B\to \Hom _{gp-sch}(\alpha,\alpha)\cong \mul$.
Say $\ker\phi =H$. Then $H$ contains $\alpha$ but not $G$. So either 
$H=\alpha$ or $H$ has order $4$. Let $\rho:S_1\to S_1/\alpha =X$
denote the quotient map;
then $X$ is smooth, by \ref{alpha free}, and so is K3-Enriques
with $B/\alpha\subset\Aut^0_X$.
 
Suppose first that $H$ has order $4$. 
Then $G,H$ are distinct normal order $2$ subgroupschemes of $B$
and $G\cap H=\alpha$. Choose generators $\xi$, resp. $\eta$, of
$\Lie(G/\alpha)$, resp. $\Lie(H/\alpha)$. Then $\ad(\xi)(\eta)\in k\eta$
and $\ad(\eta)(\xi)\in k \xi$. So $[\xi,\eta]=0$, and we deduce
that $B/\alpha =G/\alpha \times H/\alpha$. Then $\xi$, $\eta$
give linearly independent elements of $H^0(X,T_X)$,
which is impossible. 

Hence $H=\alpha$ and there is an exact sequence
$$ 1 \to \alpha \to B \to \mu_4 \to 1.$$
In particular, $\mu_4\inj\Aut^0_X$.
Moreover, the sequence splits,
as follows. The adjoint representation
gives an embedding $B\inj GL_2$, and $B$ preserves a line, namely
$\Lie(\alpha)$. So $B$ lies in an upper triangular group and the intersection
of $B$ with the diagonal subgroup is a copy $\mu$ of $\mu_4$. 
We consider the action of $\mu$ on $S_1$. Denote by $\mu_2$
the unique copy of $\mu_2$ in $\mu$.
 
\begin{lemma}
If $\mu_2\inj\mu$ acts freely at a point $s\in S_1$, then so does $\mu$.
\begin{proof} Since $\mu_2$ is the unique proper subgroup scheme of $\mu$,
$\Stab(s)$ is trivial if it does not contain $\mu_2$.
\end{proof}
\end{lemma}

\begin{lemma} $B$ acts freely on $S_1\setminus \Sing (S_1)$.
\begin{proof} $G$ acts freely on $S_1\setminus \Sing (S_1)$ and meets every
non-trivial subgroup scheme of $B$ non-trivially.
\end{proof}
\end{lemma}

Now let $\rho:S_1\to Y$ denote the quotient by $\mu_2$ and $\sigma:Y\to Z$
the quotient by $\mu/\mu_2$. Put $\rho(P_i)=Q_i$ and $\sigma(Q_i)=R_i$.


\begin{lemma} If $\mu_2$ acts freely at $P_i$, then 
$Y$ is smooth at $Q_i$ and $Z$ is smooth at $R_i$.
\begin{proof}
Same arguments as before.
\end{proof}
\end{lemma}

\begin{notation} We let $[a_1,\ldots,a_n]$ denote a chain of $n$ transverse
smooth rational curves with self-intersections successively $a_1,\ldots,a_n$.
\end{notation}

\begin{lemma} If $\mu_2$ fixes $P_i$, then the minimal resolution of $(Y,Q_i)$ is 
either $[-2,-2,-2]$, in which case $(Y,Q_i)$ is an $A_3$ singularity, or $[-4]$.
\begin{proof}
Since $\mu_2$ fixes $P_i$ it acts on its blow-up, which is the minimal resolution.
Let $E$ denote the exceptional curve. Either $\mu_2$ acts trivially on $E$,
giving the first case, or it acts with just two fixed points,
which gives the second.
\end{proof}
\end{lemma}

Hence $\Sing Y$ consists of $a$ of the $[-4]$ singularities and $b$ of the
$A_3$ singularities, where $a+b\le 12$. Let $f:\tY\to Y$ be the minimal
resolution. Recall that $K_Y$ is numerically trivial and that 
the irregularity $q(\tY)=0$.

Suppose that $a=0$. Then $\tY$ is $K3$ or Enriques and $c_2(\tY)=12+3b$.
So $b=0$ or $4$. If $b=0$, then $Y$ is a smooth $\Z/2$-Enriques surface
whose canonical double cover is $S_1$. So $H^0(T_Y)=0$. However, $\mu_4/\mu_2$
acts on $Y$, and so $b=4$. Then (Noether)
$\chi(\sO_Y)=2$, while $Y$ can be deformed,
by varying a vector field on $S_1$, to a smooth Enriques surface with
$\chi(\sO)=1$. (Note that this relies upon the linear
reductivity of $\mu_2$, which ensures that taking invariants commutes with
specialization.)

Hence $a>0$. Then $\tY$ is rational, $K_{\tY}^2 =-a$ and $c_2(\tY)=12+a+3b$.
So $b=0$.

\begin{lemma} The action of $\mu_4/\mu_2$ on $Y$ lifts to $\tY$.
\begin{proof} Pick a singular point $Q_i$ on $Y$ and put $f^{-1}(Q_i)=E$.
The action of $\mu_4/\mu_2$ on $Y$ defines a $2$-closed $1$-foliation 
$\sF$ on $\tY$, with $\sF\cong\sO(aE)$ near $E$. Comparing the induced
map $\sF\vert_E\to T_{\tY}\vert_E$ with the normal bundle sequence of
the embedding $E\inj\tY$ shows that $a\ge 0$.
\end{proof}
\end{lemma}

Recall that $\mu_4/\mu_2$ acts freely on $Y\setminus\Sing Y$. The
fixed locus of the action on $\tY$ is smooth, so near an exceptional curve
$E$ consists
either of $E$ or of two points on $E$. It is then easy to see that
if $\tZ\to Z=Y/(\mu_4/\mu_2)$ is the minimal resolution, the exceptional 
locus in $\tZ$ is $c\times [-2,-3,-2] + d\times [-8]$ with $c+d=a$. 
Moreover, $K_Z$ is numerically trivial, so that $\tZ$ is rational,
$K_{\tZ}^2= -c/2 -9d/2$ and $c_2(\tZ)=12+3c+d$. Then $5c=7d$,
so that $c=7$ and $d=5$. Then $r= a=12$, so that $\Sing S_1 = 12\times A_1$
and the $\mu_2$-action fixes every singular point. Then the $\mu_2$-action
lifts to the minimal resolution of $S_1$, which contradicts the Rudakov-Shafarevich
theorem.
\end{proof}
\end{theorem}

\begin{proposition}\label{kernel} Suppose that $\Xi$ is a line in $\frak g$. 
Put $Y=S_1/\Xi$.
Then the kernel $A$ of the natural homomorphism
$\Aut_{(X,\Xi)}\to\Aut_Y$ is $\exp(\Xi)$.

\begin{proof} Put $\exp(\Xi)=H$. Certainly $H\subset A$. 

Suppose first that $g\in A(k)$. Since $S_1$ and $Y$ have the same underlying topological 
space and $g$ maps to the identity in $\Aut_Y(k)$, it follows that $g=1$. 
So $A$ is connected. 

Next, suppose that $\eta\in\Lie(A)\setminus\Xi$. Then $\eta$ is a vector field on $S_1$
that commutes with $\Xi$, and so descends to a vector field $\bar\eta$ on $Y$.
By assumption, $\bar\eta=0$. Since $\Xi,\eta$ generate $T_{S_1}$, it follows that
$\sO_Y$ is invariant under $T_{S_1}$, so that $S_1\to Y$ factors through the Frobenius
$S_1\to S_1^{(1)}$, which is absurd.
Hence $\Lie(A)=\Xi$, so that $A$ has height $\ge 2$.

Since $S_1\to Y$ is a torsor under the $Y$-groupscheme $H\times Y\to Y$
and $A\times Y$ is a subgroupscheme of the $Y$-groupscheme $\Aut_{X/Y}$,
it follows that $A$ embeds into the semi-direct product $H\rtimes \Aut_{gpsch}(H)$.
Since $\Aut_{gpsch}(\mu_2)=1$, it follows that $H=\alpha_2$ and
$A$ embeds into $\alpha_2\rtimes\mul =B$, say. However, $A$ contains $H$
and has $1$-dimensional Lie algebra, so this is impossible.
\end{proof}
\end{proposition} 

\end{section}

\begin{section}{From surfaces to periods}\label{stacks of surfaces}
\label{period mapping}
In this section we shall look more closely at the passage from Enriques surfaces
to periods. This restricts us to K3-Enriques surfaces. It is clear from
what we now know that the coarse moduli space of K3-Enriques
surfaces whose $\Pic^\tau$ is unipotent is,
modulo finite (reduced) groups of automorphisms,
an open piece of something that is fibred over something that is equivalent,
under a radicial map, to an explicit open piece of the
period space of K3 surfaces, and the reduced geometric
fibres are curves of \emph{geometric}
genus zero. Obviously, this description is very crude;
for example, the geometric generic fibre just referred to might,
\emph{a priori}, fail to be reduced. 
To get a more satisfactory description, we must make a
more careful examination of the stacks involved and the relevant morphisms
between them. It turns out to be easier to go in the other direction and examine
the passage
from periods to Enriques surfaces.
We need a definition.

\begin{definition}
A family of RDP-K3 surfaces is a \Definition{Zariski family} if
the singularities in each geometric fibre are Zariski RDPs and
the induced deformation of each of these singularities
is a Zariski deformation.
\end{definition}

Here is a list of most of the stacks that we shall investigate. 
They will be algebraic and of finite type
over the base $B=\Spec \F_2$.

\noindent (1) The stack $\Est_{uni}$ of
unipotent $\sD$--polarized Enriques surfaces. The open
substack corresponding to the K3-Enriques surfaces is $\Est_{K3,uni}$.



\noindent (2)
The stack $\Kst$ of pairs
$(f\co Y\to S, \psi)$ such that $f$ is a Zariski family of
RDP--K3 surfaces of RDP rank 12,
the sheaf $f_*T_{Y/S}$ of Lie algebras is locally free of rank 2, 
the natural homomorphism
$f^*f_*T_{Y/S}\to T_{Y/S}$ is an isomorphism, every sub-bundle of $f_*T_{Y/S}$
is $2$-closed and $\psi$
is an embedding $E(2)_S\to\Pic_{Y/S}$ such that $\psi_s(\sD)$
lies in the ample cone of $Y_s$ for all geometric points
$s$ of $S$. (These we abbreviate to ``$\sD$-polarized 12-nodal RDP K3 surfaces 
with trivial tangent bundle'' or just ``12-nodal RDP K3 surface''.)

\noindent (3) The open substack $\Kst^{odd}$ of $\Kst$ where the geometric fibres 
have odd singularities.

\noindent (4) The relative stabilizer groupscheme $\sG=\sG_{\Kst/B}\to\Kst$ and its
height $1$ subrelative groupscheme ${}_1\sG\to\Kst$. Since for every
morphism $s:S=\Spec\sO\to \Kst$ of a local scheme to $\Kst$ the family
$s:X\to S$ gives a free $2$-closed $\sO$-Lie algebra of rank $2$,
the Lie algebra of $\sG\to\Kst$ or ${}_1\sG\to B$
is a rank 2 locally free sheaf of $2$-Lie algebras $\frak G\to\Kst$,

\noindent (5) The bundle $\pi:\sP=\P(\frak G^\vee)\to\Kst$ ([LMB00], p. 137). 
We identify $\sP$ with
the stack of 12-nodal RDP K3 surfaces with a choice of a
line of vector fields; recall that every such line is $2$-closed. 
So an $S$-point of $\sP$ is a pair $(X\to S,\Xi)$, where $X\to S$
is an $S$-point of $\sK$ and $\Xi$ is a line of vector fields.
Note that there is a tautological finite flat closed subrelative groupscheme
$\sH\to \sP$ of $\pi^*({}_1\sG)\to\sP$ that is of order $2$; its Lie algebra
is the kernel of the natural homomorphism $\pi^*\frak G\to\sO(1)$ and
over an $S$-point $(X\to S,\Xi)$ of $\sP$ this Lie algebra also coincides 
with $\Xi$.

\noindent (6) The open substack $\sU$ of $\sP$, the complement of the closed
substack where $\Xi$ vanishes at a singularity. Of course,
$\pi:\sU\to\Kst$ is smooth and representable.

\noindent (7) The closed substack $\sP_{nilpt}$ of $\sP$ where $\Xi$
is nilpotent, and its complement $\sP_0$. (That $\sP_{nilpt}$ is a closed substack
follows from the formula $(a+b)^{[p]}=a^{[p]}+b^{[p]} +\sum s_i(a,b)$
that holds in a $p$-Lie algebra, where the $s_i$ are universal Lie polynomials.)
Note that $\sP_{nilpt}$ is a quasi--section
of $\sP\to \Kst$, in the sense that it restricts to a section over $\Kst^{odd}$.


\begin{notation} Recall that every finite locally free bundle
$\frak H\to S$ of $p$-Lie algebras is the Lie algebra
of a unique connected height $1$ finite flat groupscheme $H\to S$; 
we shall write $H=\exp(\frak H)$.
\end{notation}

\begin{theorem}\label{phi is smooth} 
There is 
an algebraic gerbe $\phi:\sU\to \Est_{K3,uni}$ that is an extension by $\sH$.
In particular, $\phi$ is smooth. 
\begin{proof} According to \ref{key}, dividing a Zariski RDP-K3 surface
with trivial tangent sheaf by a nowhere vanishing $2$-closed vector field
gives an Enriques surface. Then
the morphism $\phi$ is given at the level of $S$-points by associating to the
point $x=(X\to S,\Xi)$ in $\sU(S)$ the point $X/\Xi\to S$ in
$\Est_{K3,uni}(S)$. Note that because $\xi$ is nowhere zero, the groupscheme
$\sH_S\to S$ acts
freely on $X\to S$. So taking the quotient
by $\Xi$ commutes with base change. 

Next, recall [Ra70], 6.2.1, that for an object \map fYS of $\Est(S)$ and a
finite flat commutative $S$-group scheme $A$ there is a natural isomorphism $R^1f_*A
\riso \Hom(A^\vee,\Picf(Y/S))$, where $A^\vee$ is the Cartier dual of
$A$. To begin with, this means that locally (in the fppf topology) there is a
$\Picf^\tau(Y/S)^\vee$-torsor which corresponds to the identity element. Over each
fibre this then is the double canonical cover. If \map fYS is an object
of $\Est_{K3,uni}$ this gives a point in 
$\sU(S)$, which shows local essential surjectivity. 

Now note that a non-trivial torsor over an Enriques surface
$Y$ (over an algebraically closed field $k$) 
under a group scheme $A/k$ of order $2$ 
is isomorphic to the canonical double cover 
of $Y$ through an isomorphism that identifies $A^\vee$ with 
$\Picf^\tau(Y)$. Indeed, the torsor corresponds to a homomorphism 
$\psi:A^\vee \to \Picf(Y)$, which is non-trivial as the torsor is, 
and so $\psi$ induces
an isomorphism $\Psi:A^\vee \riso \Picf^\tau(Y)$ as both are group schemes of 
order $2$. Then the torsor is,
again by Raynaud's isomorphism, isomorphic to the one obtained from the
canonical double cover by $\Psi$.

Assume now
that we have two objects $(X,\Xi)$ and $(X',\Xi')$ of $\sU(S)$ and an isomorphism
between the associated objects of $\Est_{K3,uni}(S)$. Composing this isomorphism 
with the quotient map gives us two non-trivial torsors, over the same base \map fYS, 
under $A=\exp(\Xi)$ and $A'=\exp(\Xi')$. 
As has just been observed they are both isomorphic, fibre by fibre, to
the canonical double cover. This shows that $A$ and $A'$ are
isomorphic to the dual of $\Picf^\tau(Y/S)$, Then, using these isomorphisms, we
get isomorphisms between $A$ and $A'$ such that the two elements of $R^1f_*A$
to which the two torsors give rise
are the same, since they are the same under the
identification of $R^1f_*A$ with $\Hom(A^\vee,\Pic^\tau(Y/S))$. This shows
that the two torsors become isomorphic after a flat base extension and hence
$\phi:\sU \to \Est_{K3,uni}$ induces surjective maps on $\Hom$-sheaves
and is therefore an algebraic gerbe. 

It remains to show that $\sG_\phi=\sH$. Certainly, $\sH\subset\sG_\phi$.
Since $\sK$ is a stack of non-ruled polarized surfaces with only RDPs,
the relative groupscheme $\sG_{\sK/B}\to\sK$ is finite, by the Matsusaka--Mumford theorem,
so the closed subrelative groupscheme $\sG_\phi\to\sK$ is finite. Since $\sH\to\sK$
is finite and flat of order $2$, it is enough to show that $\sG_\phi$
has order $2$ over every geometric point of $\sK$.
But this is the content of \ref{kernel}. Finally, the smoothness
of $\phi$ follows from \ref{gerbes with flat stabilizer are smooth}.
\end{proof}
\end{theorem}

We now focus on the diagram
$$\xymatrix{
{\sU} \ar[r]^{\phi} \ar[d]_{\pi} & {\Est_{K3,uni}}\\
{\Kst.}
}$$
The morphisms are smooth and $\pi$ is representable. The next step is to
analyze $\Kst$ in terms of Ogus' periods.

For this, return to the diagram
$$\xymatrix{
{S^0} \ar[r]^g \ar[d] & {\sX}\ar[d]^q\\
{B_N^0}\ar[r]_{\cong} & {\sM_N^0}
}$$
of \ref{5.8}. The symmetric group $\frak S_{12}$ acts on $M$, 
and so on $N$,
by permuting the $e_i$, and is the stabilizer 
in the orthogonal group $O(M)$ of a suitable chamber $\sC$.

\begin{corollary}
The family $q\co \sX\to\sM_N^0$  is Zariski.
\begin{proof}
This is an infinitesimal problem so we may assume the base is local
infinitesimal. Then by \ref{5.14} every vector field on the closed fibre can be
lifted and as there is one that acts freely at a given (in fact on all)
singular point(s) we have a Zariski deformation.
\end{proof}
\end{corollary}

\begin{lemma}
The action of $\frak S_{12}$ interchanges the two geometric components
of $\sM_N$.
\begin{proof} Trivial.
\end{proof}
\end{lemma}

\begin{proposition}\label{local to global for K3s}
Suppose that $X/k$, with $k=\bar k$, is a $k$-point of $\Kst$. 
Let $D_X$ be a hull of $X$ and $T_{D_X}$ its tangent space.

\part [i] The natural map $res: t_{D_X}\to H^0(X,T^1_X)$ given by restricting
deformations of $X$ to deformations of its singularities is an inclusion.

\part [ii] $D_X$ and $\Kst$ are each locally a hypersurface.

\begin{proof}
The computation
of the hypercohomology of the tangent complex gives a short exact sequence
\begin{displaymath}
H^1(X,T_X)\to t_{D_X}\to H^0(X,T^1_X),
\end{displaymath}
where $t_{D_X}$ is the tangent space to a hull $D_X$. 
Now $T_X$ is free, so that $H^1(X,T_X)=0$ and $res$ is an injection.

A Zariski RDP of rank $r$ has a smooth hull of dimension $2r$,
and the subspace of Zariski deformations is smooth of dimension $r$. So
the subspace $D^{Zar}_X$ of $D_X$ corresponding to Zariski deformations has embedding
dimension at most $12$. On the other hand, we know that $X$ has $9$ moduli
as a supersingular K3 surface, so that by \ref{fibre-dimension}
$\dim D^{Zar}_X \ge 11$.
\end{proof}
\end{proposition}

Denote by $\sQ$ the Deligne--Mumford stack $[\sM_N^0/\frak S_{12}]$.
There is an extension of stacks $\Kst\to\Kst/\sG$ by $\sG$.

\begin{proposition}\label{7.2}
The morphism $b:\sM_N^0\to \Kst$ that is
the family $q\co X\to\sM_N^0$ is proper and factors through $\sQ^{(1)}$.
The induced morphism of geometric quotients induces a bijection
on geometric points.
\begin{proof} The properness of $b$ is immediate from its construction
in \ref{Torelli}. 
The group $\frak S_{12}$ acts on the lattice $N$, and so on the
functor represented by $S\to B_N$. So it acts equivariantly on
$S\to B_N$. Then it acts on the restriction $S^0\to B_N^0$ so as
to preserve the exceptional locus of $g$,
so it acts on the contracted family $q\co X\to\sM_N^0$.
Regard this action as a groupoid 
$\sM_N^0\times \frak S_{12}\rightrightarrows \sM_N^0\times\sM_N^0$; this
groupoid maps naturally to the $\Isom$-scheme
$R:=\sM_N^0\times_{q,\Kst,q}\sM_N^0\rightrightarrows \sM_N^0\times\sM_N^0$,
so there are morphisms $\sM_N^0\to \sQ\to [\sM_N^0/R]$. But 
$q:\sM_N^0\to \Kst$ factors naturally through $[\sM_N^0/R]$;
this is a general statement about morphisms from schemes to algebraic stacks.
So $b:\sM_N^0\to \Kst$ factors through $\sQ$.
To show that $b$ factors through $\sQ^{(1)}$ it is enough to show
that it factors through $\sM_N^{0(1)}$.
For this, it is enough
to show that the derivative of $b$ vanishes at any point corresponding to a
surface $X=\sX_s$ with 12 $A_1$ points. 

The inclusion
$t_{D_X}\inj H^0(X,T^1_X)$ described in \ref{local to global for K3s}
shows that it is enough to show that the deformation
of each singular point induced by $q$ is trivial,
locally on $\sX$.
Consider the blowing--down
$g\co S^0\to \sX$; locally over $\sX$, the dualizing sheaf
$\omega_{S^0}$ is trivial, so that $\sX$ is singular 
along the critical locus $\Sing q$ of $q$.
In suitable local co--ordinates,
$(X,P)$ is given by an equation $xy+z^2=0$. All geometric fibres of
the deformation of this induced by $q$ are singular, so that the
deformation is Zariski, so given by
the equation $xy+z^2+F(t_1,\ldots,t_9)=0$,
where the $t_i$ are local co--ordinates on $\sM_N^0$. 
Since $\Sing \sX=\Sing q$, by construction of $g$,
$F$ is a square, which proves
the local triviality and the first part of the proposition.

The last part is an immediate consequence of Ogus' global Torelli
theorem.
\end{proof}
\end{proposition}

Let $c:\sQ^{(1)}\to\Kst$ be the morphism given by \ref{7.2}.

\begin{proposition} $\Kst$ is normal. 
\begin{proof} Consider the morphisms $\phi:\sU\to\Est_{K3,uni}$ 
and $\pi:\sU\to\Kst$. The stack
$\Est_{K3,\Z/2}$ is smooth, since its $k$-points correspond to surfaces
with $H^2(T)=0$, and so $\sU_0$ is smooth. Let $\Kst_0$ be the open
substack of $\Kst$ that is the image of $\sU_0$ (since $\sU\to\Kst$
is flat and representable, this makes sense); then $\Kst_0$ is smooth. 
Since the locus $V$ of points in $\sM_N^0$
corresponding to $K3$ surfaces with even singularities is of codimension
$2$ everywhere, by \ref{3 x D4 locus},
it follows that the complement of $\Kst_0$ in $\Kst$ has
codimension $2$. Since $\Kst$ is everywhere locally a hypersurface,
we are done.
\end{proof}
\end{proposition}

\begin{theorem}\label{composite}
The composite morphism $f:\sQ^{(1)}\to\Kst/\sG$ induced
by \ref{7.2} is a $1$-isomorphism and $\Kst/\sG$ is Deligne--Mumford.
\begin{proof} Let $\Kst_{12\times A_1}$ denote the open substack of $\Kst$ defined by
the property that all geometric fibres of its objects have $12$ $A_1$
singularities and $\sG_{12\times A_1}\to\Kst_{12\times A_1}$ the restriction of $\sG$. 
There is a Galois \'etale cover $\tKst_{12\times A_1}\to\Kst_{12\times A_1}$
with group $\frak S_{12}$ determined by ordering the singularities
and an open subscheme $\sM_{N,12\times A_1}$ obtained by deleting the closure
of the $E_7+5\times A_1$ locus and the closure of the $D_4+8\times A_1$ locus
that fit into a Cartesian diagram
$$\xymatrix{
{\sM_{N,12\times A_1}} \ar[r] \ar[d] & {\tKst_{12\times A_1}}\ar[d]\\
{\sQ_{12\times A_1}} \ar[r] & {\Kst_{12\times A_1}}
}$$
where the vertical morphisms are stack quotients by $\frak S_{12}$.
Since $\sM_N$ has two geometric components that are exchanged by $\frak S_{12}$,
it follows from Ogus' global Torelli theorem that $\Kst$ and so $\Kst_{12\times A_1}$,
is irreducible.

Now pick a smooth presentation
$v:V\to \Kst_{12\times A_1}$. Put $\tV=V\times_{v,\Kst_{12\times A_1}}\tKst_{12\times A_1}$.
Then $v$ is a family $Y\to V$ of $K3$'s that have, amongst their other properties, 
12 disjoint
$A_1$ singularities. Consider $\tY=Y\times_V\tV\to\tV$. We know, by a previous
argument, that the fibre product $\tY\times_{\tV}\tV^{(-1)}\to \tV^{(-1)}$
has a simultaneous resolution, so that there is a commutative diagram
$$\xymatrix{
{\tV^{(-1)}} \ar[r] \ar[d] & {\tV} \ar[d]\\
{\sM_{N,12\times A_1}} \ar[r] & {\tKst_{12\times A_1}.}
}$$
This shows that there is a factorization 
$\sM_{N,12\times A_1}\to\tKst_{12\times A_1}\to \sM^{(1)}_{N,12\times A_1}$
of the geometric Frobenius on $\sM_{N,12\times A_1}$, so a factorization
$\sQ_{12\times A_1}\to \Kst_{12\times A_1}\to\sQ^{(1)}_{12\times A_1}$. 
Since $\sQ_{12\times A_1}\to \Kst_{12\times A_1}$ factors through 
$\sQ^{(1)}_{12\times A_1}$, it follows that the restriction 
$f_{12\times A_1}:\sQ^{(1)}_{12\times A_1}\to \Kst_{12\times A_1}/\sG_{12\times A_1}$ 
of $f$ is split.

Now we know two things. One is Theorem \ref{height 1}, that for a surface $X/k$
corresponding to a geometric point $\Spec k\to \Kst_{12	\times A_1}$
the connected component of the automorphism groupscheme of $X$ is 
${}_1\Aff$, so is the restriction of $\sG$; this tells us that 
$\Kst_{12\times A_1}/\sG_{12\times A_1}$ is a Deligne--Mumford stack. 
The other is the part
of Ogus' global Torelli theorem, saying that the Picard group
of the minimal resolution $\tX$ of $X$ when $X/k$
corresponds to a geometric point $\Spec k\to \Kst$
rigidifies $\tX$;
since $\Aut_X(k)=\Aut_{\tX}(k)$ 
this tells us that, in the notation of \ref{injectivity},
the morphisms $f_{x,x'}$ arising from $f:\sQ^{(1)}\to\Kst/\sG$ are embeddings. 
Then $f_{12\times A_1}$ is a split morphism of irreducible normal algebraic
stacks that is representable, by \ref{injectivity}; 
it follows that $f_{12\times A_1}$ is an isomorphism. Then $f$ is an isomorphism,
by \ref{isomorphism}.
\end{proof}
\end{theorem}

\begin{corollary}\label{7.9} \part [i] If $X$ is a singular K3 surface corresponding to 
a geometric point of $\sK$, then $\Aut_X^0$ is of height $1$.

\part [ii] $\sG={}_1\sG$.

\begin{proof} This follows from the fact that $\sK/\sG$ is Deligne-Mumford.
\end{proof}
\end{corollary}

Now consider the diagram with $2$-Cartesian squares
$$\xymatrix{
{\sU'} \ar[r]^i\ar[d] & {\sU} \ar[r]^{\phi}\ar[d] & {\Est_{K3,uni}}\\
{\sP'} \ar[r]^{j} \ar[d]^{\pi'} & {\sP}\ar[d]_{\pi}\\ 
{\sQ^{(1)}} \ar[r]^{c} & {\Kst} \ar[r]^{b} & {\Kst/\sG,}
}$$
that extends our previous diagram. (The two upper
vertical arrows are open embeddings.) 
Via \ref{composite}, we identify
$\Kst/\sG$ with $\sQ^{(1)}$ and $bc$ with $1_{\sQ^{(1)}}$. Let
$\pi_0:\sU\to\Kst$ be the composite.

We can recover a result of Liedtke [L10].

\begin{corollary} 
$\sK$, $\sU$ and $\sE_{K3,uni}$ are smooth and irreducible.

\begin{proof} Since $\sQ^{(1)}$ is smooth and irreducible, the smoothness
and irreducibility of $\sK$
follows from \ref{composite} and \ref{gerbes with flat stabilizer are smooth}.
Then $\sU$ is smooth and irreducible; since $\phi$ is smooth, so is $\sE_{K3,uni}$.
\end{proof}
\end{corollary}

Consider the relative groupscheme $\sL=\sG\times_{\sK,c}\sQ^{(1)}\to \sQ^{(1)}$.
This is finite, flat, of order $4$ and of height $1$. Put $\frak L=\Lie(\sL)$;
then $\sP'=\P(\frak L^\vee)$. There is a tautological subgroupscheme $\sM=j^*\sH$
of order $2$ in
$\pi'^*\sL$, where $\sH$ is the tautological subgroupscheme of order $2$ in $\pi^*\sG$

\begin{theorem}\label{0.9}\label{8.11} 
\part[i] The morphisms
$\sQ^{(1)}\stackrel{c}{\to}\sK\stackrel{b}{\to}\sQ^{(1)}$ 
provide a $1$-isomorphism from $\sK$ to the
classifying stack $\sB\sL$ with its canonical projection $\sB\sL\to \sQ^{(1)}$
and its tautological section.

\part[ii] $j:\sP'\to\sP$ and $i:\sU'\to\sU$ induce $1$-isomorphisms 
$[\sP'/\sL]\to\sP$ and $[\sU'/\sL]\to \sU$.

\part[iii] The stacks $\sU$ and $\Est_{K3,uni}$
have Deligne--Mumford quotients $\sU_{DM}$ and
${\Est_{K3,uni}}_{DM}$, and
the morphism $\sU\to\sU'/\sL$ induced from the isomorphism
of \DHrefpart{ii} factors through $\phi:\sU\to\Est_{K3,uni}$ and induces
$1$-isomorphisms  
$\sU_{DM}\to(\sE_{K3,uni})_{DM}$ and
$(\sE_{K3,uni})_{DM}\to \sU'/\sL$.

\part[iv] $\sP'\to\sP'/\sL$ is the Frobenius relative to $\sQ^{(1)}$ 
and $\sP'/\sL\to\sQ^{(1)}$ is a $\P^1$-bundle.

\begin{proof} 
Except for \DHrefpart{iv}, these all follow from 
the general results of \S \ref{stacktology}. 
More precisely, \DHrefpart{i} follows from \ref{composite}
and \ref{neutral gerbes are classifying stacks}. \DHrefpart{ii} is a 
consequence of \ref{1.8}. For \DHrefpart{iii}, the first statement follows
from \ref{summary} and \ref{summary sequel} and the second from \ref{1.16}.

For \DHrefpart{iv} we make a calculation.

We know that outside the $3\times D_4$ locus in $\sQ^{(1)}$,
which is of pure codimension $2$, the geometric fibres of $\sL$
are ${}_1\Aff$, and over the $3\times D_4$ locus they are
$\alpha_2\times\alpha_2$. Hence, locally on $\sQ^{(1)}$ around
a point in the $3\times D_4$ locus,
there is a basis $\{e_1,e_2\}$ of $\Lie(\sL)$ and a regular
sequence $s_1,s_2$ in $\sO=\sO_{\sQ^{(1)}}$ such that
$[e_1,e_2]=s_1e_1+s_2e_2$. Let $\{e_1^*,e_2^*\}$ be
the dual basis of $\Lie(\sL)^*$ and 
$R=\sO[e_1^*,e_2^*]=\Symm^\bullet\Lie(\sL)^*$. We need to
compute the $\sO$-subalgebra of $\sL$-invariants of $R$.
Since $\sL$ is of height $1$, this is the algebra of $\Lie(\sL)$-invariants,
and is easily computed to be $\sO[(e_1^*)^2,(e_2^*)^2]$.
This completes the proof.
\end{proof}
\end{theorem}

\begin{corollary}\label{8.12} There is a period morphism 
$\psi:(\Est_{K3,uni})_{DM}\to\sQ^{(1)}$
that is an open piece of a $\P^1$-bundle.
\noproof
\end{corollary}

\begin{theorem}\label{8.13} There is a closed substack $\sE_{K3,\alpha}$ of
$\sE_{K3,uni}$ whose geometric points are the geometric points
of $\sE_{K3,uni}$ whose $\Pic^\tau$ is isomorphic to $\alpha_2$.

\begin{proof} Define $\sE_{K3,\alpha}$ to be the closed substack of
$\sE_{K3,uni}$ that corresponds under the bijection given by
\ref{closed substacks of smooth gerbes} to $\sU_{nilpt}$. For any $h:S\to\sU$,
the restriction $\sH_S\to S$ of $\sH\to\sU$ is a twisted form of $\alpha_2$
if and only if $h$ factors through $\sU_{nilpt}$. Identify $h$ with
$h:(X,\Xi)\to S$ and put $Y=X/\exp(\Xi)\to S$. Then $X\to Y$ is a torsor under
$\sH_S\to S$, so if $S$ is a geometric point, then $S\to \sE_{K3,uni}$
defined by $Y$ factors through $\sE_{K3,\alpha}$ if and only if
$\Pic^\tau(Y)\cong\alpha_2$.
\end{proof}
\end{theorem}

\begin{remark} It is not so clear how to give an \emph{a priori}
definition of $\sE_{K3,\alpha}$ even as a stack, let alone a substack
of $\sE_{K3,uni}$; the difficulty stems from the fact that $\alpha_2$
has a non-trivial automorphism groupscheme, namely $\mul$.
\end{remark}

To describe $\sE_{K3,\alpha}$ further needs some more notation.
We have the $3\times D_4$-locus $Z\subset \sM_N^0$, which is smooth and of
codimension $2$. It maps to the closed
substack $\sZ^{(1)}=[Z^{(1)}/\frak S_{12}]$ of $\sQ^{(1)}$. 
Then the inverse image $\sY$ of $\sZ^{(1)}$ in $\Kst$, via the isomorphism 
$f:\sQ^{(1)}\to \Kst/\sG$, is smooth and $\sY\to\sZ^{(1)}$ is an extension
by the restriction $\sG_{\sY}$ of $\sG$ to $\sY$. Note that the complement
of $\sY$ in $\Kst$ is the open substack $\Kst^{odd}$.
Let $\sW$ be the inverse image
of $\sY$ in $\sP$ and $\sV=\sW\cap \sU$.

\begin{proposition}\label{8.14} The image $(\Est_{K3,\alpha})_{DM}$ of
$\Est_{K3,\alpha}$ in $(\Est_{K3,uni})_{DM}$ is a section of $\psi$ over
the complement of the $3\times D_4$-locus. That is, the period morphism
is an isomorphism on the Deligne--Mumford quotient of the stack of
multiplicative $\alpha_2$-Enriques surfaces.
\noproof
\end{proposition}

\begin{theorem} $\sE_{K3,\alpha}$ is irreducible and of codimension $1$
in $\sE_{K3,uni}$.

\begin{proof} We know by \ref{alphadef} that $\sE_{K3,\alpha}$ is everywhere
of codimension $1$ in $\sE_{K3,uni}$. So
it is enough to show that $\sP_{nilpt}$ has a unique component of
codimension $1$. Now $\sP_{nilpt}$ is the union of two pieces, namely
$\pi^{-1}(\sY)$ and a section of $\pi$ over $\Kst^{odd}$. The first is
of codimension $2$ and the second is irreducible, and we are done.
\end{proof}
\end{theorem}

\begin{theorem}\label{8.16} There is an irreducible smooth open substack
$\sE_{K3,\alpha,mult}$ of $\sE_{K3,\alpha}$ whose geometric points
are multiplicative $\alpha_2$-Enriques surfaces. Its complement
is a smooth irreducible closed substack $\sE_{K3,\alpha,add}$
whose geometric points are additive $\alpha_2$-Enriques
surfaces.

\begin{proof} Define $\sE_{K3,\alpha,add}$ to be the closed substack
of $\sE_{K3,uni}$ that corresponds, according to 
\ref{closed substacks of smooth gerbes}, to $\sU\cap \pi^{-1}(\sY)$.
Then take $\sE_{K3,\alpha,mult}$ to be its complement in
$\sE_{K3,\alpha}$.
\end{proof}
\end{theorem}
\end{section}

\begin{section}{Periods of $12A_1$-Enriques surfaces}\label{12A1}

We will in this section make a special study of $12A_1$-Enriques surfaces and
their periods. We begin by showing that all multiplicative $\alpha_2$-surfaces
belong to this category.
We will now make a closer study of the subset $\P^0$ of the $\P^1$-bundle $\P
\to \sM^0_N$ or rather its complement in $\P$.
\begin{lemma}
Let \map\pi{\sX}S be a flat and proper family of RDP-K3 surfaces whose
singularities are Zariski RDP's of total index $12$ and such that $T_{\sX/S}$ is
fibrewise free so that $\pi_*T_{\sX/S}$ is locally free of rank $2$ and commutes
with base change. Let $Z$ be the singular locus of $\pi$ (i.e., locally defined
by the partial derivatives of a defining equation for $\sX$ over $S$). Then $Z$
is finite flat of degree $12$ over $S$ and there is a natural map $Z \to
\P(\pi_*T_{\sX/S})$.
\begin{proof}
The fact that $\pi_*T_{\sX/S}$ is locally free and commutes with base change as
well as the fact that $Z$ is finite flat of degree $12$ is clear. To continue we
note that we have a natural evaluation map $\Omega^1_{\sX/S} \to
\Hom_X(T_{\sX/S},\sO_X)$ which is an injection and whose cokernel $\sL_Z$ is a line
bundle over $Z$. This can be verified locally and fibrewise and is then a direct
computation. This gives us a map $\pi_*T_{\sX/S} \to \sL_Z$ and we claim that
the induced  map $\pi_*T_{\sX/S}\Tensor_{\sO_S}\sO_Z \to \sL_Z$ is surjective. This
however follows directly from the fact that $\pi^*\pi_*T_{\sX/S} \to T_{\sX/S}$
is surjective. The map $\pi_*T_{\sX/S}\Tensor\sO_Z \to \sL_Z$ now induces the
desired map.
\end{proof}
\end{lemma}
With the aid of this lemma and provided with a map \map\pi{\sX}S as in it, we
may now use 
[Mu65], \S 5.3 to define a Cartier divisor $\Div(Z)$, the
\Definition{non-free locus},
in $\P(\pi_*T_{\sX/S})$. This divisor commutes with base change by [Mu65], \S 5.3
and is hence a relative Cartier divisor as it gives a divisor on
each fibre. By construction it is of degree $12$ over the base. It is
not \'etale over a point of $S$ over which $Z$ is not \'etale and $Z\to S$ is
\'etale
precisely over the points at which $\sX$ has $12$ $A_1$-singularities. Hence,
when we start with a family of K3-Enriques surfaces and let the RDP-K3-surface
be the canonical double cover, the non-free locus has a chance of being \'etale
only over points over which the Enriques surface an $12A_1$-surface. Our
aim is not to show that in that case non-free locus is indeed \'etale or
equivalently, $Z \to \P(\pi_*(T_{\sX/S}))$ is an embedding. We will do this by
reading off the position of the points of $Z$ in $\P(\pi_*(T_{\sX/S}))$ from the
period of the surface. To do this we start by characterising those periods that
correspond to $12A_1$-surfaces.
\begin{lemma}\label{M-superlattice}
Any even superlattice of $M$ is generated by $M$ and the $-2$-elements it
contains.
\begin{proof}
It is enough to assume that such a superlattice $L$ contains $M$ as a subgroup of
index $2$ and is hence generated by $M$ and one other element. That extra
element may be assumed to have the form $\frac12(a_1,\dots,a_{12})$ where $a_i
\in \{0,1\}$ and the fact that $L$ should be even forces the number $a_i$'s that 
are equal to $1$ to be divisible by $4$. As $\frac12(1,1,\dots,1)$ belongs to
$M$ we see that $L$ must contain an element of the form
$\frac12(a_1,\dots,a_{12})$ with exactly $4$ $a_i$'s being equal to $1$ but that
is a $-2$-element.
\end{proof}
\end{lemma}
We need to relate $N$-periods to the orthogonal decomposition $N=E(2)\perp M$.
\begin{lemma}
Let $V \subset N\Tensor k$ be a period in $\sM^0_N$ where $k$ is a perfect
field of characteristic $2$.

\part[i] The projection of $V$ into $M\Tensor k$ has as image an $\F_2$-rational
subspace.

\part[ii] $V$ corresponds to an RDP-K3 surface with $12$ $A_1$-singularities if
and only if the projection of $V$ into $M\Tensor k$ is surjective.
\begin{proof}
If $\varphi$ is the Frobenius map on $N\Tensor k$ with respect to the
$\F_2$-rational structure given by $N$, then $V+\varphi(V)$ is $11$-dimensional
and so has a non-trivial intersection with $E(2)\Tensor k$. We claim that $V
\cap E(2)\Tensor k \ne (V+\varphi(V))\cap E(2)\Tensor k$. In fact if they are
equal, then $V \cap E(2)\Tensor k$ is stable under $\varphi$ and hence is
$\F_2$-rational. It is totally isotropic because it is a subspace of $V$. Hence, it
corresponds to an even superlattice of $E(2)$ but this contradicts the
definition of $\sM^0_N$ and lemma \ref{E(2)-superlattices}. Hence,
$V+\varphi(V)$ is the sum of $V$ and $(V+\varphi(V))\cap E(2)\Tensor k$. This
means that the projection onto $M\Tensor k$ is stable under $\varphi$ and is
hence $\F_2$-rational. This proves \DHrefpart{i}.

As for \DHrefpart{ii}, as $V$ is maximal isotropic, $V$ will contain the
orthogonal complement in $M\Tensor k$ to the image of $V$ in $M\Tensor k$. As
the image is $\F_2$-rational so is its orthogonal complement. That orthogonal
complement thus equals $V\cap M$ so it is also $\F_2$-rational and it is, by
\ref{M-superlattice} trivial if and only if $V$ corresponds to an RDP-K3 surface
with $12$ $A_1$-singularities.
\end{proof}
\end{lemma}
\begin{definition}
If $1\le i \ne j \le 12$ we set $e_{ij}:= e_i+e_j \in M/2M$. We let
\Definition{$\sM^1_N$} be the open subvariety of $\sM^0_N$ consisting of those
periods $V \in N\Tensor k$ for which $V \cap M/2M = \{0\}$ and which contains no 
element of the form $\alpha+e_{ij}$ with $\alpha \in E(2)/2E(2)$.
\end{definition}
\begin{remark}
The complement of $\sM^1_N$ in $\sM^0_N$ contains points which correspond to
RDP-K3 surfaces with $12$ $A_1$-singularities. Indeed, a period in $\sM^0_N$
corresponds to such a surface precisely when $V \cap M/2M = \{0\}$ and if we
pick $\alpha \in E(2)/2E(2)$ with $\alpha^2=1$ then the set of such $V$ that
contain $\alpha+e_{12}$, say, is non-empty and a generic point of it contains no
other points of $N/2N$.
\end{remark}
Our next step is to get a direct description of the global sections of the
tangent bundle of an RDP-K3 surface with only $A_1$-singularities in terms of a
resolution of singularities. 
\begin{lemma}\label{non-freeAndChern}
Let $X$ be a surface over a field $k$ whose singular locus consists of only
$A_1$-singularities, let \map \pi\tX X be a minimal resolution of singularities,
$Z$ the exceptional divisor for $\pi$, and \map jUX the inclusion of the
non-singular locus.

\part[i] We have $T_X=j_*T_U=\pi_*T_{\tX}(-\log Z)(Z)$. 

\part[ii] When the characteristic of $k$ equals $2$, the sheaf of non-free
vector fields equals the kernel of $\pi_*$ applied to the natural map
$T_{\tX}(-\log Z)(Z) \to T_Z(Z)$.

\part[iii] If $\omega_X$ is trivial then using a trivialisation of it
$T_{\tX}(\log Z)(Z)$ is identified with $\Omega^1(\log Z)$, compatibly with the
identification of $T_{\tX}$ and $\Omega^1_{\tX}$, and $T_Z(Z)$ with
$\sO_Z$. Under these identifications the map 
$T_{\tX}(-\log Z)(Z) \to T_Z(Z)$ is
identified with the Poincar\'e residue map $\Omega^1_{\tX}(\log Z) \to \sO_Z$.
\begin{proof}
That $T_X=j_*T_U$ follows directly from $j_*\sO_U=\sO_X$. Furthermore, we have
that $j_*T_U=\pi_*T_{\tX}(*Z)$ and what remains to show for \DHrefpart{i} is
that the inclusion $T_{\tX}(-\log Z)(Z) \hookrightarrow T_{\tX}(*Z)$ induces an
isomorphism upon applying $\pi$. This is local so we may assume that $X$
contains only one singular point and then $Z \iso \P^1$. Now, the short exact
sequence
\begin{displaymath}
\shex{T_{\tX}(-\log Z)}{T_{\tX}}{\sN_{Z/\tX}}
\end{displaymath}
and the fact that $\sN_{Z/\tX} \iso \sO_Z(-2)$ shows that $\pi_*T_{\tX}(-\log
Z)(Z)=\pi_*T_{\tX}(Z)$ and the facts that $\sO_Z(Z)=\sO_Z(-2)$ and 
$T_{\tX}(-\log Z)/T_{\tX}(-Z)=T_Z\iso \sO_Z(2)$ shows that
$\pi_*T_{\tX}(iZ)=\pi_*T_{\tX}((i+1)Z)$ for all $i >1$.

Continuing with \DHrefpart{ii}, the kernel of 
$T_{\tX}(-\log Z)(Z) \to T_Z(Z)$ is $T_{\tX}$ so that the kernel of 
the map induced by applying 
$\pi_*$ to it consists of those vector fields that lift to the resolution and it 
is a simple local calculation that these are exactly the non-free vector fields.

Finally, to prove \DHrefpart{iii} as the statement is independent of the choice
of trivialisation of $\omega_X$ we may work locally and choose coordinates such
that $Z$ is defined by the vanishing of one of them and the trivialisation is
given by $dx \wedge dy$. Then the result is a simple calculation.
\end{proof}
\end{lemma}
We are now prepared to state and prove the main result of this section.
\begin{theorem}
Suppose that a period point $x$ in $\sM^0_N$ corresponds to the RDP-K3 surface $X$. Then the
non-free locus in $\P(H^0(X,T_X))$ is \'etale precisely when $x$ lies in $\sM^1_N$.
\begin{proof}
As we have seen for the non-free locus to be \'etale it is necessary that $X$
has $12$ $A_1$-singularities and if the period lies in $\sM^1_N$ this is always
the case. We may and will therefore assume that $X$ has $12$
$A_1$-singularities.

If $W$ is the singular locus of $X$ then, as has already been noted, the
non-free locus is \'etale precisely when the map $W \to \P(H^0(X,T_X))$ is
injective. For each $x \in W$ let $T_X \to \ell_x$ be the quotient by the non-free
subsheaf with respect to $x$. Then this means that for any two $x \ne y \in W$,
the map $H^0(X,T_X) \to \ell_x\Dsum\ell_y$ is surjective. Let $\tX \to X$ be a minimal
resolution of singularities and $Z$ its exceptional divisor. We have a short
exact sequence
\begin{displaymath}
\shex{\Omega^1_{\tX}}{\Omega^1_{\tX}(\log Z)}{\sO_Z},
\end{displaymath}
and using the fact that $h^0(\tX,\Omega^1_{\tX})=0$ we see that
$H^0(\tX,\Omega^1_{\tX}(\log Z))$ can be identified with the kernel of the map
$H^0(Z,\sO_Z) \to H^1(\tX,\Omega^1_{\tX})$. Now, $Z$ is the disjoint union of
its irreducible components $Z_x$ for $x \in W$ and thus 
$H^0(Z,\sO_Z)=\Dsum_{x \in W}k1_{Z_x}$. 
Furthermore, the image in $H^1(\tX,\Omega^1_{\tX})$ of $1_{Z_x}$
equals the de Rham cycle class of $Z_x$. Hence, we may identify
$H^0(\tX,\Omega^1_{\tX}(\log Z))$ with the kernel of the map 
$\sum^{12}_{i=1}k e_i \to H^1(\tX,\Omega^1_{\tX})$ 
which sends $e_i$ to the first Chern class of
$\sO(Z_x)$. Using \ref{non-freeAndChern} we can identify  $H^0(X,T_X)$ and
$H^0(\tX,\Omega^1_{\tX}(\log Z))$ and then the map $H^0(X,T_X) \to \ell_x$ can be
identified with the projection of $H^0(\tX,\Omega^1_{\tX}(\log Z))$ onto the 
$k e_i$-factor where $x \in W$ is the $i$'th point (for some fixed
ordering). Furthermore, if $V$ is the period of $\tX$ then the kernel of the map 
$N\Tensor k \to H^1(\tX,\Omega^1_{\tX})$ is $V+\varphi(V)$, where $\varphi$ is
the Frobenius map with respect to $N/2N$ (cf. [Og79], 3.20.1). Hence, if
we consider $V$ as a subspace of $E(2)\Tensor k\perp M'\Tensor k$, where
$M'=\sum_i\Z e_i$ then $H^0(X,T_X)$ can be identified with 
$(V+\varphi(V))\cap M'\Tensor k$.

We have thus reduced the theorem to proving that $V \in \sM^1_N$ precisely when
$(V+\varphi(V))\cap M'\Tensor k$ maps onto $k e_i\Dsum k e_j$ for all 
$i \ne j$. As $\frac12\sum^{12}_{i=1}e_i \in M$ we see that the diagonal 
$\lambda e_{ij}$ will always lie in this image and hence the map fails to be surjective
exactly when the image is this diagonal, i.e., when 
$(V+\varphi(V))\cap M\Tensor k$ is contained in $M_{ij}\Tensor k$, 
where $M_{ij}:=\set{\sum_ra_re_r\in M}{a_i=a_j}$. 
If $V$ contains an element of the form $\alpha+e_{ij}$,
$\alpha \in E(2)/2E(2)$, then $V$ and therefore $\varphi(V)$ as well as
$V+\varphi(V)$ is orthogonal to $\alpha+e_{ij}$ and thus any element in
$(V+\varphi(V))\cap M\Tensor k$ lies in $M_{ij}\Tensor k$. Conversely, let us
assume that $(V+\varphi(V))\cap M\Tensor k$ lies in $M_{ij}\Tensor k$. Then
$V+M_{ij}$ is stable under $\varphi$ and is hence $\F_2$-rational. As $V$ is the
graph of an isomorphism between $E(2)\Tensor k$ and $M\Tensor k$ the sum
$V+M_{ij}$ is direct and its dimension is therefore equal to $19$. Its
annihilator is a $1$-dimensional subspace of $N/2N$ and is contained in $V$ as
$V$ is maximal totally isotropic. Furthermore, it does not lie in $E(2)/2E(2)$
as $V \in \sM^0_N$ and as it is orthogonal to $M_{ij}$ a non-zero element in it
must therefore be of the form $\alpha+e_{ij}$.
\end{proof}
\end{theorem}
\end{section}
\bibliography{alggeom,ekedahl}
\bibliographystyle{pretex}
\end{document}